\documentclass[10pt]{article}
\usepackage{times}

\DeclareMathVersion{varnormal}
\DeclareMathVersion{varbold}

\SetSymbolFont{letters}{normal}{OML}{txmi}{m}{it}
\SetSymbolFont{letters}{bold}{OML}{txmi}{bx}{it}
\SetSymbolFont{letters}{varnormal}{OML}{mdugm}{m}{it}
\SetSymbolFont{letters}{varbold}{OML}{mdugm}{b}{it}
\SetSymbolFont{operators}{normal}{OT1}{txr}{m}{n}
\SetSymbolFont{operators}{bold}{OT1}{txr}{bx}{n}
\SetSymbolFont{operators}{varnormal}{OT1}{mdugm}{m}{n}
\SetSymbolFont{operators}{varbold}{OT1}{mdugm}{b}{n}
\SetSymbolFont{symbols}{normal}{OMS}{txsy}{m}{n}
\SetSymbolFont{symbols}{bold}{OMS}{txsy}{bx}{n}
\SetSymbolFont{symbols}{varnormal}{OMS}{mdugm}{m}{n}
\SetSymbolFont{symbols}{varbold}{OMS}{mdugm}{b}{n}
\SetSymbolFont{largesymbols}{normal}{OMX}{txex}{m}{n}
\SetSymbolFont{largesymbols}{bold}{OMX}{txex}{bx}{n}
\SetSymbolFont{largesymbols}{varnormal}{OMX}{mdugm}{m}{n}
\SetSymbolFont{largesymbols}{varbold}{OMX}{mdugm}{b}{n}

\SetMathAlphabet{\mathrm}{normal}{OT1}{txr}{m}{n}
\SetMathAlphabet{\mathrm}{bold}{OT1}{txr}{bx}{n}
\SetMathAlphabet{\mathrm}{varnormal}{OT1}{mdugm}{m}{n}
\SetMathAlphabet{\mathrm}{varbold}{OT1}{mdugm}{b}{n}

\SetMathAlphabet{\mathit}{normal}{OT1}{txr}{m}{it}
\SetMathAlphabet{\mathit}{bold}{OT1}{txr}{bx}{it}
\SetMathAlphabet{\mathit}{varnormal}{OT1}{mdugm}{m}{it}
\SetMathAlphabet{\mathit}{varbold}{OT1}{mdugm}{b}{it}

\usepackage{amsmath}
\usepackage{amsthm}\theoremstyle{plain}
\usepackage{mathrsfs}
\usepackage{amssymb}
\usepackage{empheq}
\usepackage[labelformat=simple]{subcaption}  
\usepackage[paper=a4paper,top=36pt,left=100pt,right=100pt,bottom=50pt,foot=25pt]{geometry}
\usepackage{titlesec}
\usepackage{lineno}
\usepackage{hyperref}
\usepackage{nicefrac}
\usepackage{float}
\usepackage{upgreek}
\usepackage[dvipsnames]{xcolor}
\usepackage[numbers,compress]{natbib}

\titleformat{\section}{\normalfont\bfseries}{\thesection}{1em}{}
\titleformat{\subsection}{\normalfont}{\thesubsection}{1em}{}

\renewcommand\thesubfigure{(\alph{subfigure})}

\renewcommand\appendix{\par
  \setcounter{section}{0}
  \setcounter{subsection}{0}
  \setcounter{figure}{0}
  \setcounter{table}{0}
  \setcounter{equation}{0}
  \renewcommand\thesection{Appendix \Alph{section}}
  \renewcommand\thefigure{\Alph{section}\arabic{figure}}
  \renewcommand\thesubfigure{(\alph{subfigure})}
  \renewcommand\thetable{\Alph{section}\arabic{table}}
  \renewcommand\theequation{\Alph{section}\arabic{equation}}  
}

\newcommand\supp{\par
  \setcounter{section}{0}
  \setcounter{subsection}{0}
  \setcounter{figure}{0}
  \setcounter{table}{0}
  \setcounter{equation}{0}  
  \renewcommand\thesection{S\arabic{section}}
  \renewcommand\thefigure{S-\arabic{figure}}
  \renewcommand\thetable{S-\arabic{table}}
  \renewcommand\theequation{S-\arabic{equation}}  
}

\usepackage{xr}

\usepackage[normalem]{ulem}

\title{\Large Nonuniform 3D finite difference elastic wave simulation on staggered grids}
\author{Longfei Gao\thanks{Oden Institute for Computational Engineering and Sciences, The University of Texas at Austin, Austin, TX 78712, USA, Email address: longfei.gao@austin.utexas.edu}\quad
Omar Ghattas\thanks{Oden Institute for Computational Engineering and Sciences, Jackson School of Geosciences, and Department of Mechanical Engineering, The University of Texas at Austin, Austin, TX 78712, USA}\quad
David Keyes\thanks{Applied Mathematics and Computational Science Program and Extreme Computing Research Center, King Abdullah University of Science and Technology, Thuwal 23955-6900, Saudi Arabia}
}

\date{}

\begin{document}


\maketitle

\begin{abstract}
We present an approach to simulate the 3D isotropic elastic wave propagation using nonuniform finite difference discretization on staggered grids.
Specifically, we consider simulation domains composed of layers of uniform grids with different grid spacings, separated by nonconforming interfaces.
We demonstrate that this layer-wise finite difference discretization has the potential to significantly reduce the simulation cost, compared to its fully uniform counterpart.
Stability of such a discretization is achieved by using specially designed difference operators, which are variants of the standard difference operators with adaptations near boundaries or interfaces, and penalty terms, which are appended to the discretized wave system to weakly impose boundary or interface conditions.
Combined with specially designed interpolation operators, the discretized wave system is shown to preserve the energy conserving property of the continuous elastic wave equation, and {\it a fortiori} ensure the stability of the simulation. 
Numerical examples are presented to demonstrate the efficacy of the proposed simulation approach.
\end{abstract}

\section{Introduction}

Simulating waves propagating through earth media is a routine task in seismic exploration.
Elastic and acoustic assumptions of earth media are common choices in practice.
The elastic model offers a more accurate depiction since, 
in addition to the compressional wave supported by the simpler acoustic model,
it also supports shear and surface waves, which are commonly observed in seismic survey data and can contain significant energy and information.
For applications involving salt bodies where significant conversions between compressional and shear waves occur at the salt boundaries (see, e.g., \cite{ogilvie1996effects}) or applications requiring an accurate account of the near surface region where surface wave has a strong presence (see, e.g., \cite{miller1955partition}), the importance of elastic model can be particularly pronounced. 
However, elastic wave simulation is significantly more expensive than its acoustic counterpart because of the increased number of variables and derivatives in the equations as well as the need to resolve the lower shear wave speed.

For a fixed frequency, wavelengths of the seismic waves are proportional to the wave speeds of the media.
The spatial grid spacing in a wave simulation is usually dictated by the minimum wavelength of the propagating waves, and hence is determined by the minimum wave speed of the media.
Based on the observation that wave speeds of earth media generally increase with depth due to geological sedimentation and consolidation, 
it is desirable to have progressively coarsened spatial discretization from earth surface to deeper depth in order to alleviate the computational burden.
To illustrate, soil and dry sand on earth's surface can have shear wave speeds as low as 100~m/s while hard rocks 
deep beneath can have shear wave speeds above 3000~m/s and compressional wave speeds around 6000~m/s; see \cite[p.~240]{thierry1987acoustics}.
Using a uniform grid with such a large contrast in wave speeds would lead to wasteful oversampling in the discretization (see \ref{appendix_thought_experiment} for a detailed explanation).

Fully unstructured finite element based discretizations can address the variation in wave speed by allowing variation in the sizes of the individual elements. 
Such an approach has been thoroughly investigated in previous works in the context of seismic wave propagation; see, e.g., \cite{lysmer1972finite, BaoBielakGhattasEtAl98, komatitsch1998spectral, komatitsch1999introduction, kaser2006arbitrary, dumbser2006arbitrary, chaljub2007spectral, WilcoxStadlerBursteddeEtAl10, bui2012analysis}.
However, these discretizations tend to be more complex to implement and more expensive to compute on a per grid point basis than those using finite differences.

\begin{figure}[h]
\captionsetup{width=0.9\textwidth, font=small,labelfont=small}
\centering\includegraphics[scale=0.05]{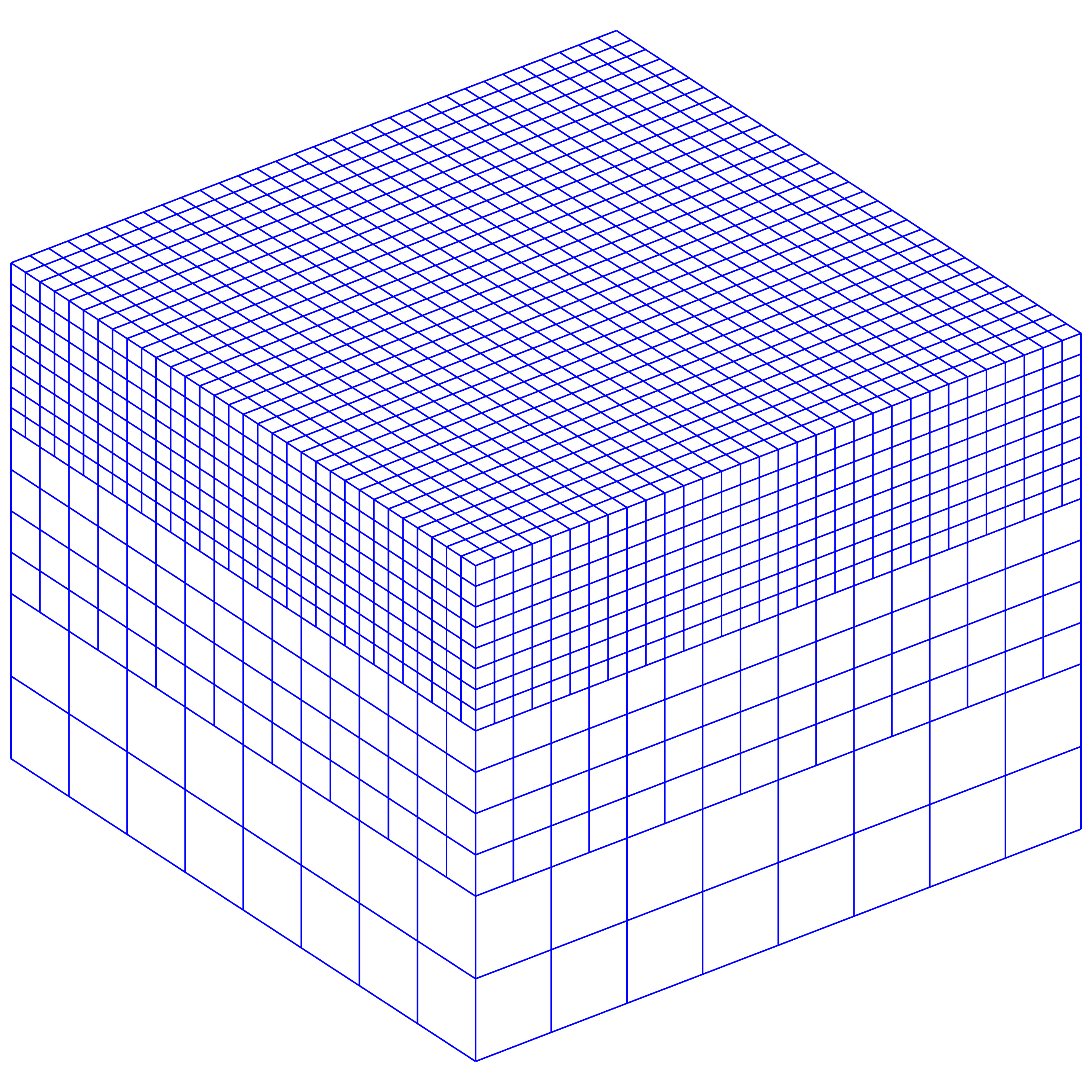}
\caption{ 
Illustration of the layer-wise uniform grid discretization setting: 
Each individual layer is discretized uniformly; 
different layers are discretized using different grid spacings.
}
\label{fig_discretization_setting}
\end{figure}

For earth media, at least at the exploration scale, finite difference discretizations on layer-wise uniform grids as illustrated in Figure~\ref{fig_discretization_setting} seem to be an attractive compromise between the complex fully unstructured finite element schemes and the rigid fully uniform finite difference schemes.
In the context of seismic wave propagation, previous works adopting this approach can be found in, e.g., \cite{jastram1994elastic, aoi19993d, hayashi2001discontinuous, kristek2010stable, zhang2013stable, nie2017fourth, gao2018long}.
However, these earlier attempts are often affected by instability issues associated with the nonconforming interfaces; see, e.g., \cite{gao2018long} for a clear demonstration.

In this work, we adopt this layer-wise finite difference discretization approach and 
make use of the summation by parts (SBP)-simultaneous approximation terms (SATs) technique to address the nonconforming interfaces and ensure stable simulations. 
Mechanisms of the proposed discretization are detailed in Section \ref{section_methodology}.
For general review of the SBP-SAT technique,
we refer the readers to \cite{fernandez2014review} and \cite{svard2014review}.
Development and application of the SBP-SAT technique on wave propagation problem can be found in \cite{appelo2009stable, sjogreen2012fourth, petersson2015wave, wang2016high, o2017energy, gao2019combining, gao2019sbp, gao2020explicit, gao2020simultaneous}, among others.
In particular, \cite{o2017energy, gao2019combining, gao2019sbp, gao2020explicit, gao2020simultaneous} utilize staggered grids.
This work extends these earlier efforts and demonstrates that the SBP-SAT technique can be effectively applied to large-scale 3D elastic wave simulations on staggered grids with nonconforming interfaces.

\section{Methodology}\label{section_methodology}

\subsection{Problem Description} \label{sub_section_problem_description}
Omitting the external source terms, elastic wave propagation can be described using the following velocity-stress formulation: 
%
\begin{subequations}
\label{Elastic_wave_equation_succinct}
\begin{empheq}[left=\empheqlbrace\,]{alignat = 2}
\displaystyle \rho \dot{v}_{i} & \ = \ \displaystyle \sigma_{ij,j} \, ; \label{Elastic_wave_equation_succinct_a}
\\
\displaystyle \dot{\sigma}_{ij}
&\ = \ \displaystyle C_{ijkl} \left( \tfrac{v_{k,l} + v_{l,k}}{2} \right) , \label{Elastic_wave_equation_succinct_b}
\end{empheq}
\end{subequations}
%
where the Einstein summation convention applies to subscript indices.
In equation \eqref{Elastic_wave_equation_succinct}, $v_i$ and $\sigma_{ij}$ are the sought solution variables, standing for components of the velocity vector and the stress tensor, respectively; $\rho$ and $C_{ijkl}$ are given physical parameters, standing for density and components of the stiffness tensor, respectively.
In the case of isotropic media, equation \eqref{Elastic_wave_equation_succinct_b} can be simplified as 
\begin{equation}
\label{Elastic_wave_equation_succinct_isotropic}
\dot{\sigma}_{ij} = \displaystyle \lambda \delta_{ij} v_{k,k} + \mu \left( v_{i,j} + v_{j,i} \right),
\end{equation}
where $\lambda$ and $\mu$ are the Lam\'e parameters, and $\delta_{ij}$ stands for the Kronecker delta.

The above equations can also be expressed in the following form:
\begin{linenomath}
\begin{subequations}
\label{Elastic_wave_equation_verbose}
\begin{empheq}[left=\empheqlbrace\,]{alignat = 2}
\displaystyle \frac{\partial v_x}{\partial t} & \ = \  
\displaystyle \frac{1}{\rho} \left( \frac{\partial \sigma_{xx}}{\partial x} \ + \ \frac{\partial \sigma_{xy}}{\partial y} \ + \ \frac{\partial \sigma_{xz}}{\partial z} \right) \, ;
\\[-0.5mm]
\displaystyle \frac{\partial v_y}{\partial t} & \ = \  
\displaystyle \frac{1}{\rho} \left( \frac{\partial \sigma_{xy}}{\partial x} \ + \ \frac{\partial \sigma_{yy}}{\partial y} \ + \ \frac{\partial \sigma_{yz}}{\partial z} \right) \, ;
\\[-0.5mm]
\displaystyle \frac{\partial v_z}{\partial t} & \ = \  
\displaystyle \frac{1}{\rho} \left( \frac{\partial \sigma_{xz}}{\partial x} \ + \ \frac{\partial \sigma_{yz}}{\partial y} \ + \ \frac{\partial \sigma_{zz}}{\partial z} \right) \, ;
\\[-0.5mm]
\displaystyle \frac{\partial \sigma_{xx}}{\partial t} & \ = \  
\displaystyle (\lambda + 2\mu) \frac{\partial v_x}{\partial x} \ + \ \lambda \frac{\partial v_y}{\partial y} \ + \ \lambda \frac{\partial v_z}{\partial z} \, ;
\\[-0.5mm]
\displaystyle \frac{\partial \sigma_{yy}}{\partial t} & \ = \  
\displaystyle \lambda \frac{\partial v_x}{\partial x} \ + \ (\lambda + 2\mu) \frac{\partial v_y}{\partial y} \ + \ \lambda \frac{\partial v_z}{\partial z} \, ;
\\[-0.5mm]
\displaystyle \frac{\partial \sigma_{zz}}{\partial t} & \ = \  
\displaystyle \lambda \frac{\partial v_x}{\partial x} \ + \ \lambda \frac{\partial v_y}{\partial y} \ + \ (\lambda + 2\mu) \frac{\partial v_z}{\partial z} \, ;
\\[-0.5mm]
\displaystyle \frac{\partial \sigma_{xy}}{\partial t} & \ = \  
\displaystyle \mu \left( \frac{\partial v_y}{\partial x} \ + \ \frac{\partial v_x}{\partial y} \right) \, ;
\\[-0.5mm]
\displaystyle \frac{\partial \sigma_{yz}}{\partial t} & \ = \  
\displaystyle \mu \left( \frac{\partial v_z}{\partial y} \ + \ \frac{\partial v_y}{\partial z} \right) \, ;
\\[-0.5mm]
\displaystyle \frac{\partial \sigma_{xz}}{\partial t} & \ = \  
\displaystyle \mu \left( \frac{\partial v_x}{\partial z} \ + \ \frac{\partial v_z}{\partial x} \right) \, ,
\end{empheq}
\end{subequations}
\end{linenomath}
which is cumbersome for the ensuing discussion, but can be a convenient reference for practical implementation. %

\begin{figure}[h]
\captionsetup{width=0.9\textwidth, font=small,labelfont=small}
\centering\includegraphics[scale=0.075]{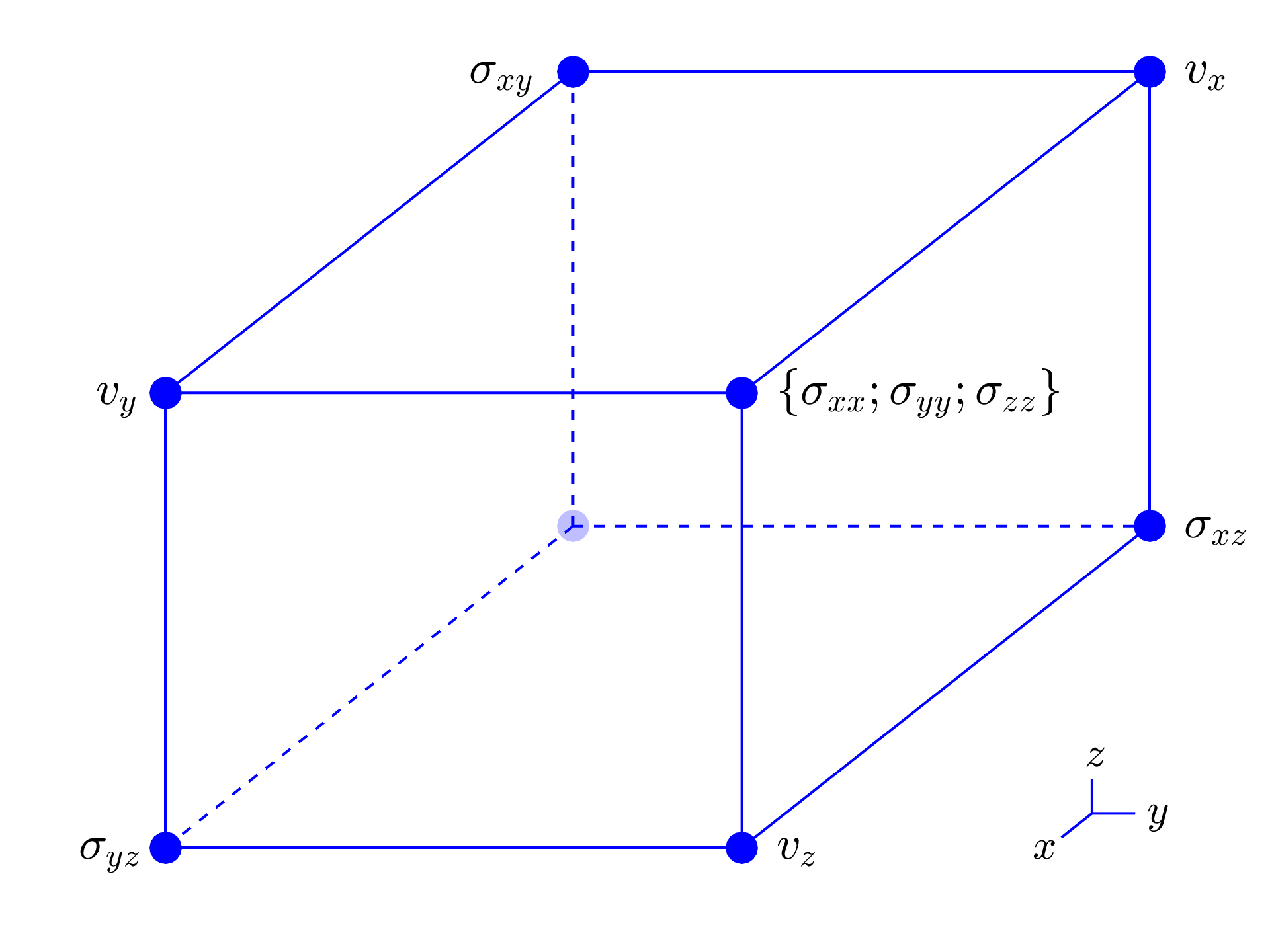}
\caption{Illustration of the grid staggering.}
\label{fig_grid_staggering}
\end{figure}

Within each layer, the spatial discretization is carried out on the staggered grids (see, e.g., \cite{yee1966numerical, virieux1986p, levander1988fourth}), whose layout is sketched in Figure~\ref{fig_grid_staggering}.
On each axial direction, two sets of grid points interlacing each other in a staggering manner, which leads to eight sub-grids in total for a 3D staggered grid discretization. 
Specifically, the three normal stress components ($\sigma_{xx}$, $\sigma_{yy}$, and $\sigma_{zz}$) reside on one sub-grid together; each shear stress component ($\sigma_{xy}$, $\sigma_{xz}$, or $\sigma_{yz}$) resides on its own sub-grid; each velocity component ($v_x$, $v_y$, or $v_z$) also resides on its own sub-grid. These nine solution fields occupy seven sub-grids in total. One sub-grid is unoccupied.

\subsection{Overview of the discretization} \label{subsection_overview}
We use the SBP-SAT technique to address the nonconforming interfaces appearing in the layer-wise discretization setting illustrated in Figure~\ref{fig_discretization_setting}. 
This technique is guided by energy analysis to construct the spatial discretizations.
The concept of SBP finite difference operators dates back to \cite{kreiss1974finite}, which adapt standard finite difference operators near boundaries or interfaces to mimic the integration-by-parts property in the continuous operators of the partial differential equations. 
The concept of SATs dates back to \cite{carpenter1994time}, which append penalty terms (i.e., SATs) to the discretized system to weakly impose boundary or interface conditions.

Applied to the layer-wise discretization setting considered in this work, the individual layers are first discretized independently using SBP finite difference operators and then combined together through SATs.
An important component of the SATs is the compatible interpolation operators that 
transfer the solution at layer interfaces. 
A collection of such interpolation operators corresponding to a variety of grid spacing ratios is included in \ref{appendix_interpolation_operators}.
The grid spacing ratios are allowed to be rational numbers, instead of just integers, which improves the application flexibility of the proposed discretization approach. 
With carefully designed SBP operators, SATs, and interpolation operators, the resulting spatial discretization preserves the energy conserving property of the continuous elastic wave equation, and {\it a fortiori} ensures the stability of the simulation.

We note here that such a discretization procedure bears a similarity to the discontinuous Galerkin (DG) method in the sense that an entire layer can be viewed as the analogue of an element in the DG method, while the SATs resemble the numerical fluxes in the DG method; see \cite[p.~13]{hesthaven2008nodal}. 
However, the `discontinuity' of discretization in the presented method occurs at a much coarser granular level than in the DG method.

\subsection{SBP operators on 1D staggered grids} \label{SBP_operators_1D}
To motivate the concept of SBP operators, let us briefly consider the case of 1D wave propagation on staggered grids.
The 1D wave propagation can be described by the following system: 
\begin{linenomath}
\begin{subequations}
\label{1D_wave_system}
\begin{empheq}[left=\empheqlbrace\,]{alignat = 2}
\displaystyle \rho \frac{\partial v}{\partial t} &\ = \ \displaystyle \frac{\partial \sigma}{\partial x} \, ; 
\label{1D_wave_system_v} \\[0.5ex]
\displaystyle \beta \frac{\partial \sigma}{\partial t} & \ = \ \displaystyle \frac{\partial v}{\partial x} \, , 
\label{1D_wave_system_p}
\end{empheq}
\end{subequations}
\end{linenomath}
where $\beta = \frac{1}{\rho c^2}$ stands for compressibility.

Supposing the above wave system is defined over interval $\left(x_L,x_R\right)$, its associated physical energy can be expressed as:
\begin{equation}
\label{1D_continuous_energy}
\mathscr E \ = \ \frac{1}{2} \int_{x_L}^{x_R} \left( \rho v^2 + \beta \sigma^2 \right) dx \, . 
\end{equation}
Taking the time derivative of equation \eqref{1D_continuous_energy}, substituting the equations \eqref{1D_wave_system_v} and \eqref{1D_wave_system_p}, and applying the integration-by-parts formula, we arrive at:
\begin{equation}
\label{1D_continuous_energy_time_derivative}
\frac{d \mathscr E}{d t} \ = \ - \ \sigma(x_L) \cdot v(x_L) \ + \ \sigma(x_R) \cdot v(x_R) \, .
\end{equation}
In other words, time derivative of the physical energy $\mathscr E$ reduces to boundary data only.
If free surface boundary condition is associated with both boundaries, i.e., $\sigma(x_L) = \sigma(x_R) = 0$, we have that $\tfrac{d \mathscr E}{d t} = 0$, i.e., the physical energy is conservative.

\begin{figure}[H]
\captionsetup{width=0.9\textwidth, font=small,labelfont=small}
\centering\includegraphics[scale=0.075]{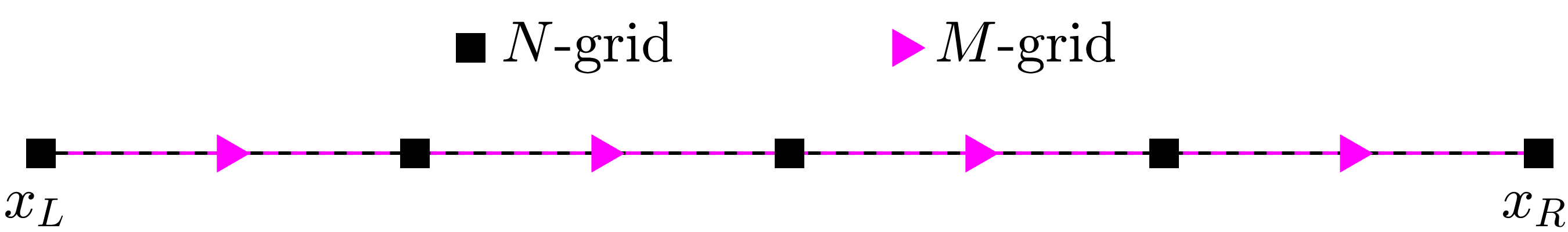}
\caption{Illustration of a set of 1D staggered grids.}
\label{fig_grid_staggering_1D}
\end{figure}

Supposing that the wave system in equation \eqref{1D_wave_system} is discretized on a set of staggered grids illustrated in Figure~\ref{fig_grid_staggering_1D} with $\sigma$ occupying the $N$-grid, which includes the two boundary points, and $v$ occupying the $M$-grid, which consists of only interior grid points, the semi-discretized system can be expressed as:
\begin{linenomath}
\begin{subequations}
\label{1D_semi_discretized_system}
\begin{empheq}[left=\empheqlbrace\,]{alignat = 2}
\displaystyle \mathcal A^M \boldsymbol \rho^M \frac{d V}{d t} & \ = \ \displaystyle \mathcal A^M \mathcal D^N \Sigma \, ; 
\label{1D_semi_discretized_system_a} \\[0.5ex]
\displaystyle \mathcal A^N \boldsymbol \beta^N \frac{d \Sigma}{d t} & \ = \ \displaystyle \mathcal A^N \mathcal D^M V \, , 
\label{1D_semi_discretized_system_b}
\end{empheq}
\end{subequations}
\end{linenomath}
where $V$ and $\Sigma$ are discretized solution vectors; 
diagonal matrices $\boldsymbol \rho^M$ and $\boldsymbol \beta^N$ are the discrete correspondences of $\rho$ and $\beta$ in equation \eqref{1D_wave_system}, respectively; 
matrices $\mathcal D^N$ and $\mathcal D^M$ are finite difference operators approximating $\frac{\partial}{\partial x}$, which apply on vectors defined on the $N$- and $M$-grid, respectively; 
finally, diagonal matrices $\mathcal A^M$ and $\mathcal A^N$ are referred to as norm matrices, which are formally redundant in the above equations, but will play an important role in the upcoming discussion.

Analogously to equation \eqref{1D_continuous_energy}, the discrete energy associated with the above semi-discretized system is defined as:
\begin{equation}
\label{1D_discrete_energy}
E \ = \ \tfrac{1}{2} V^T \! \Big(\mathcal A^M \boldsymbol \rho^M \Big) V 
\ + \ \tfrac{1}{2} \Sigma^T \! \Big(\mathcal A^N \boldsymbol \beta^N \Big) \Sigma \, .
\end{equation}
Taking the time derivative of equation \eqref{1D_discrete_energy} and substituting the equations \eqref{1D_semi_discretized_system_a} and \eqref{1D_semi_discretized_system_b}, we arrive at:
\begin{equation}
\label{1D_discrete_energy_time_derivative}
\frac{d E}{d t} \ = \ \Sigma^T \left[ \mathcal A^N \mathcal D^M \ + \ \left( \mathcal A^M \mathcal D^N \right)^T \right] V \, .
\end{equation} 
To simplify the notation, we introduce matrix $Q$, which is defined as 
\begin{equation}
\label{1D_Q_definition}
Q \ = \ \mathcal A^N \mathcal D^M \ + \ ( \mathcal A^M \mathcal D^N )^T \, .
\end{equation}

For $\tfrac{d E}{d t}$ in equation \eqref{1D_discrete_energy_time_derivative} to mimic $\tfrac{d \mathscr E}{d t}$ in equation \eqref{1D_continuous_energy_time_derivative}, we ask $Q$ to satisfy the following property:
\begin{equation}
\label{1D_Q_property}
Q \ = \ - \ \mathcal E^L \left(\mathcal P^L\right)^T + \ \mathcal E^R \left(\mathcal P^R\right)^T,
\end{equation}
where $\mathcal E^L$ and $\mathcal E^R$ are the canonical basis vectors whose first and last components are 1, respectively;
$\mathcal P^L$ and $\mathcal P^R$ are projection operators that can be applied on vector $V$ to provide approximations to $v(x_L)$ and $v(x_R)$, respectively.

With the above property, equation \eqref{1D_discrete_energy_time_derivative} can be written as:
\begin{equation}
\label{1D_discrete_energy_time_derivative_boundary}
\frac{d E}{d t} \ = \ - \ \left( \mathcal E^L \right)^T \! \Sigma \cdot \left( \mathcal P^L \right)^T \! V 
                  \ + \ \left( \mathcal E^R \right)^T \! \Sigma \cdot \left( \mathcal P^R \right)^T \! V \, .
\end{equation} 
Realizing that 
$\left( \mathcal E^L \right)^T \! \Sigma$, 
$\left( \mathcal P^L \right)^T \! V$, 
$\left( \mathcal E^R \right)^T \! \Sigma$,
and $\left( \mathcal P^R \right)^T \! V$
are approximations to $\sigma(x_L)$, $v(x_L)$, $\sigma(x_R)$, and $v(x_R)$, respectively, we have that the discrete energy $E$ in equation \eqref{1D_discrete_energy_time_derivative} mimics the behavior of $\mathscr E$ in equation \eqref{1D_continuous_energy_time_derivative}.

The property of $Q$ in equation \eqref{1D_Q_property} is attained by carefully constructing the SBP operators $\mathcal A^N$, $\mathcal D^M$, $\mathcal A^M$, and $\mathcal D^N$ appearing in its definition; see equation \eqref{1D_Q_definition}.
In practice, such SBP operators are often constructed with the assistance of symbolic computing software.
One may refer to \cite[p.~6]{gao2019sbp} for an example of such SBP operators, which will be used throughout this study. 
Validation of these operators can be found in, e.g., \cite[Figure 3]{gao2020explicit} and \cite[Appendix A]{gao2019sbp}.
These operators are included in \ref{appendix_SBP_operators} to make this work self-contained.
From equation \eqref{SBP_matrices_1D_D}, we can see that in the interior, the finite difference operators $\mathcal D^M$ and $\mathcal D^N$ are the same as the standard fourth-order staggered grid stencil
$\nicefrac{[\nicefrac{1}{24},\ \nicefrac{-9}{8},\ \nicefrac{9}{8},\ \nicefrac{-1}{24}]}{\Delta x}$ (see, e.g., \cite{levander1988fourth}), 
and that only the stencils near the boundaries are adapted.

The norm matrices $\mathcal A^N$ and $\mathcal A^M$, seemingly redundant in the semi-discretized system \eqref{1D_semi_discretized_system}, introduce extra degrees of freedom to allow $Q$ to have the desired property. 
They reduce to the grid spacing $\Delta x$ in the interior and vary only near the boundaries. 
In the above energy analysis, $\mathcal A^N$ and $\mathcal A^M$ play the role of the integral operator 
$\displaystyle \int_{x_L}^{x_R}$ in equation \eqref{1D_continuous_energy}, and can be understood as providing the quadrature rules; see \cite{hicken2013summation}. 
In this work, we only consider the case where $\mathcal A^N$ and $\mathcal A^M$ are diagonal, so that the products $\big( \mathcal A^M \boldsymbol \rho^M \big)$ and $\big( \mathcal A^N \boldsymbol \beta^N \big)$ are symmetric, and, consequently, the expression in equation \eqref{1D_discrete_energy} provides a valid definition of discrete energy.
A drawback of this choice of diagonal norm matrices is that the formal order of the stencils near the boundary can only be as high as half of that in the interior; see \cite{kreiss1974finite, strand1994summation}.

In the following, we use the term ``discrete energy analysis'' to refer to the procedure of taking the time derivative of the discrete energy, substituting the semi-discretized system to replace the temporal derivatives of the solution variables, and finally, invoking the properties built in the SBP operators to reduce the result to boundaries or interfaces.

The free surface boundary conditions $\sigma(x_L)=\sigma(x_R)=0$ can be imposed by appending penalty terms to the semi-discretized system \eqref{1D_semi_discretized_system}, leading to the following modified system:
\begin{linenomath}
\begin{subequations}
\label{1D_semi_discretized_system_modified}
\begin{empheq}[left=\empheqlbrace]{alignat = 2}
\displaystyle \enskip \mathcal A^M \boldsymbol \rho^M \frac{d V}{d t} 
\enskip &= \enskip \displaystyle 
\! \mathcal A^M \mathcal D^N \Sigma 
\ + \ \mathcal P^L \left[ \left( \mathcal E^L \right)^T \! \Sigma - \bf 0 \right] 
\ - \ \mathcal P^R \left[ \left( \mathcal E^R \right)^T \! \Sigma - \bf 0 \right] \, ; \label{1D_semi_discretized_system_modified_a} \\[0.25ex]
\displaystyle \enskip \mathcal A^N \boldsymbol \beta^N \frac{d \Sigma}{d t} 
\enskip &= \enskip \displaystyle 
\! \mathcal A^N \mathcal D^M V \, .
\label{1D_semi_discretized_system_modified_b}
\end{empheq}
\end{subequations}
\end{linenomath}
Following the aforementioned procedure of discrete energy analysis, one can verify that the discrete energy $E$ as defined in equation \eqref{1D_discrete_energy} is conservative, i.e., $\tfrac{dE}{dt} = 0$, for the above modified system.

\subsection{SBP finite difference operators on 3D staggered grids} \label{SBP_operators_3D}
The 1D SBP operators introduced above serve as the building blocks of the 3D SBP operators used in this work.
We append subscripts $x$, $y$, and $z$ to the 1D operators $\mathcal A^N$, $\mathcal D^M$, $\mathcal A^M$, and $\mathcal D^N$ to distinguish the directions.
By design, these 1D operators satisfy the following properties:
\begin{linenomath}
\begin{subequations}
\label{Q_definitions}
\begin{empheq}{alignat = 3}
\displaystyle Q_x \ &= \ \mathcal A^N_x \mathcal D^M_x + ( \mathcal A^M_x \mathcal D^N_x )^T
\displaystyle \enskip &= \ - \ \mathcal E^L_x \left(\mathcal P^L_x\right)^T \, + \ \mathcal E^R_x \left(\mathcal P^R_x\right)^T; 
\label{Q_definitions_x} \\
\displaystyle Q_y \ &= \ \mathcal A^N_y \mathcal D^M_y + ( \mathcal A^M_y \mathcal D^N_y )^T
\displaystyle \enskip &= \ - \ \mathcal E^L_y \left(\mathcal P^L_y\right)^T \, + \ \mathcal E^R_y \left(\mathcal P^R_y\right)^T; 
\label{Q_definitions_y} \\
\displaystyle Q_z \ &= \ \mathcal A^N_z \mathcal D^M_z + ( \mathcal A^M_z \mathcal D^N_z )^T
\displaystyle \enskip &= \ - \ \mathcal E^L_z \left(\mathcal P^L_z\right)^T \, + \ \mathcal E^R_z \left(\mathcal P^R_z\right)^T. 
\label{Q_definitions_z}
\end{empheq}
\end{subequations}
\end{linenomath}
Furthermore, we use $\mathcal I^N_x$, $\mathcal I^M_x$, $\mathcal I^N_y$, $\mathcal I^M_y$, $\mathcal I^N_z$, and $\mathcal I^M_z$ to denote the identity matrices corresponding to the $N$- and $M$-grids of the $x$-, $y$-, and $z$-directions, respectively.

The 3D norm matrices are constructed as the tensor products of the 1D norm matrices as follows:
\begin{equation}
\label{SBP_norm_matrices_3D}
\arraycolsep=.75em\def\arraystretch{1}
\hspace{-.25em}
\begin{array}{lll}
\mathcal A^{V_x} = \mathcal A^M_x \otimes \mathcal A^N_y \otimes \mathcal A^N_z, & 
\mathcal A^{\Sigma_{xy}} = \mathcal A^M_x \otimes \mathcal A^M_y \otimes \mathcal A^N_z, &
\mathcal A^{\Sigma_{xx}} = \mathcal A^N_x \otimes \mathcal A^N_y \otimes \mathcal A^N_z, \\
\mathcal A^{V_y} = \mathcal A^N_x \otimes \mathcal A^M_y \otimes \mathcal A^N_z, & 
\mathcal A^{\Sigma_{yz}} = \mathcal A^N_x \otimes \mathcal A^M_y \otimes \mathcal A^M_z, &
\mathcal A^{\Sigma_{yy}} = \mathcal A^N_x \otimes \mathcal A^N_y \otimes \mathcal A^N_z, \\
\mathcal A^{V_z} = \mathcal A^N_x \otimes \mathcal A^N_y \otimes \mathcal A^M_z, & 
\mathcal A^{\Sigma_{xz}} = \mathcal A^M_x \otimes \mathcal A^N_y \otimes \mathcal A^M_z, &
\mathcal A^{\Sigma_{zz}} = \mathcal A^N_x \otimes \mathcal A^N_y \otimes \mathcal A^N_z. \\
\end{array}
\end{equation}
The 3D finite difference operators are constructed as the tensor products of the 1D finite difference operators and the 1D identity matrices as follows:
\begin{equation}
\label{SBP_difference_operators_3D}
\arraycolsep=.75em\def\arraystretch{1}
\hspace{-.25em}
\begin{array}{lll}
\mathcal D^{\Sigma_{xx}}_x = \mathcal D^N_x \otimes \mathcal I^N_y \otimes \mathcal I^N_z, & 
\mathcal D^{\Sigma_{xy}}_y = \mathcal I^M_x \otimes \mathcal D^M_y \otimes \mathcal I^N_z, & 
\mathcal D^{\Sigma_{xz}}_z = \mathcal I^M_x \otimes \mathcal I^N_y \otimes \mathcal D^M_z, \\
\mathcal D^{\Sigma_{xy}}_x = \mathcal D^M_x \otimes \mathcal I^M_y \otimes \mathcal I^N_z, & 
\mathcal D^{\Sigma_{yy}}_y = \mathcal I^N_x \otimes \mathcal D^N_y \otimes \mathcal I^N_z, & 
\mathcal D^{\Sigma_{yz}}_z = \mathcal I^N_x \otimes \mathcal I^M_y \otimes \mathcal D^M_z, \\
\mathcal D^{\Sigma_{xz}}_x = \mathcal D^M_x \otimes \mathcal I^N_y \otimes \mathcal I^M_z, & 
\mathcal D^{\Sigma_{yz}}_y = \mathcal I^N_x \otimes \mathcal D^M_y \otimes \mathcal I^M_z, & 
\mathcal D^{\Sigma_{zz}}_z = \mathcal I^N_x \otimes \mathcal I^N_y \otimes \mathcal D^N_z, \\[1.ex]
\mathcal D^{V_x}_x = \mathcal D^M_x \otimes \mathcal I^N_y \otimes \mathcal I^N_z, & 
\mathcal D^{V_x}_y = \mathcal I^M_x \otimes \mathcal D^N_y \otimes \mathcal I^N_z, & 
\mathcal D^{V_x}_z = \mathcal I^M_x \otimes \mathcal I^N_y \otimes \mathcal D^N_z, \\
\mathcal D^{V_y}_x = \mathcal D^N_x \otimes \mathcal I^M_y \otimes \mathcal I^N_z, & 
\mathcal D^{V_y}_y = \mathcal I^N_x \otimes \mathcal D^M_y \otimes \mathcal I^N_z, & 
\mathcal D^{V_y}_z = \mathcal I^N_x \otimes \mathcal I^M_y \otimes \mathcal D^N_z, \\
\mathcal D^{V_z}_x = \mathcal D^N_x \otimes \mathcal I^N_y \otimes \mathcal I^M_z, & 
\mathcal D^{V_z}_y = \mathcal I^N_x \otimes \mathcal D^N_y \otimes \mathcal I^M_z, & 
\mathcal D^{V_z}_z = \mathcal I^N_x \otimes \mathcal I^N_y \otimes \mathcal D^M_z. \\
\end{array}
\end{equation}

Supposing that the 3D wave system described in equation \eqref{Elastic_wave_equation_succinct} is defined over a cuboid $\Omega = (x_L,x_R) \times (y_L,y_R) \times (z_L,z_R)$, its associated physical energy can be expressed as:
\begin{equation}
\label{3D_continuous_energy}
\mathscr E = \frac{1}{2} \int_{\Omega} \left( \rho v_i v_i + \sigma_{ij} S_{ijkl} \sigma_{kl} \right) d_\Omega \, , 
\end{equation}
where $S$ stands for the compliance tensor, which is the inverse of the stiffness tensor $C$ from equation \eqref{Elastic_wave_equation_succinct_b}. 
Taking the time derivative of equation \eqref{3D_continuous_energy}, substituting the equations \eqref{Elastic_wave_equation_succinct_a} and \eqref{Elastic_wave_equation_succinct_b}, and applying the divergence theorem, we arrive at 
\begin{equation}
\label{3D_continuous_energy_time_derivative}
\frac{d \mathscr E}{d t} = \int_{\partial \Omega} v_i \sigma_{ij} n_j \ d_{\partial \Omega}\, ,
\end{equation}
where $n_j$ stands for component of the outward normal vector on the boundary.
From equation \eqref{3D_continuous_energy_time_derivative}, we have that the time derivative of the physical energy $\mathscr E$ reduces to boundary data only, analogously to the 1D case as in equation \eqref{1D_continuous_energy_time_derivative}.

Replacing the spatial derivatives in equation \eqref{Elastic_wave_equation_succinct} with the 3D SBP operators introduced above, the resulting discretized system can be expressed as:
\begin{linenomath}
\begin{subequations}
\label{Semi_discretized_system_succinct}
\begin{empheq}[left=\empheqlbrace\,]{alignat = 2}
\displaystyle \mathcal A^{V_i} \boldsymbol \rho^{V_i} \frac{d V_i}{d t} & \ = \ \displaystyle \mathcal A^{V_i} \mathcal D^{\Sigma_{ij}}_j \Sigma_{ij} \, ; \label{Semi_discretized_system_succinct_a}
\\[0.25ex]
\displaystyle \mathcal A^{\Sigma_{ij}} \boldsymbol S^{\Sigma_{ij}}_{ijkl} \frac{d \Sigma_{kl}}{d t}
&\ = \ \mathcal A^{\Sigma_{ij}} \tfrac{1}{2} \left( \mathcal D^{V_i}_j V_i + \mathcal D^{V_j}_i V_j \right) ,\label{Semi_discretized_system_succinct_b}
\end{empheq}
\end{subequations}
\end{linenomath}
where the Einstein summation convention applies to subscript indices only 
(e.g., for the expression $\mathcal A^{V_i} \mathcal D^{\Sigma_{ij}}_j \Sigma_{ij}$ on the right-hand side of equation \eqref{Semi_discretized_system_succinct_a}, only the index $j$ is summed over).

The discrete energy associated with system \eqref{Semi_discretized_system_succinct} is defined as
\begin{equation}
\label{3D_discrete_energy}
E = \tfrac{1}{2} V_i^T \left(\mathcal A^{V_i} \boldsymbol \rho^{V_i} \right) V_i 
\ + \ \tfrac{1}{2} \Sigma_{ij}^T \left( \mathcal A^{\Sigma_{ij}} \boldsymbol S^{\Sigma_{ij}}_{ijkl} \right) \Sigma_{kl} \, .
\end{equation}
Following the procedure of discrete energy analysis introduced in the previous section, it can be shown that 
\begin{equation}
\label{3D_discrete_energy_time_derivative_boundary}
\arraycolsep=0.5pt\def\arraystretch{0.875}
\begin{array}{l}
\left.
\begin{array}{rrrllcrr}
\displaystyle \frac{d E}{d t}
& \ = & \ - \,\, & \displaystyle 
\big[ \Sigma_{xx}^T \left(\mathcal E^L_x \otimes \mathcal I^N_y \otimes \mathcal I^N_z \right) \big] 
& \cdot & \big[\mathcal A^N_y \otimes \mathcal A^N_z \big] & \cdot & 
\big[ \left(\mathcal P^L_x \otimes \mathcal I^N_y \otimes \mathcal I^N_z \right)^T V_x \big] \phantom{,} \\
& & \ + \,\, & \displaystyle 
\big[ \Sigma_{xx}^T \left(\mathcal E^R_x \otimes \mathcal I^N_y \otimes \mathcal I^N_z \right) \big] 
& \cdot & \big[\mathcal A^N_y \otimes \mathcal A^N_z \big] & \cdot & 
\big[ \left(\mathcal P^R_x \otimes \mathcal I^N_y \otimes \mathcal I^N_z \right)^T V_x \big] \phantom{,} \\
& & \ - \,\, & \displaystyle 
\big[ \Sigma_{xy}^T \left(\mathcal P^L_x \otimes \mathcal I^M_y \otimes \mathcal I^N_z \right) \big] 
& \cdot & \big[\mathcal A^M_y \otimes \mathcal A^N_z \big] & \cdot & 
\big[ \left(\mathcal E^L_x \otimes \mathcal I^M_y \otimes \mathcal I^N_z \right)^T V_y \big] \phantom{,} \\
& & \ + \,\, & \displaystyle 
\big[ \Sigma_{xy}^T \left(\mathcal P^R_x \otimes \mathcal I^M_y \otimes \mathcal I^N_z \right) \big] 
& \cdot & \big[\mathcal A^M_y \otimes \mathcal A^N_z \big] & \cdot & 
\big[ \left(\mathcal E^R_x \otimes \mathcal I^M_y \otimes \mathcal I^N_z \right)^T V_y \big] \phantom{,} \\
& & \ - \,\, & \displaystyle 
\big[ \Sigma_{xz}^T \left(\mathcal P^L_x \otimes \mathcal I^N_y \otimes \mathcal I^M_z \right) \big] 
& \cdot & \big[\mathcal A^N_y \otimes \mathcal A^M_z \big] & \cdot & 
\big[ \left(\mathcal E^L_x \otimes \mathcal I^N_y \otimes \mathcal I^M_z \right)^T V_z \big] \phantom{,} \\
& & \ + \,\, & \displaystyle 
\big[ \Sigma_{xz}^T \left(\mathcal P^R_x \otimes \mathcal I^N_y \otimes \mathcal I^M_z \right) \big] 
& \cdot & \big[\mathcal A^N_y \otimes \mathcal A^M_z \big] & \cdot & 
\big[ \left(\mathcal E^R_x \otimes \mathcal I^N_y \otimes \mathcal I^M_z \right)^T V_z \big] \phantom{,}
\end{array}
\right\} yz\text{-faces}
\\
\left.
\begin{array}{rrrllcrr}
\displaystyle \phantom{\frac{d E}{d t}}
& \ \phantom{=} & \ - \,\, & \displaystyle
\big[ \Sigma_{xy}^T \left(\mathcal I^M_x \otimes \mathcal P^L_y \otimes \mathcal I^N_z \right) \big] 
& \cdot & \big[\mathcal A^M_x \otimes \mathcal A^N_z \big] & \cdot & 
\big[ \left(\mathcal I^M_x \otimes \mathcal E^L_y \otimes \mathcal I^N_z \right)^T V_x \big] \phantom{,} \\
& & \ + \,\, & \displaystyle 
\big[ \Sigma_{xy}^T \left(\mathcal I^M_x \otimes \mathcal P^R_y \otimes \mathcal I^N_z \right) \big] 
& \cdot & \big[\mathcal A^M_x \otimes \mathcal A^N_z \big] & \cdot & 
\big[ \left(\mathcal I^M_x \otimes \mathcal E^R_y \otimes \mathcal I^N_z \right)^T V_x \big] \phantom{,} \\
& & \ - \,\, & \displaystyle 
\big[ \Sigma_{yy}^T \left(\mathcal I^N_x \otimes \mathcal E^L_y \otimes \mathcal I^N_z \right) \big] 
& \cdot & \big[\mathcal A^N_x \otimes \mathcal A^N_z \big] & \cdot & 
\big[ \left(\mathcal I^N_x \otimes \mathcal P^L_y \otimes \mathcal I^N_z \right)^T V_y \big] \phantom{,} \\
& & \ + \,\, & \displaystyle 
\big[ \Sigma_{yy}^T \left(\mathcal I^N_x \otimes \mathcal E^R_y \otimes \mathcal I^N_z \right) \big] 
& \cdot & \big[\mathcal A^N_x \otimes \mathcal A^N_z \big] & \cdot & 
\big[ \left(\mathcal I^N_x \otimes \mathcal P^R_y \otimes \mathcal I^N_z \right)^T V_y \big] \phantom{,} \\
& & \ - \,\, & \displaystyle 
\big[ \Sigma_{yz}^T \left(\mathcal I^N_x \otimes \mathcal P^L_y \otimes \mathcal I^M_z \right) \big] 
& \cdot & \big[\mathcal A^N_x \otimes \mathcal A^M_z \big] & \cdot & 
\big[ \left(\mathcal I^N_x \otimes \mathcal E^L_y \otimes \mathcal I^M_z \right)^T V_z \big] \phantom{,} \\
& & \ + \,\, & \displaystyle 
\big[ \Sigma_{yz}^T \left(\mathcal I^N_x \otimes \mathcal P^R_y \otimes \mathcal I^M_z \right) \big] 
& \cdot & \big[\mathcal A^N_x \otimes \mathcal A^M_z \big] & \cdot & 
\big[ \left(\mathcal I^N_x \otimes \mathcal E^R_y \otimes \mathcal I^M_z \right)^T V_z \big] \phantom{,}
\end{array}
\right\} xz\text{-faces}
\\
\left.
\begin{array}{rrrllcrr}
\phantom{\displaystyle \frac{d E}{d t}}
& \ \phantom{=} & \ - \,\, & \displaystyle 
\big[ \Sigma_{xz}^T \left(\mathcal I^M_x \otimes \mathcal I^N_y \otimes \mathcal P^L_z \right) \big] 
& \cdot & \big[\mathcal A^M_x \otimes \mathcal A^N_y \big] & \cdot & 
\big[ \left(\mathcal I^M_x \otimes \mathcal I^N_y \otimes \mathcal E^L_z \right)^T V_x \big] \phantom{,} \\
& & \ + \,\, & \displaystyle 
\big[ \Sigma_{xz}^T \left(\mathcal I^M_x \otimes \mathcal I^N_y \otimes \mathcal P^R_z \right) \big] 
& \cdot & \big[\mathcal A^M_x \otimes \mathcal A^N_y \big] & \cdot & 
\big[ \left(\mathcal I^M_x \otimes \mathcal I^N_y \otimes \mathcal E^R_z \right)^T V_x \big] \phantom{,} \\
& & \ - \,\, & \displaystyle 
\big[ \Sigma_{yz}^T \left(\mathcal I^N_x \otimes \mathcal I^M_y \otimes \mathcal P^L_z \right) \big] 
& \cdot & \big[\mathcal A^N_x \otimes \mathcal A^M_y \big] & \cdot & 
\big[ \left(\mathcal I^N_x \otimes \mathcal I^M_y \otimes \mathcal E^L_z \right)^T V_y \big] \phantom{,} \\
& & \ + \,\, & \displaystyle 
\big[ \Sigma_{yz}^T \left(\mathcal I^N_x \otimes \mathcal I^M_y \otimes \mathcal P^R_z \right) \big] 
& \cdot & \big[\mathcal A^N_x \otimes \mathcal A^M_y \big] & \cdot & 
\big[ \left(\mathcal I^N_x \otimes \mathcal I^M_y \otimes \mathcal E^R_z \right)^T V_y \big] \phantom{,} \\
& & \ - \,\, & \displaystyle 
\big[ \Sigma_{zz}^T \left(\mathcal I^N_x \otimes \mathcal I^N_y \otimes \mathcal E^L_z \right) \big] 
& \cdot & \big[\mathcal A^N_x \otimes \mathcal A^N_y \big] & \cdot & 
\big[ \left(\mathcal I^N_x \otimes \mathcal I^N_y \otimes \mathcal P^L_z \right)^T V_z \big] \phantom{,} \\
& & \ + \,\, & \displaystyle 
\big[ \Sigma_{zz}^T \left(\mathcal I^N_x \otimes \mathcal I^N_y \otimes \mathcal E^R_z \right) \big] 
& \cdot & \big[\mathcal A^N_x \otimes \mathcal A^N_y \big] & \cdot & 
\big[ \left(\mathcal I^N_x \otimes \mathcal I^N_y \otimes \mathcal P^R_z \right)^T V_z \big] \phantom{,}
\end{array}
\right\} xy\text{-faces,}
\end{array}
\end{equation} 
which again reduces to boundary data only.

In order to see this more clearly, let us take the first line of equation \eqref{3D_discrete_energy_time_derivative_boundary} as an example.
Recognizing that $\big[ \Sigma_{xx}^T \left(\mathcal E^L_x \otimes \mathcal I^N_y \otimes \mathcal I^N_z \right) \big]$ is the restriction of $\Sigma_{xx}$ on the $yz$-face at $x \!=\! x_L$, and that $\big[ \left(\mathcal P^L_x \otimes \mathcal I^N_y \otimes \mathcal I^N_z \right)^T V_x \big]$ is the projection of $V_x$ on the same face, we have that the expression in this line is an approximation to 
\begin{equation}
\label{yz_face_xL}
\displaystyle - \int_{y_L}^{y_R} \int_{z_L}^{z_R} \sigma_{xx}(x_L,y,z)\, v_x(x_L,y,z)\ dy dz \, ,
\end{equation}
with $\big[\mathcal A^N_y \otimes \mathcal A^N_z \big]$ playing the role of the integral operators in equation \eqref{yz_face_xL}.
Similar interpretations apply to the other lines in equation \eqref{3D_discrete_energy_time_derivative_boundary}.
Altogether, we have that the right-hand side of equation \eqref{3D_discrete_energy_time_derivative_boundary} is an approximation to the right-hand side of equation \eqref{3D_continuous_energy_time_derivative}.
In other words, the discrete energy $E$ as defined in equation \eqref{3D_discrete_energy} mimics the behavior of the continuous energy $\mathscr E$ in equation \eqref{3D_continuous_energy_time_derivative}.

Moreover, similar to the 1D case, the free surface boundary condition can also be imposed by appending penalty terms to equation \eqref{Semi_discretized_system_succinct}.
Taking the two $xy$-faces at $z \!=\! z_L$ and $z \!=\! z_R$ for example, the free surface boundary condition states that $\sigma_{xz} \!=\! \sigma_{yz} \!=\! \sigma_{zz} \!=\! 0$.
These conditions can be imposed by modifying the updating formulas for the velocity components in equation \eqref{Semi_discretized_system_succinct} as follows: 
\begin{linenomath}
\begin{subequations}
\label{3D_semi_discretized_system_modified_velocities}
\begin{empheq}[left=\empheqlbrace\,]{alignat = 2}
\displaystyle \mathcal A^{V_x} \boldsymbol \rho^{V_x} \frac{d V_x}{d t} & \ = \ \displaystyle \mathcal A^{V_x} 
\left( 
   \mathcal D^{\Sigma_{xx}}_x \Sigma_{xx} 
+ \mathcal D^{\Sigma_{xy}}_y \Sigma_{xy} 
+ \mathcal D^{\Sigma_{xz}}_z \Sigma_{xz} 
\right) \label{3D_semi_discretized_system_modified_velocities_Vx}
\\
& \ + \ \displaystyle 
\left(\mathcal I^M_x \otimes \mathcal I^N_y \otimes \mathcal E^L_z \right) 
\cdot \left(\mathcal A^M_x \otimes \mathcal A^N_y \right) \cdot
\left[ \big(\mathcal I^M_x \otimes \mathcal I^N_y \otimes \mathcal P^L_z \big)^T \Sigma_{xz} - \boldsymbol 0 \right] \nonumber
\\
& \ - \ \displaystyle 
\left( \mathcal I^M_x \otimes \mathcal I^N_y \otimes \mathcal E^R_z \right)
\cdot \left( \mathcal A^M_x \otimes \mathcal A^N_y \right) \cdot
\left[ \big(\mathcal I^M_x \otimes \mathcal I^N_y \otimes \mathcal P^R_z \big)^T \Sigma_{xz} - \boldsymbol 0 \right] ; \nonumber
\\[0.25ex]
\displaystyle \mathcal A^{V_y} \boldsymbol \rho^{V_y} \frac{d V_y}{d t} & \ = \ \displaystyle \mathcal A^{V_y} 
\left( 
   \mathcal D^{\Sigma_{xy}}_x \Sigma_{xy} 
+ \mathcal D^{\Sigma_{yy}}_y \Sigma_{yy} 
+ \mathcal D^{\Sigma_{yz}}_z \Sigma_{yz} 
\right) \label{3D_semi_discretized_system_modified_velocities_Vy}
\\
& \ + \ \displaystyle 
\left( \mathcal I^N_x \otimes \mathcal I^M_y \otimes \mathcal E^L_z \right)
\cdot \left( \mathcal A^N_x \otimes \mathcal A^M_y \right) \cdot 
\left[ \big( \mathcal I^N_x \otimes \mathcal I^M_y \otimes \mathcal P^L_z \big)^T \Sigma_{yz} - \boldsymbol 0 \right] \nonumber
\\
& \ - \ \displaystyle 
\left( \mathcal I^N_x \otimes \mathcal I^M_y \otimes \mathcal E^R_z \right)
\cdot \left( \mathcal A^N_x \otimes \mathcal A^M_y \right) \cdot 
\left[ \big( \mathcal I^N_x \otimes \mathcal I^M_y \otimes \mathcal P^R_z \big)^T \Sigma_{yz} - \boldsymbol 0 \right] ; \nonumber
\\[0.25ex]
\displaystyle \mathcal A^{V_z} \boldsymbol \rho^{V_z} \frac{d V_z}{d t} & \ = \ \displaystyle \mathcal A^{V_z} 
\left( 
   \mathcal D^{\Sigma_{xz}}_x \Sigma_{xz} 
+ \mathcal D^{\Sigma_{yz}}_y \Sigma_{yz} 
+ \mathcal D^{\Sigma_{zz}}_z \Sigma_{zz} 
\right) \label{3D_semi_discretized_system_modified_velocities_Vz}
\\
& \ + \ \displaystyle 
\left( \mathcal I^N_x \otimes \mathcal I^N_y \otimes \mathcal P^L_z \right)
\cdot \left( \mathcal A^N_x \otimes \mathcal A^N_y \right) \cdot
\left[ \big( \mathcal I^N_x \otimes \mathcal I^N_y \otimes \mathcal E^L_z \big)^T \Sigma_{zz} - \boldsymbol 0 \right] \nonumber
\\
& \ - \ \displaystyle 
\left( \mathcal I^N_x \otimes \mathcal I^N_y \otimes \mathcal P^R_z \right)
\cdot \left( \mathcal A^N_x \otimes \mathcal A^N_y \right) \cdot
\left[ \big( \mathcal I^N_x \otimes \mathcal I^N_y \otimes \mathcal E^R_z \big)^T \Sigma_{zz} - \boldsymbol 0 \right] , \nonumber
\end{empheq}
\end{subequations}
\end{linenomath}
where the penalty terms appended to equations \eqref{3D_semi_discretized_system_modified_velocities_Vx}, \eqref{3D_semi_discretized_system_modified_velocities_Vy}, and \eqref{3D_semi_discretized_system_modified_velocities_Vz} correspond to conditions $\sigma_{xz} \!=\! 0$, $\sigma_{yz} \!=\! 0$, and $\sigma_{zz} \!=\! 0$, respectively.
It can be verified that the penalty terms introduced above would cancel out the last six lines in equation \eqref{3D_discrete_energy_time_derivative_boundary} (i.e., those included in the last bracket that corresponds to the $xy$-faces) in discrete energy analysis.

We note here that appending penalty terms to the discretized system can be viewed as modifying the corresponding spatial derivative approximations. 
This realization will help simplify the upcoming discussion.
For example, appending the penalty terms in equations \eqref{3D_semi_discretized_system_modified_velocities_Vx}, \eqref{3D_semi_discretized_system_modified_velocities_Vy}, and \eqref{3D_semi_discretized_system_modified_velocities_Vz} can be viewed as replacing 
$\mathcal D^{\Sigma_{xz}}_z \Sigma_{xz}$, $\mathcal D^{\Sigma_{yz}}_z \Sigma_{yz}$, and $\mathcal D^{\Sigma_{zz}}_z \Sigma_{zz}$, respectively, with their modified spatial derivative approximations defined as follows:
\begin{linenomath}
\begin{subequations}
\label{3D_modified_spatial_derivative_approximations_velocities}
\begin{empheq}[left=\!\!\!\!,right=\!\!\!\!]{alignat = 2}
\overline{\mathcal D^{\Sigma_{xz}}_z \Sigma_{xz}} 
\, = \, 
\mathcal D^{\Sigma_{xz}}_z \Sigma_{xz} 
& \, + \, 
\left\{ \mathcal I^M_x \otimes \mathcal I^N_y \otimes \left[ \big( \mathcal A^N_z \big)^{-1} \mathcal E^L_z \right] \right\}
\cdot
\left[ \big(\mathcal I^M_x \otimes \mathcal I^N_y \otimes \mathcal P^L_z \big)^T \Sigma_{xz} - \boldsymbol 0 \right]
\\
& \, - \,
\left\{ \mathcal I^M_x \otimes \mathcal I^N_y \otimes \left[ \big( \mathcal A^N_z \big)^{-1} \mathcal E^R_z \right] \right\}
\cdot
\left[ \big(\mathcal I^M_x \otimes \mathcal I^N_y \otimes \mathcal P^R_z \big)^T \Sigma_{xz} - \boldsymbol 0 \right] ; \nonumber
\\[0.5ex]
\overline{\mathcal D^{\Sigma_{yz}}_z \Sigma_{yz}} 
\, = \, 
\mathcal D^{\Sigma_{yz}}_z \Sigma_{yz} 
& \, + \, 
\left\{ \mathcal I^N_x \otimes \mathcal I^M_y \otimes \left[ \big( \mathcal A^N_z \big)^{-1} \mathcal E^L_z \right] \right\}
\cdot
\left[ \big(\mathcal I^N_x \otimes \mathcal I^M_y \otimes \mathcal P^L_z \big)^T \Sigma_{yz} - \boldsymbol 0 \right]
\\
& \, - \, 
\left\{ \mathcal I^N_x \otimes \mathcal I^M_y \otimes \left[ \big( \mathcal A^N_z \big)^{-1} \mathcal E^R_z \right] \right\}
\cdot
\left[ \big(\mathcal I^N_x \otimes \mathcal I^M_y \otimes \mathcal P^R_z \big)^T \Sigma_{yz} - \boldsymbol 0 \right] ; \nonumber
\\[0.5ex]
\overline{\mathcal D^{\Sigma_{zz}}_z \Sigma_{zz}} 
\, = \, 
\mathcal D^{\Sigma_{zz}}_z \Sigma_{zz} 
& \, + \, 
\left\{ \mathcal I^N_x \otimes \mathcal I^N_y \otimes \left[ \big( \mathcal A^M_z \big)^{-1} \mathcal P^L_z \right] \right\}
\cdot
\left[ \big(\mathcal I^N_x \otimes \mathcal I^N_y \otimes \mathcal E^L_z \big)^T \Sigma_{zz} - \boldsymbol 0 \right]
\\
& \, - \, 
\left\{ \mathcal I^N_x \otimes \mathcal I^N_y \otimes \left[ \big( \mathcal A^M_z \big)^{-1} \mathcal P^R_z \right] \right\}
\cdot
\left[ \big(\mathcal I^N_x \otimes \mathcal I^N_y \otimes \mathcal E^R_z \big)^T \Sigma_{zz} - \boldsymbol 0 \right] , \nonumber
\end{empheq}
\end{subequations}
\end{linenomath}
where inverses of the norm matrices $\mathcal A^{V_x}$, $\mathcal A^{V_y}$, and $\mathcal A^{V_z}$ have been applied to the penalty terms in equations \eqref{3D_semi_discretized_system_modified_velocities_Vx}, \eqref{3D_semi_discretized_system_modified_velocities_Vy}, and \eqref{3D_semi_discretized_system_modified_velocities_Vz}, respectively.
Using the modified spatial derivative approximations introduced above, the norm matrices $\mathcal A^{V_x}$, $\mathcal A^{V_y}$, and $\mathcal A^{V_z}$ can now be omitted from the discretized system, though their impact remains through the modified spatial derivative approximations.

\subsection{Interface conditions and interpolation operators} \label{interface_and_interpolation}
At the interface between two neighboring layers, we seek to impose the following interface conditions (see, e.g., \cite[p.~52]{stein2009introduction}): 
\begin{linenomath}
\begin{subequations}
\label{interface_conditions}
\begin{empheq}[left=\empheqlbrace\,]{alignat = 3}
 v_i^+ & = v_i^- & ; \label{interface_conditions_a}
 \\
 \sigma^+_{ij} n^+_j + \sigma^-_{ij} n^-_j & = 0 & , \label{interface_conditions_b}
\end{empheq}
\end{subequations}
\end{linenomath}
where symbols $^+$ and $^-$ are appended to distinguish quantities from the plus and minus sides of the interface.
equation \eqref{interface_conditions_a} stems from continuity of the elastic medium, i.e., no overlap or tear; equation \eqref{interface_conditions_b} stems from Newton’s third law. 
For the discretization setting considered in this work, where the layer interfaces are horizontal, the outward normal vectors $n^+$ and $n^-$ in equation \eqref{interface_conditions} are simply $[0,0,-1]^T$ and $[0,0,1]^T$, respectively, and the interface conditions in equation \eqref{interface_conditions} reduce to continuity in solution variables $v_x$, $v_y$, $v_z$, $\sigma_{xz}$, $\sigma_{yz}$, and $\sigma_{zz}$, with the remaining solution variables uninvolved.

These interface conditions can also be imposed by appending proper penalty terms to the discretized system described in equation \eqref{Semi_discretized_system_succinct}. Effectively, this boils down to modifying the spatial derivative approximations as shown in equations \eqref{3D_modified_spatial_derivative_approximations_plus} and \eqref{3D_modified_spatial_derivative_approximations_minus}. 

\noindent (On the plus side:)
\begin{subequations}
\label{3D_modified_spatial_derivative_approximations_plus}
\begin{empheq}[left=\hspace{-5em},right=]{alignat = 3}
\mathcal D^{\Sigma^+_{xz}}_z \Sigma^+_{xz} 
& \, \longrightarrow \, &
\overline{\mathcal D^{\Sigma^+_{xz}}_z \Sigma^+_{xz}} 
\, = \, &
\mathcal D^{\Sigma^+_{xz}}_z \Sigma^+_{xz} 
\, + \, \tfrac{1}{2}
\left\{ \mathcal I^M_{x^+} \otimes \mathcal I^N_{y^+} \otimes \left[ \big( \mathcal A^N_{z^+} \big)^{-1} \mathcal E^L_{z^+} \right] \right\}
\label{3D_modified_spatial_derivative_approximations_plus_Sxz}
\\ 
& & & 
\cdot
\left\{ \big(\mathcal I^M_{x^+} \otimes \mathcal I^N_{y^+} \otimes \mathcal P^L_{z^+} \big)^T \Sigma^+_{xz} \right. 
- \, 
\left.\mathcal T^{\Sigma^-_{xz}} \left[ \big(\mathcal I^M_{x^-} \otimes \mathcal I^N_{y^-} \otimes \mathcal P^R_{z^-} \big)^T \Sigma^-_{xz} \right] \right\} ; \nonumber
\\
\mathcal D^{V^+_x}_z V^+_x 
& \, \longrightarrow \, &
\overline{\mathcal D^{V^+_x}_z V^+_x} 
\, = \, &
\mathcal D^{V^+_x}_z V^+_x
\, + \, \tfrac{1}{2}
\left\{ \mathcal I^M_{x^+} \otimes \mathcal I^N_{y^+} \otimes \left[ \big( \mathcal A^M_{z^+} \big)^{-1} \mathcal P^L_{z^+} \right] \right\}
\label{3D_modified_spatial_derivative_approximations_plus_Vx}
\\ 
& & & 
\cdot
\left\{ \big(\mathcal I^M_{x^+} \otimes \mathcal I^N_{y^+} \otimes \mathcal E^L_{z^+} \big)^T V^+_x \right. 
- \, 
\left.\mathcal T^{V^-_x} \left[ \big(\mathcal I^M_{x^-} \otimes \mathcal I^N_{y^-} \otimes \mathcal E^R_{z^-} \big)^T V^-_x \right] \right\} ; \nonumber
\\
\mathcal D^{\Sigma^+_{yz}}_z \Sigma^+_{yz} 
& \, \longrightarrow \, &
\overline{\mathcal D^{\Sigma^+_{yz}}_z \Sigma^+_{yz}} 
\, = \, &
\mathcal D^{\Sigma^+_{yz}}_z \Sigma^+_{yz} 
\, + \, \tfrac{1}{2}
\left\{ \mathcal I^N_{x^+} \otimes \mathcal I^M_{y^+} \otimes \left[ \big( \mathcal A^N_{z^+} \big)^{-1} \mathcal E^L_{z^+} \right] \right\}
\label{3D_modified_spatial_derivative_approximations_plus_Syz}
\\
& & &
\cdot
\left\{ \big(\mathcal I^N_{x^+} \otimes \mathcal I^M_{y^+} \otimes \mathcal P^L_{z^+} \big)^T \Sigma^+_{yz} \right. 
- \,
\left. \mathcal T^{\Sigma^-_{yz}} \left[ \big(\mathcal I^N_{x^-} \otimes \mathcal I^M_{y^-} \otimes \mathcal P^R_{z^-} \big)^T \Sigma^-_{yz} \right] \right\} ; \nonumber
\\
\mathcal D^{V^+_y}_z V^+_y 
& \, \longrightarrow \, &
\overline{\mathcal D^{V^+_y}_z V^+_y} 
\, = \, &
\mathcal D^{V^+_y}_z V^+_y 
\, + \, \tfrac{1}{2}
\left\{ \mathcal I^N_{x^+} \otimes \mathcal I^M_{y^+} \otimes \left[ \big( \mathcal A^M_{z^+} \big)^{-1} \mathcal P^L_{z^+} \right] \right\}
\label{3D_modified_spatial_derivative_approximations_plus_Vy}
\\
& & &
\cdot
\left\{ \big(\mathcal I^N_{x^+} \otimes \mathcal I^M_{y^+} \otimes \mathcal E^L_{z^+} \big)^T V^+_y \right. 
 - \,
\left. \mathcal T^{V^-_y} \left[ \big(\mathcal I^N_{x^-} \otimes \mathcal I^M_{y^-} \otimes \mathcal E^R_{z^-} \big)^T V^-_y \right] \right\} ; \nonumber
\\ 
\mathcal D^{\Sigma^+_{zz}}_z \Sigma^+_{zz} 
& \, \longrightarrow \, &
\overline{\mathcal D^{\Sigma^+_{zz}}_z \Sigma^+_{zz}} 
\, = \, &
\mathcal D^{\Sigma^+_{zz}}_z \Sigma^+_{zz} 
\, + \, \tfrac{1}{2}
\left\{ \mathcal I^N_{x^+} \otimes \mathcal I^N_{y^+} \otimes \left[ \big( \mathcal A^M_{z^+} \big)^{-1} \mathcal P^L_{z^+} \right] \right\}
\label{3D_modified_spatial_derivative_approximations_plus_Szz}
\\
& & &
\cdot
\left\{ \big(\mathcal I^N_{x^+} \otimes \mathcal I^N_{y^+} \otimes \mathcal E^L_{z^+} \big)^T \Sigma^+_{zz} \right. 
 - \,
\left. \mathcal T^{\Sigma^-_{zz}} \left[ \big(\mathcal I^N_{x^-} \otimes \mathcal I^N_{y^-} \otimes \mathcal E^R_{z^-} \big)^T \Sigma^-_{zz} \right] \right\} ; \nonumber
\\ 
\mathcal D^{V^+_z}_z V^+_z 
& \, \longrightarrow \, &
\overline{\mathcal D^{V^+_z}_z V^+_z} 
\, = \, &
\mathcal D^{V^+_z}_z V^+_z 
\, + \, \tfrac{1}{2}
\left\{ \mathcal I^N_{x^+} \otimes \mathcal I^N_{y^+} \otimes \left[ \big( \mathcal A^N_{z^+} \big)^{-1} \mathcal E^L_{z^+} \right] \right\}
\label{3D_modified_spatial_derivative_approximations_plus_Vz}
\\
& & &
\cdot
\left\{ \big(\mathcal I^N_{x^+} \otimes \mathcal I^N_{y^+} \otimes \mathcal P^L_{z^+} \big)^T V^+_z \right. 
 - \,
\left. \mathcal T^{V^-_z} \left[ \big(\mathcal I^N_{x^-} \otimes \mathcal I^N_{y^-} \otimes \mathcal P^R_{z^-} \big)^T V^-_z \right] \right\} . \nonumber
\end{empheq}
\end{subequations}
\noindent (On the minus side:)
\begin{subequations}
\label{3D_modified_spatial_derivative_approximations_minus}
\begin{empheq}[left=\hspace{-5em},right=]{alignat = 3}
\mathcal D^{\Sigma^-_{xz}}_z \Sigma^-_{xz} 
& \, \longrightarrow \, &
\overline{\mathcal D^{\Sigma^-_{xz}}_z \Sigma^-_{xz}} 
\, = \, &
\mathcal D^{\Sigma^-_{xz}}_z \Sigma^-_{xz} 
\, - \, \tfrac{1}{2}
\left\{ \mathcal I^M_{x^-} \otimes \mathcal I^N_{y^-} \otimes \left[ \big( \mathcal A^N_{z^-} \big)^{-1} \mathcal E^R_{z^-} \right] \right\}
\label{3D_modified_spatial_derivative_approximations_minus_Sxz}
\\ 
& & &
\cdot
\left\{ \big(\mathcal I^M_{x^-} \otimes \mathcal I^N_{y^-} \otimes \mathcal P^R_{z^-} \big)^T \Sigma^-_{xz} \right. 
 - \, 
\left.\mathcal T^{\Sigma^+_{xz}} \left[ \big(\mathcal I^M_{x^+} \otimes \mathcal I^N_{y^+} \otimes \mathcal P^L_{z^+} \big)^T \Sigma^+_{xz} \right] \right\} ; \nonumber
\\
\mathcal D^{V^-_x}_z V^-_x 
& \, \longrightarrow \, &
\overline{\mathcal D^{V^-_x}_z V^-_x} 
\, = \, &
\mathcal D^{V^-_x}_z V^-_x
\, - \, \tfrac{1}{2}
\left\{ \mathcal I^M_{x^-} \otimes \mathcal I^N_{y^-} \otimes \left[ \big( \mathcal A^M_{z^-} \big)^{-1} \mathcal P^R_{z^-} \right] \right\}
\label{3D_modified_spatial_derivative_approximations_minus_Vx}
\\ 
& & &
\cdot
\left\{ \big(\mathcal I^M_{x^-} \otimes \mathcal I^N_{y^-} \otimes \mathcal E^R_{z^-} \big)^T V^-_x \right. 
 - \, 
\left.\mathcal T^{V^+_x} \left[ \big(\mathcal I^M_{x^+} \otimes \mathcal I^N_{y^+} \otimes \mathcal E^L_{z^+} \big)^T V^+_x \right] \right\} ; \nonumber
\\
\mathcal D^{\Sigma^-_{yz}}_z \Sigma^-_{yz} 
& \, \longrightarrow \, &
\overline{\mathcal D^{\Sigma^-_{yz}}_z \Sigma^-_{yz}} 
\, = \, &
\mathcal D^{\Sigma^-_{yz}}_z \Sigma^-_{yz} 
\, - \, \tfrac{1}{2}
\left\{ \mathcal I^N_{x^-} \otimes \mathcal I^M_{y^-} \otimes \left[ \big( \mathcal A^N_{z^-} \big)^{-1} \mathcal E^R_{z^-} \right] \right\}
\label{3D_modified_spatial_derivative_approximations_minus_Syz}
\\
& & &
\cdot
\left\{ \big(\mathcal I^N_{x^-} \otimes \mathcal I^M_{y^-} \otimes \mathcal P^R_{z^-} \big)^T \Sigma^-_{yz} \right. 
 - \,
\left. \mathcal T^{\Sigma^+_{yz}} \left[ \big(\mathcal I^N_{x^+} \otimes \mathcal I^M_{y^+} \otimes \mathcal P^L_{z^+} \big)^T \Sigma^+_{yz} \right] \right\} ; \nonumber
\\
\mathcal D^{V^-_y}_z V^-_y 
& \, \longrightarrow \, &
\overline{\mathcal D^{V^-_y}_z V^-_y} 
\, = \, &
\mathcal D^{V^-_y}_z V^-_y 
\, - \, \tfrac{1}{2}
\left\{ \mathcal I^N_{x^-} \otimes \mathcal I^M_{y^-} \otimes \left[ \big( \mathcal A^M_{z^-} \big)^{-1} \mathcal P^R_{z^-} \right] \right\}
\label{3D_modified_spatial_derivative_approximations_minus_Vy}
\\
& & &
\cdot
\left\{ \big(\mathcal I^N_{x^-} \otimes \mathcal I^M_{y^-} \otimes \mathcal E^R_{z^-} \big)^T V^-_y \right. 
 - \,
\left. \mathcal T^{V^+_y} \left[ \big(\mathcal I^N_{x^+} \otimes \mathcal I^M_{y^+} \otimes \mathcal E^L_{z^+} \big)^T V^+_y \right] \right\} ; \nonumber
\\ 
\mathcal D^{\Sigma^-_{zz}}_z \Sigma^-_{zz} 
& \, \longrightarrow \, &
\overline{\mathcal D^{\Sigma^-_{zz}}_z \Sigma^-_{zz}} 
\, = \, &
\mathcal D^{\Sigma^-_{zz}}_z \Sigma^-_{zz} 
\, - \, \tfrac{1}{2}
\left\{ \mathcal I^N_{x^-} \otimes \mathcal I^N_{y^-} \otimes \left[ \big( \mathcal A^M_{z^-} \big)^{-1} \mathcal P^R_{z^-} \right] \right\}
\label{3D_modified_spatial_derivative_approximations_minus_Szz}
\\
& & &
\cdot
\left\{ \big(\mathcal I^N_{x^-} \otimes \mathcal I^N_{y^-} \otimes \mathcal E^R_{z^-} \big)^T \Sigma^-_{zz} \right. 
 - \,
\left. \mathcal T^{\Sigma^+_{zz}} \left[ \big(\mathcal I^N_{x^+} \otimes \mathcal I^N_{y^+} \otimes \mathcal E^L_{z^+} \big)^T \Sigma^+_{zz} \right] \right\} ; \nonumber
\\ 
\mathcal D^{V^-_z}_z V^-_z 
& \, \longrightarrow \, &
\overline{\mathcal D^{V^-_z}_z V^-_z} 
\, = \, &
\mathcal D^{V^-_z}_z V^-_z 
\, - \, \tfrac{1}{2}
\left\{ \mathcal I^N_{x^-} \otimes \mathcal I^N_{y^-} \otimes \left[ \big( \mathcal A^N_{z^-} \big)^{-1} \mathcal E^R_{z^-} \right] \right\}
\label{3D_modified_spatial_derivative_approximations_minus_Vz}
\\
& & &
\cdot
\left\{ \big(\mathcal I^N_{x^-} \otimes \mathcal I^N_{y^-} \otimes \mathcal P^R_{z^-} \big)^T V^-_z \right. 
 - \,
\left. \mathcal T^{V^+_z} \left[ \big(\mathcal I^N_{x^+} \otimes \mathcal I^N_{y^+} \otimes \mathcal P^L_{z^+} \big)^T V^+_z \right] \right\} . \nonumber
\end{empheq}
\end{subequations}

In equations \eqref{3D_modified_spatial_derivative_approximations_plus} and \eqref{3D_modified_spatial_derivative_approximations_minus}, 
the continuity in $\sigma_{xz}$, $v_x$, $\sigma_{yz}$, $v_y$, $\sigma_{zz}$, and $v_z$ are accounted for by 
equations \eqref{3D_modified_spatial_derivative_approximations_plus_Sxz} - \eqref{3D_modified_spatial_derivative_approximations_plus_Vz} on the plus side, respectively, and by equations \eqref{3D_modified_spatial_derivative_approximations_minus_Sxz} - \eqref{3D_modified_spatial_derivative_approximations_minus_Vz} on the minus side, respectively.
Moreover, we used $\mathcal T$ to denote a 2D interpolation operator acting on the interface. 
For example, $\mathcal T^{V^+_z}$ in equation \eqref{3D_modified_spatial_derivative_approximations_minus_Vz} operates on the projection of $V^+_z$ on the interface, i.e., $\big(\mathcal I^N_{x^+} \otimes \mathcal I^N_{y^+} \otimes \mathcal P^L_{z^+} \big)^T V^+_z$, providing the approximations to $v_z$ at locations matching those of the projection of $V^-_z$ on the interface, i.e., $\big(\mathcal I^N_{x^-} \otimes \mathcal I^N_{y^-} \otimes \mathcal P^R_{z^-} \big)^T V^-_z$.

The 2D interpolation operators appearing in equations \eqref{3D_modified_spatial_derivative_approximations_plus} and \eqref{3D_modified_spatial_derivative_approximations_minus} are constructed as tensor products of 1D interpolation operators as follows:
\begin{equation}
\label{2D_interpolation_operators_construction}
\arraycolsep=1em\def\arraystretch{0.75}
\begin{array}{ll}
\mathcal T^{\Sigma^+_{xz}} \, = \, \mathcal T^{V^+_x} \, = \, \mathcal T^{M}_{x^+} \otimes \mathcal T^{N}_{y^+} \, , 
&
\mathcal T^{\Sigma^-_{xz}} \, = \, \mathcal T^{V^-_x} \, = \, \mathcal T^{M}_{x^-} \otimes \mathcal T^{N}_{y^-} \, ;  
\\
\mathcal T^{\Sigma^+_{yz}} \, = \, \mathcal T^{V^+_y} \, = \, \mathcal T^{N}_{x^+} \otimes \mathcal T^{M}_{y^+} \, ,
&
\mathcal T^{\Sigma^-_{yz}} \, = \, \mathcal T^{V^-_y} \, = \, \mathcal T^{N}_{x^-} \otimes \mathcal T^{M}_{y^-} \, ;
\\
\mathcal T^{\Sigma^+_{zz}} \, = \, \mathcal T^{V^+_z} \, = \, \mathcal T^{N}_{x^+} \otimes \mathcal T^{N}_{y^+} \, ,
&
\mathcal T^{\Sigma^-_{zz}} \, = \, \mathcal T^{V^-_z} \, = \, \mathcal T^{N}_{x^-} \otimes \mathcal T^{N}_{y^-} \ .
\end{array}
\end{equation}
To properly handle the nonconforming interface, the 1D building blocks in equation \eqref{2D_interpolation_operators_construction} are required to satisfy the following properties:
\begin{equation}
\label{1D_interpolation_operators_properties}
\arraycolsep=1em\def\arraystretch{1}
\begin{array}{ll}
\mathcal A^{N}_{x^+} \mathcal T^{N}_{x^-} \, = \, \big( \mathcal A^{N}_{x^-} \mathcal T^{N}_{x^+} \big)^T \, ,
&
\mathcal A^{M}_{x^+} \mathcal T^{M}_{x^-} \, = \, \big( \mathcal A^{M}_{x^-} \mathcal T^{M}_{x^+} \big)^T \, ;
\\
\mathcal A^{N}_{y^+} \mathcal T^{N}_{y^-} \, = \, \big( \mathcal A^{N}_{y^-} \mathcal T^{N}_{y^+} \big)^T \, ,
&
\mathcal A^{M}_{y^+} \mathcal T^{M}_{y^-} \, = \, \big( \mathcal A^{M}_{y^-} \mathcal T^{M}_{y^+} \big)^T \, .
\end{array}
\end{equation}

A collection of 1D interpolation operators satisfying the properties in equation \eqref{1D_interpolation_operators_properties} are presented in \ref{appendix_interpolation_operators}, corresponding to a variety of grid spacing ratios.
Such ratios are allowed to be rational numbers, instead of just integers, which improves the application flexibility of the proposed discretization setting. 
Similar to the construction of SBP operators, these interpolation operators are also constructed with the assistance of symbolic computing software.

If the discrete energy analysis is carried out without regard to the interface conditions,  
it can be shown that the remaining terms associated with the interface are: 
\begin{equation}
\label{3D_discrete_energy_time_derivative_interface}
\arraycolsep=0.5pt\def\arraystretch{1}
\begin{array}{l}
\left.
\arraycolsep=0.5pt\def\arraystretch{1}
\begin{array}{rrrllcrr}
\displaystyle \frac{d E}{d t}
& \ = & \ - \,\, & \displaystyle 
\Big[ \big(\Sigma^+_{xz}\big)^T \big(\mathcal I^M_{x^+} \otimes \mathcal I^N_{y^+} \otimes \mathcal P^L_{z^+} \big) \Big] \,
& \cdot & \,
\Big[\mathcal A^M_{x^+} \otimes \mathcal A^N_{y^+} \Big] 
& \cdot & 
\Big[ \big(\mathcal I^M_{x^+} \otimes \mathcal I^N_{y^+} \otimes \mathcal E^L_{z^+} \big)^T V^+_x \Big] \phantom{.}
\\
& & \ - \,\, & \displaystyle 
\Big[ \big(\Sigma^+_{yz}\big)^T \big(\mathcal I^N_{x^+} \otimes \mathcal I^M_{y^+} \otimes \mathcal P^L_{z^+} \big) \Big] 
& \cdot & 
\Big[\mathcal A^N_{x^+} \otimes \mathcal A^M_{y^+} \Big] 
& \cdot & 
\Big[ \big(\mathcal I^N_{x^+} \otimes \mathcal I^M_{y^+} \otimes \mathcal E^L_{z^+} \big)^T V^+_y \Big] \phantom{.}
\\
& & \ - \,\, & \displaystyle 
\Big[ \big(\Sigma^+_{zz}\big)^T \big(\mathcal I^N_{x^+} \otimes \mathcal I^N_{y^+} \otimes \mathcal E^L_{z^+} \big) \Big] 
& \cdot & 
\Big[\mathcal A^N_{x^+} \otimes \mathcal A^N_{y^+} \Big] \,
& \cdot &
\Big[ \big(\mathcal I^N_{x^+} \otimes \mathcal I^N_{y^+} \otimes \mathcal P^L_{z^+} \big)^T V^+_z \Big] \phantom{.}
\end{array}
\right\} \text{($+$ side)}
\\[5.ex]
\left.
\arraycolsep=0.5pt\def\arraystretch{1}
\begin{array}{rrrllcrr}
\displaystyle \phantom{\frac{d E}{d t}}
& \ \phantom{=} & \ + \,\, & \displaystyle 
\Big[ \big(\Sigma^-_{xz}\big)^T \big(\mathcal I^M_{x^-} \otimes \mathcal I^N_{y^-} \otimes \mathcal P^R_{z^-} \big) \Big] 
& \cdot & 
\Big[\mathcal A^M_{x^-} \otimes \mathcal A^N_{y^-} \Big] 
& \cdot & 
\Big[ \big(\mathcal I^M_{x^-} \otimes \mathcal I^N_{y^-} \otimes \mathcal E^L_{z^-} \big)^T V^-_x \Big] \phantom{.}\hspace{0.2em}
\\
& & \ + \,\, & \displaystyle 
\Big[ \big(\Sigma^-_{yz}\big)^T \big(\mathcal I^N_{x^-} \otimes \mathcal I^M_{y^-} \otimes \mathcal P^R_{z^-} \big) \Big] 
& \cdot & 
\Big[\mathcal A^N_{x^-} \otimes \mathcal A^M_{y^-} \Big] 
& \cdot & 
\Big[ \big(\mathcal I^N_{x^-} \otimes \mathcal I^M_{y^-} \otimes \mathcal E^L_{z^-} \big)^T V^-_y \Big] \phantom{.}\hspace{0.2em}
\\
& & \ + \,\, & \displaystyle 
\Big[ \big(\Sigma^-_{zz}\big)^T \big(\mathcal I^N_{x^-} \otimes \mathcal I^N_{y^-} \otimes \mathcal E^R_{z^-} \big) \Big] 
& \cdot & 
\Big[\mathcal A^N_{x^-} \otimes \mathcal A^N_{y^-} \Big] 
& \cdot & 
\Big[ \big(\mathcal I^N_{x^-} \otimes \mathcal I^N_{y^-} \otimes \mathcal P^L_{z^-} \big)^T V^-_z \Big] .\hspace{0.2em}
\end{array}
\right\} \text{($-$ side)}
\end{array}
\end{equation} 
Instead, carrying out the discrete energy analysis with the modifications presented in equations \eqref{3D_modified_spatial_derivative_approximations_plus} - \eqref{3D_modified_spatial_derivative_approximations_minus} and the properties in equation \eqref{1D_interpolation_operators_properties}, it can be shown that the remaining terms in equation \eqref{3D_discrete_energy_time_derivative_interface} are cancelled out without introducing any extra term, i.e., the energy conserving property is preserved across the interface.

To illustrate, let us consider the pairs $\left\{\Sigma_{zz}^+; V_z^+\right\}$ and $\left\{\Sigma_{zz}^-; V_z^-\right\}$ as an example. 
The remaining terms in equation \eqref{3D_discrete_energy_time_derivative_interface} associated with these pairs are at the third and the sixth lines, respectively.
Now, carrying out the discrete energy analysis with the modified derivative approximations presented in equations \eqref{3D_modified_spatial_derivative_approximations_plus_Szz}, \eqref{3D_modified_spatial_derivative_approximations_plus_Vz},
\eqref{3D_modified_spatial_derivative_approximations_minus_Szz}, and 
\eqref{3D_modified_spatial_derivative_approximations_minus_Vz}, 
the remaining terms in $\tfrac{d E}{d t}$ associated with these pairs become

\begin{linenomath}
\begin{subequations}
\label{3D_energy_conserving_illustration_Szz_Vz}
\begin{empheq}[left=\hspace{-2em},right=]{alignat = 1}
& 
- \
\uline{ 
\left[ \big(\Sigma^+_{zz}\big)^T \big(\mathcal I^N_{x^+} \otimes \mathcal I^N_{y^+} \otimes \mathcal E^L_{z^+} \big) \right] 
\cdot \left[\mathcal A^N_{x^+} \otimes \mathcal A^N_{y^+} \right] \cdot 
\left[ \big(\mathcal I^N_{x^+} \otimes \mathcal I^N_{y^+} \otimes \mathcal P^L_{z^+} \big)^T V^+_z \right] 
}
\label{3D_energy_conserving_illustration_Szz_Vz_a}
\\ 
& + \ 
\uwave{
\left[ \big(\Sigma^-_{zz}\big)^T \big(\mathcal I^N_{x^-} \otimes \mathcal I^N_{y^-} \otimes \mathcal E^R_{z^-} \big) \right] 
\cdot \left[\mathcal A^N_{x^-} \otimes \mathcal A^N_{y^-} \right] \cdot 
\left[ \big(\mathcal I^N_{x^-} \otimes \mathcal I^N_{y^-} \otimes \mathcal P^L_{z^-} \big)^T V^-_z \right] 
}
\label{3D_energy_conserving_illustration_Szz_Vz_b}
\\
& + \ 
\tfrac{1}{2} \,
\uline{ 
\left[ \big( V^+_z \big)^T \big(\mathcal I^N_{x^+} \otimes \mathcal I^N_{y^+} \otimes \mathcal P^L_{z^+} \big) \right] 
\cdot \left[\mathcal A^N_{x^+} \otimes \mathcal A^N_{y^+} \right] \cdot 
\left[ \big(\mathcal I^N_{x^+} \otimes \mathcal I^N_{y^+} \otimes \mathcal E^L_{z^+} \big)^T \Sigma^+_{zz} \right] 
}
\label{3D_energy_conserving_illustration_Szz_Vz_c}
\\
& - \ 
\tfrac{1}{2} \, 
\dashuline{
\left[ \big( V^+_z \big)^T \big(\mathcal I^N_{x^+} \otimes \mathcal I^N_{y^+} \otimes \mathcal P^L_{z^+} \big) \right] 
\cdot \left[ \big( \mathcal A^N_{x^+} \mathcal T^N_{x^-} \big) \otimes \big( \mathcal A^N_{y^+} \mathcal T^N_{y^-} \big) \right] \cdot 
\left[ \big(\mathcal I^N_{x^-} \otimes \mathcal I^N_{y^-} \otimes \mathcal E^L_{z^-} \big)^T \Sigma^-_{zz} \right] 
}
\label{3D_energy_conserving_illustration_Szz_Vz_d}
\\
& + \ 
\tfrac{1}{2} \, 
\uline{
\left[ \big(\Sigma^+_{zz}\big)^T \big(\mathcal I^N_{x^+} \otimes \mathcal I^N_{y^+} \otimes \mathcal E^L_{z^+} \big) \right] 
\cdot \left[\mathcal A^N_{x^+} \otimes \mathcal A^N_{y^+} \right] \cdot 
\left[ \big(\mathcal I^N_{x^+} \otimes \mathcal I^N_{y^+} \otimes \mathcal P^L_{z^+} \big)^T V^+_z \right] 
}
\label{3D_energy_conserving_illustration_Szz_Vz_e}
\\
& - \ 
\tfrac{1}{2} \, 
\dotuline{
\left[ \big(\Sigma^+_{zz}\big)^T \big(\mathcal I^N_{x^+} \otimes \mathcal I^N_{y^+} \otimes \mathcal E^L_{z^+} \big) \right] 
\cdot \left[ \big( \mathcal A^N_{x^+} \mathcal T^N_{x^-} \big) \otimes \big( \mathcal A^N_{y^+} \mathcal T^N_{y^-} \big) \right] \cdot
\left[ \big(\mathcal I^N_{x^-} \otimes \mathcal I^N_{y^-} \otimes \mathcal P^L_{z^-} \big)^T V^-_z \right] 
}
\label{3D_energy_conserving_illustration_Szz_Vz_f}
\\
& - \ 
\tfrac{1}{2} \,
\uwave{ 
\left[ \big( V^-_z \big)^T \big(\mathcal I^N_{x^-} \otimes \mathcal I^N_{y^-} \otimes \mathcal P^L_{z^-} \big) \right] 
\cdot \left[\mathcal A^N_{x^-} \otimes \mathcal A^N_{y^-} \right] \cdot 
\left[ \big(\mathcal I^N_{x^-} \otimes \mathcal I^N_{y^-} \otimes \mathcal E^L_{z^-} \big)^T \Sigma^-_{zz} \right] 
}
\label{3D_energy_conserving_illustration_Szz_Vz_g}
\\
& + \ 
\tfrac{1}{2} \, 
\dotuline{
\left[ \big( V^-_z \big)^T \big(\mathcal I^N_{x^-} \otimes \mathcal I^N_{y^-} \otimes \mathcal P^L_{z^-} \big) \right] 
\cdot \left[ \big( \mathcal A^N_{x^-} \mathcal T^N_{x^+} \big) \otimes \big( \mathcal A^N_{y^-} \mathcal T^N_{y^+} \big) \right] \cdot 
\left[ \big(\mathcal I^N_{x^+} \otimes \mathcal I^N_{y^+} \otimes \mathcal E^L_{z^+} \big)^T \Sigma^+_{zz} \right] 
}
\label{3D_energy_conserving_illustration_Szz_Vz_h}
\\
& - \ 
\tfrac{1}{2} \, 
\uwave{
\left[ \big(\Sigma^-_{zz}\big)^T \big(\mathcal I^N_{x^-} \otimes \mathcal I^N_{y^-} \otimes \mathcal E^L_{z^-} \big) \right] 
\cdot \left[\mathcal A^N_{x^-} \otimes \mathcal A^N_{y^-} \right] \cdot 
\left[ \big(\mathcal I^N_{x^-} \otimes \mathcal I^N_{y^-} \otimes \mathcal P^L_{z^-} \big)^T V^-_z \right] 
}
\label{3D_energy_conserving_illustration_Szz_Vz_i}
\\
& + \ 
\tfrac{1}{2} \, 
\dashuline{
\left[ \big(\Sigma^-_{zz}\big)^T \big(\mathcal I^N_{x^-} \otimes \mathcal I^N_{y^-} \otimes \mathcal E^L_{z^-} \big) \right] 
\cdot \left[ \big( \mathcal A^N_{x^-} \mathcal T^N_{x^+} \big) \otimes \big( \mathcal A^N_{y^-} \mathcal T^N_{y^+} \big) \right] \cdot
\left[ \big(\mathcal I^N_{x^+} \otimes \mathcal I^N_{y^+} \otimes \mathcal P^L_{z^+} \big)^T V^+_z \right] 
},
\label{3D_energy_conserving_illustration_Szz_Vz_j}
\end{empheq}
\end{subequations}
\end{linenomath}
where equations \eqref{3D_energy_conserving_illustration_Szz_Vz_a} and \eqref{3D_energy_conserving_illustration_Szz_Vz_b} are copied from the third and the sixth lines of equation \eqref{3D_discrete_energy_time_derivative_interface}, respectively; 
equations \eqref{3D_energy_conserving_illustration_Szz_Vz_c} and \eqref{3D_energy_conserving_illustration_Szz_Vz_d} stem from the modification in equation \eqref{3D_modified_spatial_derivative_approximations_plus_Szz}; 
equations \eqref{3D_energy_conserving_illustration_Szz_Vz_e} and \eqref{3D_energy_conserving_illustration_Szz_Vz_f} stem from the modification in equation \eqref{3D_modified_spatial_derivative_approximations_plus_Vz}; 
equations \eqref{3D_energy_conserving_illustration_Szz_Vz_g} and \eqref{3D_energy_conserving_illustration_Szz_Vz_h} stem from the modification in equation \eqref{3D_modified_spatial_derivative_approximations_minus_Szz}; 
equations \eqref{3D_energy_conserving_illustration_Szz_Vz_i} and \eqref{3D_energy_conserving_illustration_Szz_Vz_j} stem from the modification in equation \eqref{3D_modified_spatial_derivative_approximations_minus_Vz}.
It can be verified easily that in equation \eqref{3D_energy_conserving_illustration_Szz_Vz}, terms underlined with lines of the same style cancel out.
The properties in equation \eqref{1D_interpolation_operators_properties} are needed for the terms underlined with the dashed and dotted lines to cancel out, respectively.

To sum up, all remaining terms in equation \eqref{3D_discrete_energy_time_derivative_interface} associated with the pairs $\left\{\Sigma_{zz}^+; V_z^+\right\}$ and $\left\{\Sigma_{zz}^-; V_z^-\right\}$ are cancelled out without introducing new terms because of the modifications made in equations \eqref{3D_modified_spatial_derivative_approximations_plus_Szz}, \eqref{3D_modified_spatial_derivative_approximations_plus_Vz},
\eqref{3D_modified_spatial_derivative_approximations_minus_Szz}, and \eqref{3D_modified_spatial_derivative_approximations_minus_Vz}.
Similar observations can be made for the pairs $\left\{\Sigma_{xz}^+; V_x^+\right\}$ and $\left\{\Sigma_{xz}^-; V_x^-\right\}$, which correspond to the first and fourth lines in equation \eqref{3D_discrete_energy_time_derivative_interface}, respectively, 
and for the pairs $\left\{\Sigma_{yz}^+; V_y^+\right\}$ and $\left\{\Sigma_{yz}^-; V_y^-\right\}$, which correspond to the second and fifth lines in equation \eqref{3D_discrete_energy_time_derivative_interface}, respectively.

In practical implementation, the following procedure can be used to arrive at the modified derivative approximations presented in equations \eqref{3D_modified_spatial_derivative_approximations_plus} and \eqref{3D_modified_spatial_derivative_approximations_minus}, illustrated using equation \eqref{3D_modified_spatial_derivative_approximations_minus_Vz} as a particular example, which concerns the approximation of $\frac{\partial v_z}{\partial z}$ on the minus side of the interface.
\begin{itemize}
\item[1)] Apply $\mathcal D^{V^-_z}_z$ to $V^-_z$ and update the derivative approximation; essentially, this boils down to applying the 1D difference operator $\mathcal D^{M}_{z^-}$ to each $z$-direction grid column on the minus side of the interface (see equation \eqref{SBP_difference_operators_3D} for the composition of $\mathcal D^{V^-_z}_z$).
\item[2)] Project $V^-_z$ to the interface from the minus side; essentially, this boils down to applying the 1D projection operator $\mathcal P^R_{z^-}$ to each $z$-direction grid column on the minus side of the interface.
\item[3)] Project $V^+_z$ to the interface from the plus side; essentially, this boils down to applying the 1D projection operator $\mathcal P^L_{z^+}$ to each $z$-direction grid column on the plus side of the interface.
\item[4)] Interpolate the projected values from step 3) to match the positions of the projected values from step 2); this step requires the 2D interpolation operator $\mathcal T^{V^+_z}$ presented in equation \eqref{2D_interpolation_operators_construction}.
\item[5)] Penalize the differences of the values obtained from steps 2) and 4) and update the derivative approximations; essentially, this boils down to applying $\big( \mathcal A^N_{z^-} \big)^{-1} \mathcal E^R_{z^-}$ on each difference.
\end{itemize}
We note here that order of the above steps does not need to be followed strictly in actual implementation. For example, step 1 and steps 2-5 can be carried out independently. Their relative order can be arranged flexibly to improve the computation and communication sequence in parallel computing environment.

\section{Numerical Examples}\label{numerical_examples}
In this section, two numerical examples are presented to corroborate the above analytical development, one concerning a layer-wise homogeneous model and the other concerning the overthrust model; see \cite{aminzadeh19973}.

To demonstrate the energy conserving property, periodic boundary conditions are considered for the two horizontal directions (i.e., the $x$- and $y$-directions) while free surface boundary conditions are considered for the top and bottom boundaries since there is no energy loss associated with these boundary conditions. 
The periodic boundary conditions are addressed using the standard practice of wrapping the stencils around the boundaries.
The staggered leapfrog scheme is used for temporal discretization.

The numerical experiments are conducted on a Cray XC-40 distributed-memory system, with each node containing dual 16-core Intel Haswell processors running at 2.3 GHz and 128 GB of DDR4 memory, connected via an Aries network.

\subsection{Layered model} \label{layered_model}
In this example, we consider a model consisting of four homogeneous layers.  
The relevant physical and numerical parameters associated with these layers (from top to bottom) are displayed in Table~\ref{layered_model_parameters}.
\begin{table}[H]
\captionsetup{width=0.9\textwidth, font=small,labelfont=small}
\caption{Physical and numerical parameters associated with the layered model.}
\label{layered_model_parameters}
\centering
\begin{tabular}{| c | c | c | c | c | c |}
\hline
             & $c_s$ & $c_p$ & $\rho$ & $\Delta x$ & $N_\text{ppw}$ \\ \hline
Layer 1 &  300 m/s   &  800 m/s  & 1600 kg/m$^3$ & 1 m & 12 \\ \hline
Layer 2 &  600 m/s   & 1800 m/s & 2100 kg/m$^3$ & 2 m & 12 \\ \hline
Layer 3 &  900 m/s   & 2400 m/s & 2300 kg/m$^3$ & 3 m & 12 \\ \hline 
Layer 4 & 2700 m/s  & 5000 m/s & 2500 kg/m$^3$ & 9 m & 12 \\ \hline
\end{tabular}
\end{table}

The source profile is specified as the Ricker wavelet with central frequency at 10 Hz. 
The maximum frequency in the wave content is assumed to be 25 Hz, which is used to calculate the $N_\text{ppw}$ (i.e., number of discretization grid points per minimum wavelength) shown in Table~\ref{layered_model_parameters}.
The grid spacing ratios at the layer interfaces are 1:2, 2:3, and 1:3, respectively, from top to bottom.
The temporal step size $\Delta t$ is specified as $\sim$$5.132 \times 10^{-4}$ s, which corresponds to equation \eqref{appendix_Delta_t_i} with $\mathbb{C}_{\text{CFL}} = 0.8$. 
(We note here that the maximum constant $\mathbb{C}_{\text{CFL}}$ associated with the interior stencil $\nicefrac{[\nicefrac{1}{24},\ \nicefrac{-9}{8},\ \nicefrac{9}{8},\ \nicefrac{-1}{24}]}{\Delta x}$ is $\nicefrac{6}{7} \approx 0.857$; see~\cite{levander1988fourth}.)

For comparison, the uniform grid simulation with $\Delta x = 1$ for all layers is performed.
The temporal step size used in the uniform grid simulation is $\sim$$9.2376 \times 10^{-5}$ s, which corresponds to equation \eqref{appendix_Delta_t} with $\mathbb{C}_{\text{CFL}} = 0.8$. 
Between the two simulations, the ratio of spatial discretization points is $\sim$$3.438$; the ratio of temporal discretization points is $\sim$$5.556$; the ratio of total spatial-temporal discretization points is $\sim$$19.1$.

The size of the simulation domain is 1080 m $\times$ 1080 m $\times$ 1080 m, with the depth for each layer being 270 m.
A compressional source is applied on the normal stress components. The $V_z$ component of the velocity is recorded for comparison.
The source and receiver are placed at 540 m apart from each other on both the $x$- and $y$-directions, and at 5 m (i.e., 5 grid points) and 4.5 m (i.e., 4.5 grid points) below the top boundary, respectively.

\begin{figure}[H]
\captionsetup{width=0.9\textwidth, font=small,labelfont=small}
\centering\includegraphics[scale=0.1]{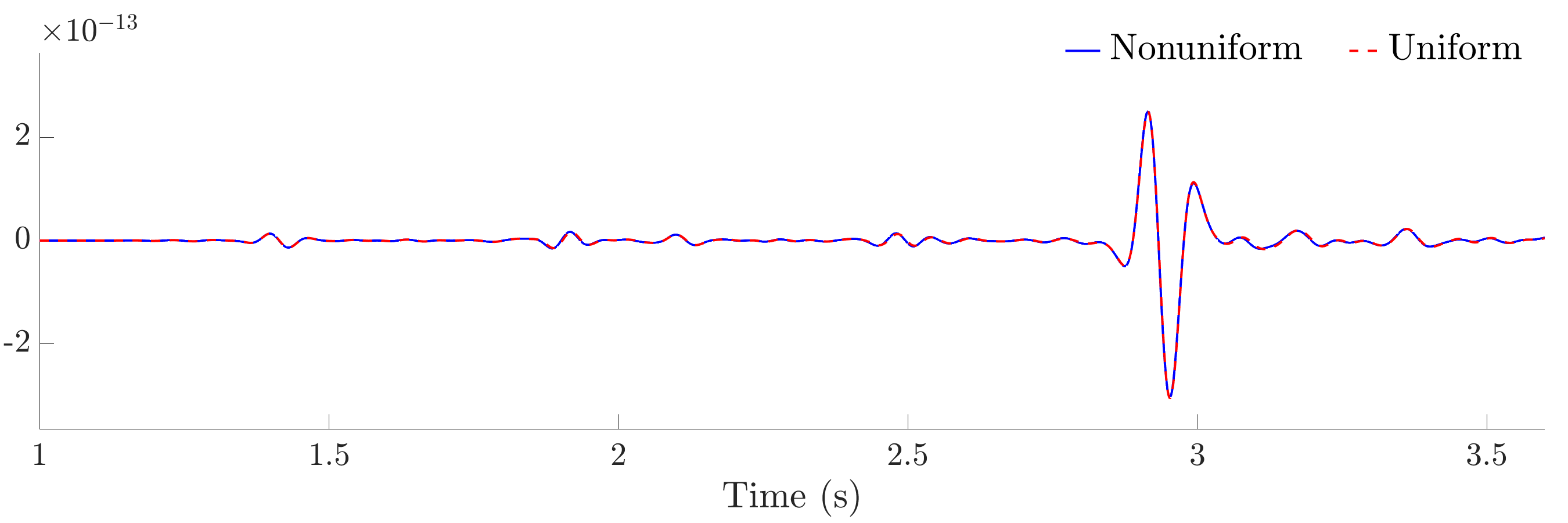}
\caption{ 
Layered model. Recorded seismograms from nonuniform and uniform grid simulations.
}
\label{fig_layered_seismogram}
\end{figure}

\begin{figure}[H]
\captionsetup{width=0.9\textwidth, font=small,labelfont=small}
\centering\includegraphics[scale=0.1]{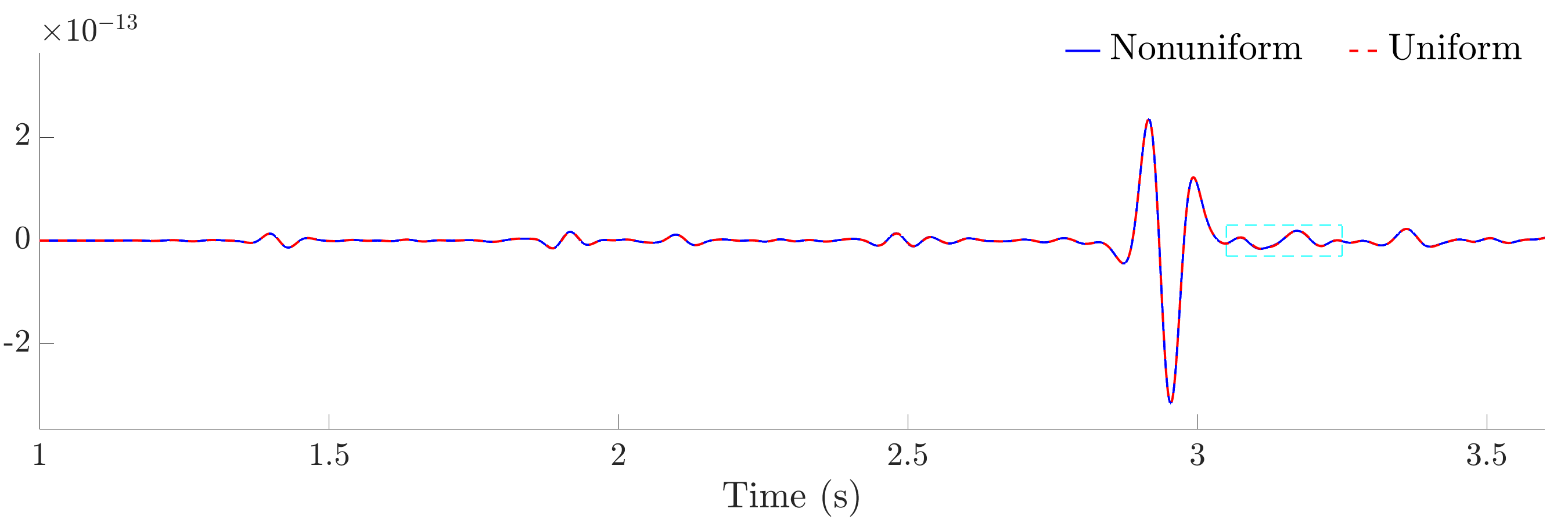}
\caption{ 
Layered model. Recorded seismograms from nonuniform and uniform grid simulations using grid spacings three times finer than those used for Figure~\ref{fig_layered_seismogram}.
Slight improvement can be noticed inside the cyan box.
}
\label{fig_layered_seismogram_3x}
\end{figure}

The recorded seismograms are displayed in Figure~\ref{fig_layered_seismogram}, where we observe satisfactory agreement between these two simulations. 
(We note here that the first second is omitted in Figure~\ref{fig_layered_seismogram} since solutions at the receiver location are still at rest.)
After projecting these seismograms onto a common set of equidistant time instances spanning the interval $[1\text{ s}, 3.5\text{ s}]$ with increment 0.001 s using cubic spline interpolation, the relative difference measured in $\ell_2$ norm is $\sim$$3.01\times 10^{-2}$.

A further verification is presented in Figure~\ref{fig_layered_seismogram_3x}, where three times finer grid spacing is used in both simulations. Slightly improved agreement can be observed in Figure~\ref{fig_layered_seismogram_3x} (e.g., inside the cyan box).
Segments of Figures~\ref{fig_layered_seismogram}~and~\ref{fig_layered_seismogram_3x} between 3.05 s and 3.25 s (corresponding to the cyan box in Figure~\ref{fig_layered_seismogram_3x}) are plotted in isolation and provided in the {\it Supplementary Materials} \ref{supp_additional_plots}, 
where the improved agreement can be observed more easily.
After applying the same interpolation as did for the seismograms in Figure 4, the relative difference measured in $\ell_2$ norm is $\sim$$8.24\times 10^{-3}$.

\begin{figure}[H]
\captionsetup{width=0.9\textwidth, font=small,labelfont=small}
\centering\includegraphics[scale=0.1]{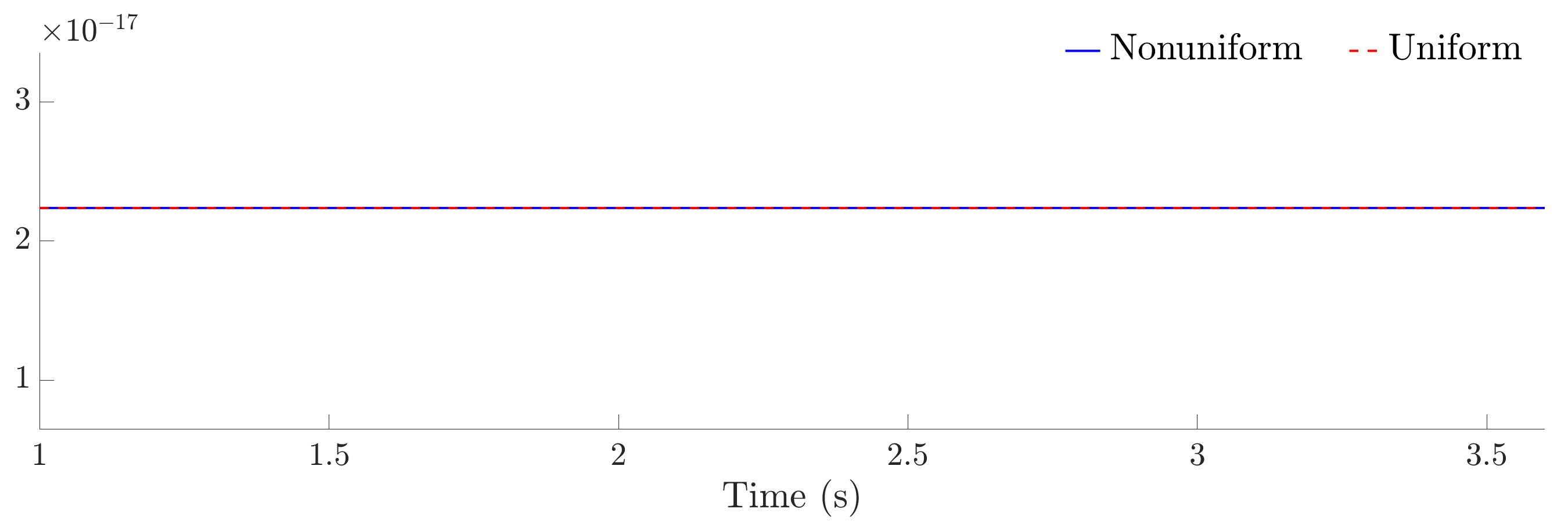}
\caption{ 
Layered model. Evolution of the discrete energy in nonuniform and uniform grid simulations.
}
\label{fig_layered_energy}
\end{figure}

Moreover, evolutions of the discrete energy in both simulations are displayed in Figure~\ref{fig_layered_energy}, which remain constant after the source term tampers off and hence corroborate the energy conserving property.
We note here that as in Figures~\ref{fig_layered_seismogram} and \ref{fig_layered_seismogram_3x}, the first second is omitted in Figure~\ref{fig_layered_energy}, which includes the period when the source terms take effect and inject energy to the simulations.

\begin{figure}[h]
\captionsetup{width=0.9\textwidth, font=small,labelfont=small}
\centering
\begin{subfigure}[b]{0.495\textwidth}
\captionsetup{width=1\textwidth, font=small,labelfont=small}
\centering\includegraphics[scale=0.125]{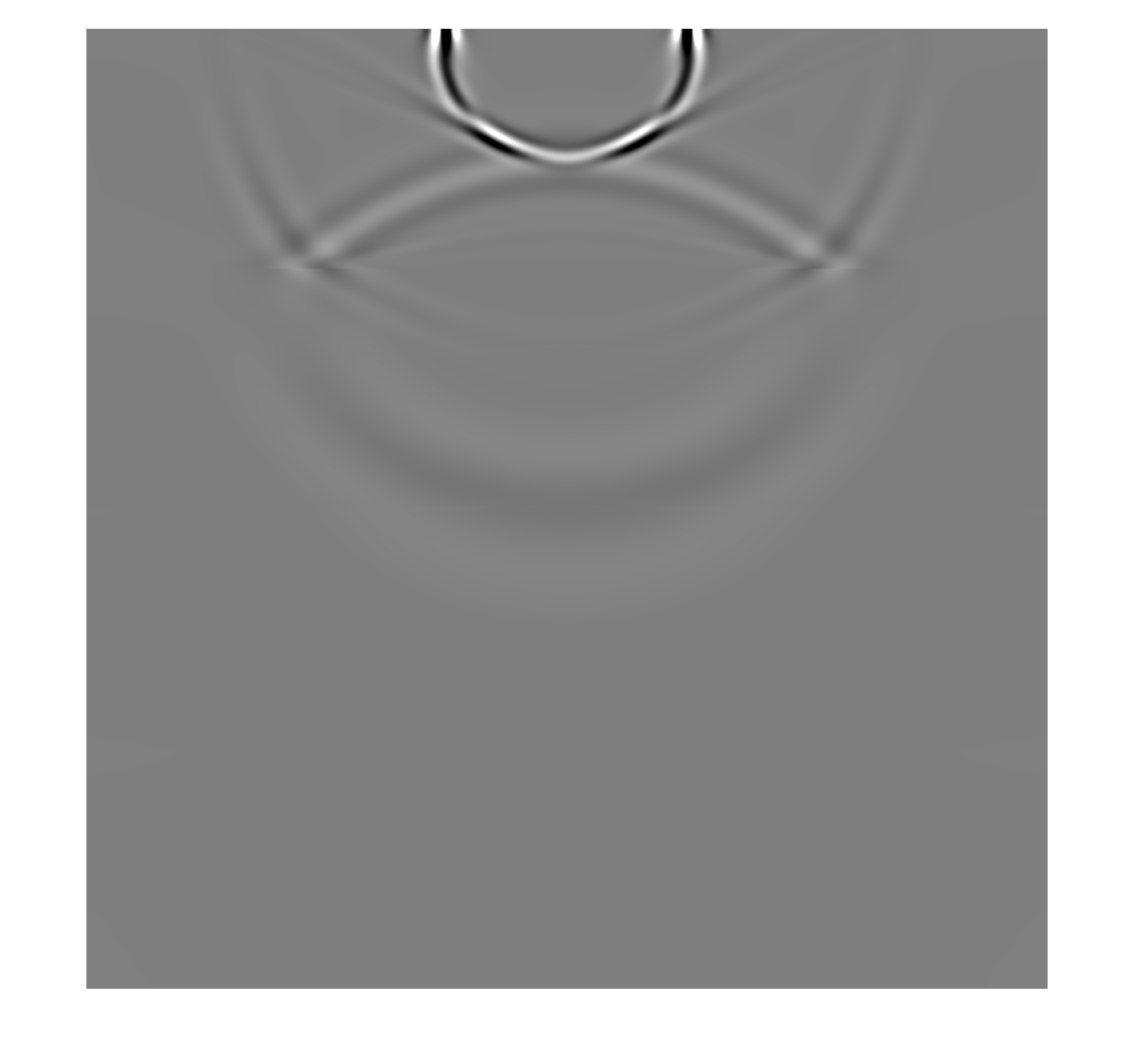}
\label{snapshot_nonuniform}
\caption{}
\end{subfigure}
\begin{subfigure}[b]{0.495\textwidth}
\captionsetup{width=1\textwidth, font=small,labelfont=small}
\centering\includegraphics[scale=0.125]{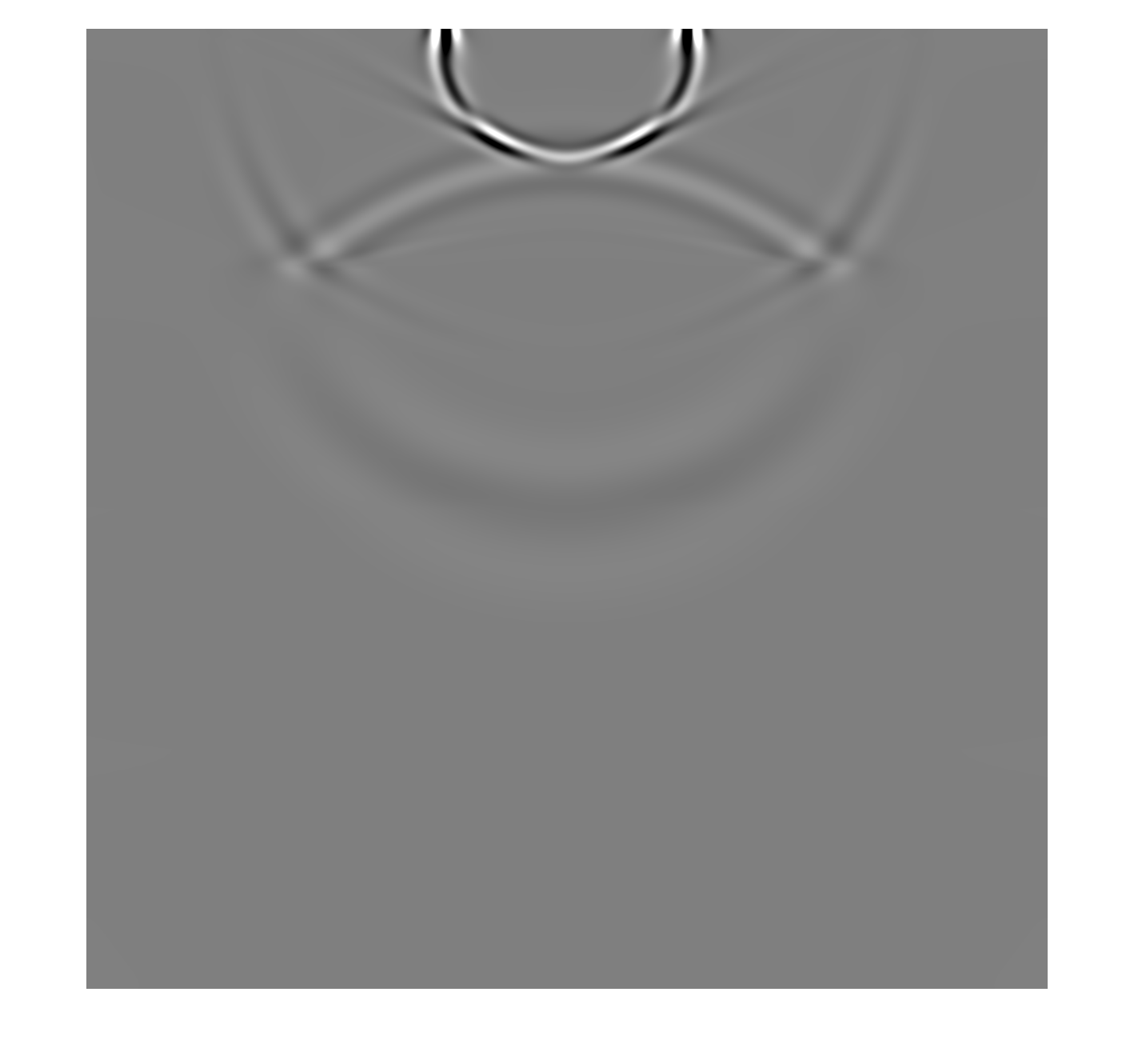}
\label{snapshot_uniform}
\caption{}
\end{subfigure} 
\caption{
Layered model. Snapshots of the vertical cross sections of the wavefield ($V_z$ component) taken at the middle of the $x$-direction and $\sim$$0.739$ s. (a): nonuniform simulation; (b): uniform simulation. These snapshots are produced with the sources placed in the center of the $xy$-plane, near the top surface.
}
\label{fig_layered_snapshot}
\end{figure}

Additionally, snapshots of the vertical cross sections of the wavefield ($V_z$ component) at the middle of the $x$-direction and parallel to the $yz$-plane are illustrated in Figure~\ref{fig_layered_snapshot} for both nonuniform (left) and uniform (right) simulations. 
We discern no obvious difference between the two snapshots.
These snapshots are produced with the sources placed in the center of the $xy$-plane, near the top surface, and taken at $\sim$$0.739$ s.

In this example, the nonuniform and uniform simulations are performed with 6785 and 23328 processor cores, respectively. 
In the nonuniform simulations, the number of processor cores assigned to each layer (top to bottom) are $18^3$, $9^3$, $6^3$, and $2^3$, respectively.
In the uniform simulations, the number of processor cores are 18, 18, and 72 along the $x$-, $y$-, and $z$-directions, respectively.

For the two simulations whose results are shown in Figure~\ref{fig_layered_seismogram}, each processor core is responsible for a subdomain with $60 \times 60 \times 15$ $M$-grid points.
The elapsed wall clock times (averaged over 5 runs) are $\sim$160 s and $\sim$872 s, respectively, translating to a ratio of $\sim$18.7 in consumed computational resources (i.e., number of processor cores multiplied by elapsed time), which agrees well with the ratio of the total number of spatial-temporal discretization points in these two simulations (19.1).
This suggests that the overhead due to the nonconforming interfaces in the nonuniform simulation is little, at $\sim$2\% on a per core basis in this case. 
The above expenditure information on computing resources are collected in Table \ref{computing_resources_figure_4}.

\begin{table}[H]
\captionsetup{width=0.9\textwidth, font=small,labelfont=small}
\caption{
Expenditure on computing resources corresponding to Figure 4. The expected ratio in total computational resources, i.e., core hours, is $\sim$$19.1$, based on the amount of spatial-temporal discretization points involved. The measured ratio is $\sim$$18.7$, which suggests a $\sim$$2\%$ overhead on a per core basis in our implementation.
}
\label{computing_resources_figure_4}
\centering
\begin{tabular}{| r | c | c | c |}
\hline
	     	  & core count & elapsed time & total core hours \\ \hline
uniform 	  & 23328   &  872 seconds  & $5.651 \times 10^3$  \\ \hline
nonuniform & 6785    &  160 seconds &  $3.016 \times 10^2$  \\ \hline
\end{tabular}
\end{table}

We note here that aside from the algorithms, simulation performance also depends on implementation quality and runtime hardware environment. 
The 
elapsed times reported here are meant only to provide supporting evidence for our algorithmic intuition that the nonuniform simulations should incur little-to-no overhead. 
However, and in general, they should not be taken as strict proof.
A few more performance related tests are included in 
the {\it Supplementary Material} \ref{supp_performance_tests}, 
which provides further support.

For the two simulations whose results are shown in Figure~\ref{fig_layered_seismogram_3x}, each processor core is responsible for a subdomain with $180 \times 180 \times 45$ $M$-grid points.
The elapsed wall clock times are $\sim$11805 s and $\sim$63482 s, respectively, 
which translates to a ratio of $\sim$18.5 in consumed computational resources, suggesting an overhead for the nonuniform simulation at $\sim$3\%.
(We note here that the two numbers on elapsed wall clock times reported in this paragraph, and similar numbers reported later for the overthrust model, are not averaged over repeated runs due to the large amount of computational resources required.)
%
A table collecting the above expenditure information on computing resources can be found in the {\it Supplementary Material} \ref{supp_additional_tables}.
It is omitted here to conserve space.

\subsection{Overthrust model} \label{overthrust_model}
\begin{figure}[H]
\captionsetup{width=0.9\textwidth, font=small,labelfont=small}
\centering\includegraphics[scale=0.125]{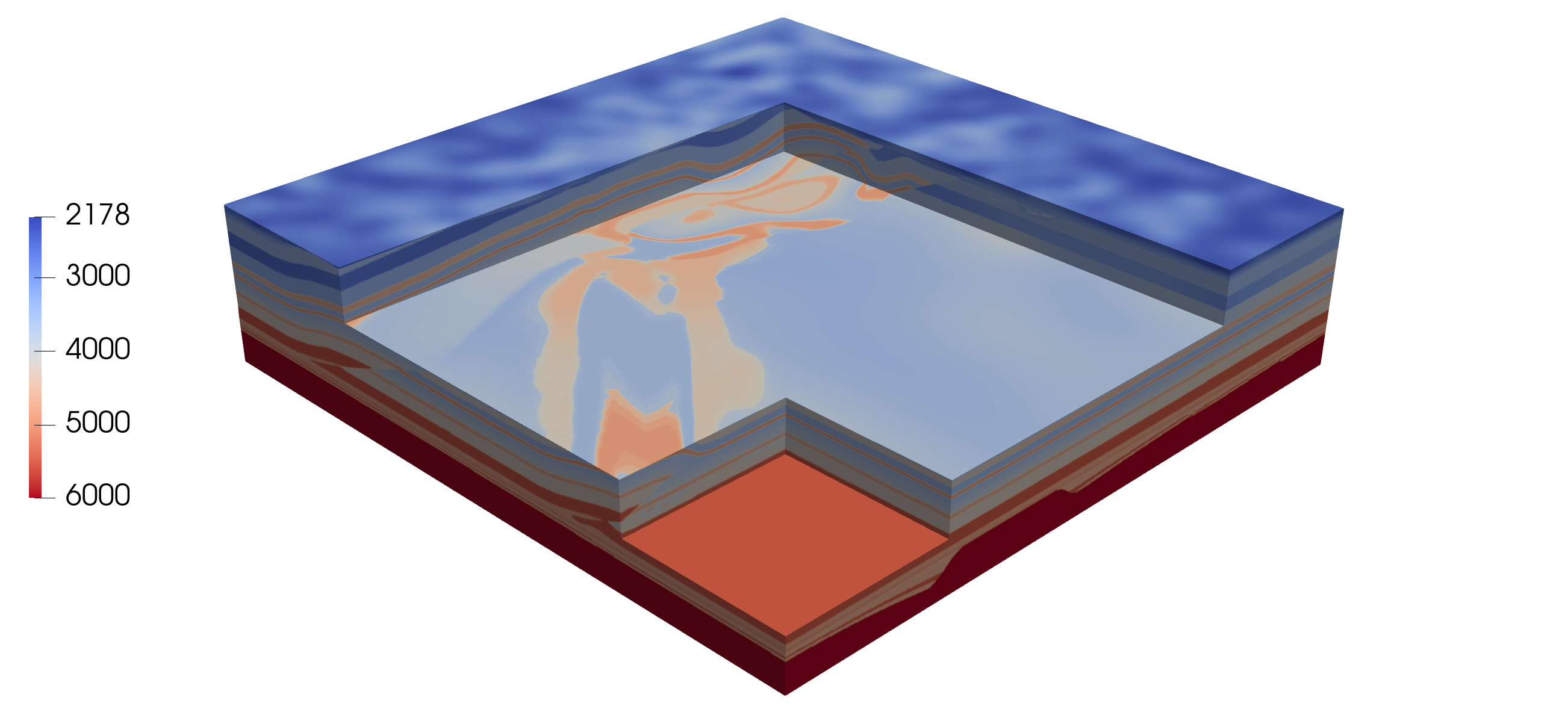}
\caption{ 
Compressional wave speed of the overthrust model. Portions of the model are carved out to illustrate the internal structures. Colorbar represents value of the compressional wave speed in unit m/s. 
}
\label{fig_overthrust_model}
\end{figure}

In this example, we test the proposed layer-wise discretization with the overthrust model, whose compressional wave speed is illustrated in Figure~\ref{fig_overthrust_model}. 
The model is described by 801 $\times$ 801 $\times$ 181 grid points. 
%
%
Grid spacing between two neighboring parameter grid points is assigned to be 10 m.
Since the overthrust model contains only the information of the compressional wave speed $c_p$, we synthesize the shear wave speed $c_s$ by prescribing the $\nicefrac{c_p}{c_s}$ ratio, which decreases from top (4.5) to bottom (1.5) linearly. 
(We note here that the $\nicefrac{c_p}{c_s}$ ratio tends to decrease with depth and increased degree of consolidation in rocks; see~\cite{gregory1976fluid, hamilton1979v}, for example.)
Moreover, the density $\rho$ is also synthesized, which increases from top (2000 kg/m$^3$) to bottom (2700 kg/m$^3$) linearly.
%

In nonuniform grid simulation, the entire simulation domain is separated into three layers of equal depth (600 m). The most relevant physical and numerical parameters associated with these layers are displayed in Table~\ref{overthrust_model_parameters}.
The source profile is specified as the Ricker wavelet with central frequency at 5 Hz. The maximum frequency in the wave content is assumed to be 12.5 Hz, which is used to calculate the $N_\text{ppw}$ shown in Table~\ref{overthrust_model_parameters}.
The grid spacing ratios at both layer interfaces are 1:2.
The temporal step size $\Delta t$ is specified as $\sim$$2.3072 \times 10^{-4}$ s, which corresponds to equation \eqref{appendix_Delta_t_i} with $\mathbb{C}_{\text{CFL}} = 0.8$.
For comparison, the uniform grid simulation with $\Delta x = 2.5$ for all layers is performed.
The temporal step size used in the uniform grid simulation is $\sim$$1.9245 \times 10^{-4}$ s, which corresponds to equation \eqref{appendix_Delta_t} with $\mathbb{C}_{\text{CFL}} = 0.8$.
Between the two simulations, the ratio of spatial discretization points is $\sim$$2.63$; the ratio of temporal discretization points is $\sim$$1.2$; the ratio of total spatial-temporal discretization points is $\sim$$3.153$.

\begin{table}[H]
\captionsetup{width=0.9\textwidth, font=small,labelfont=small}
\caption{Physical and numerical parameters associated with the overthrust model.}
\label{overthrust_model_parameters}
\centering
\begin{tabular}{| c | c | c | c | c |}
\hline
             & min. $c_s$          & max. $c_p$           & $\Delta x$ & $N_\text{ppw}$ \\ \hline
Layer 1 & $\sim$484.2 m/s & $\sim$5004.7 m/s & 2.5 m        & $\sim$15.5        \\ \hline
Layer 2 & 958.5 m/s           & 5500 m/s               & 5 m           & $\sim$15.3        \\ \hline
Layer 3 & 1920 m/s            & 6000 m/s               & 10 m         & $\sim$15.4        \\ \hline 
\end{tabular}
\end{table}

In physical units, the size of the simulation domain is 8000 m $\times$ 8000 m $\times$ 1800 m.
A compressional source is applied on the normal stress components. The $V_z$ component of the velocity is recorded for comparison.
The source and receiver are placed at 4000 m apart from each other on both the $x$- and $y$-directions, and at 12.5 m (i.e., 5 grid points) and 11.25 m (i.e., 4.5 grid points) below the top boundary, respectively.

\begin{figure}[H]
\captionsetup{width=0.9\textwidth, font=small,labelfont=small}
\centering\includegraphics[scale=0.1]{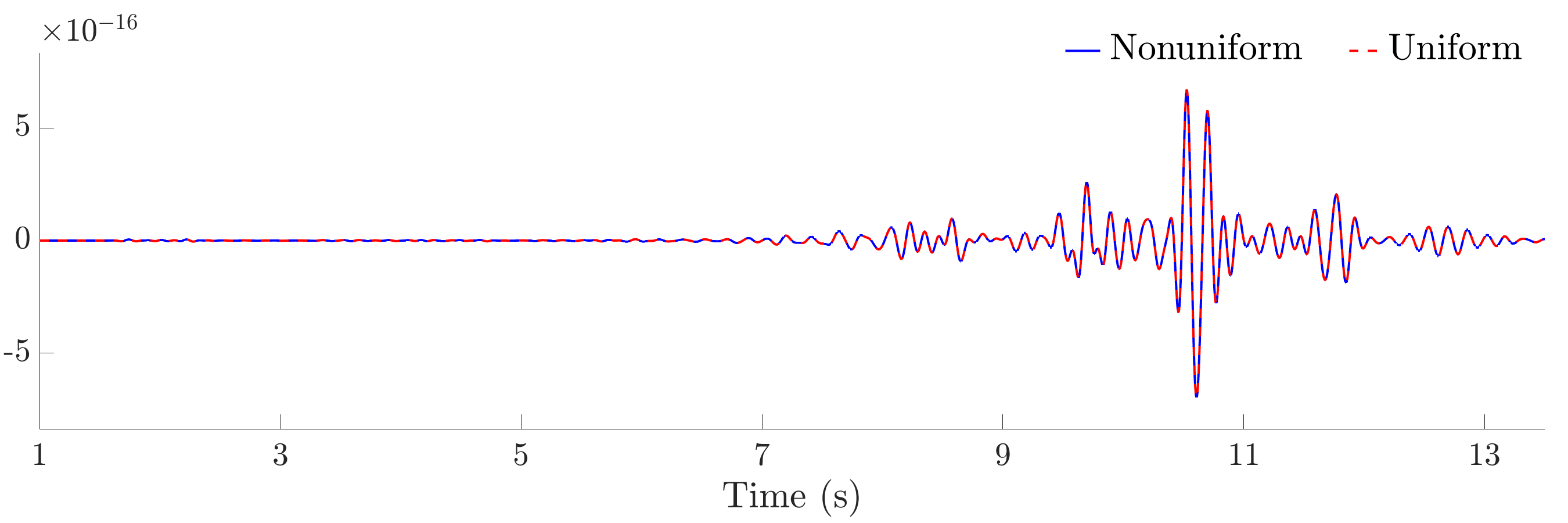}
\caption{ 
Overthrust model. Recorded seismograms from nonuniform and uniform grid simulations.
}
\label{fig_overthrust_seismogram}
\end{figure}

\begin{figure}[H]
\captionsetup{width=0.9\textwidth, font=small,labelfont=small}
\centering\includegraphics[scale=0.1]{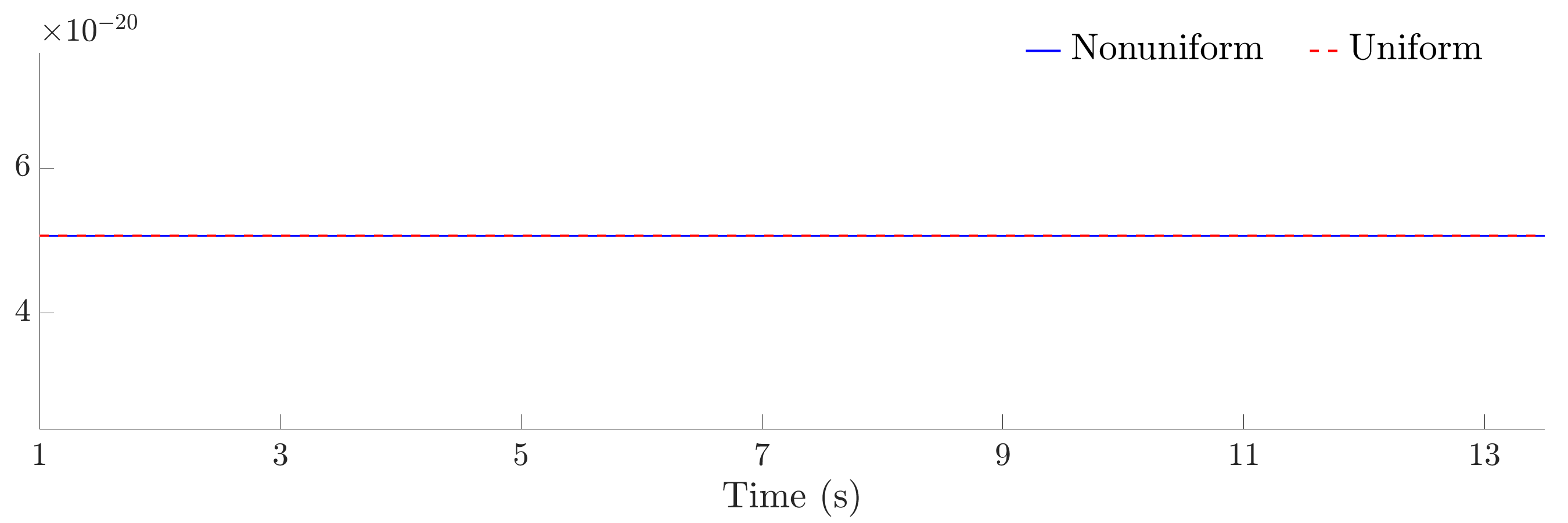}
\caption{ 
Overthrust model. Evolution of the discrete energy in nonuniform and uniform grid simulations.
}
\label{fig_overthrust_energy}
\end{figure}

The recorded seismograms are displayed in Figure~\ref{fig_overthrust_seismogram}, where we observe satisfactory agreement between these two simulations. 
After projecting these seismograms onto a common set of equidistant time instances spanning the interval $[1\text{ s},13.5\text{ s}]$ with increment 0.001 s using cubic spline interpolation, the relative difference measured in $\ell_2$ norm is $\sim$$5.69\times 10^{-4}$.
Moreover, evolutions of the discrete energy in both simulations are displayed in Figure~\ref{fig_overthrust_energy}, which remain constant after the source term tampers off and hence corroborate the energy conserving property.

For this example, the nonuniform and uniform simulations are run with 7300 and 19200 processor cores, respectively, each responsible for a subdomain with $80 \times 80 \times 60$ $M$-grid points.
The elapsed wall clock times are $\sim$10860 s and $\sim$12787 s, respectively, translating to a ratio of $\sim$3.095 in consumed computational resources (i.e., number of processors multiplied by elapsed time), which agrees well with the ratio of the total number of spatial-temporal discretization points in these two simulations (3.153), suggesting an overhead for the nonuniform simulation at $\sim$2\%.
A table collecting the above expenditure information on computing resources can be found in the {\it Supplementary Material} \ref{supp_additional_tables}.
It is omitted here to conserve space.

\section{Conclusions} \label{section_conclusion}
Layer-wise finite difference discretization of the 3D isotropic elastic wave propagation on staggered grids with nonconforming interfaces between layers is considered in this work.
It is shown that energy conserving simulation across the nonconforming layer interfaces can be achieved through the use of summation-by-parts finite difference operators, simultaneous approximation terms, and compatible interpolation operators, which ensure the stability of the simulation. 
Moreover, by examining the standard process of determining the discretization parameters for seismic wave simulations under the usual constraints, as presented in \ref{appendix_thought_experiment}, and by the numerical examples, it is demonstrated that, compared to finite difference discretization on uniform grids, significant reduction of computational resources can be achieved by employing the proposed layer-wise discretization approach, thus enabling more efficient processing of seismic data.

\appendix

\section{Impact of wave speed variation on discretization parameters}\label{appendix_thought_experiment}
In the following, we use a simple thought experiment to illustrate that, for a simulation domain composed of different subterranean media, a layer-wise spatial discretization grid can significantly reduce the expenditure in computational resources compared to its fully uniform counterpart.

First, considering a uniformly discretized simulation domain associated with a medium whose wave speed ranges from $c_{\min}$ to $c_{\max}$,
the minimum wavelength of the simulated waves can be expressed as
\begin{equation}
\label{appendix_lambda}
\lambda_{\min} \ = \ \frac{c_{\min}}{f}\, ,
\end{equation}
where $f$ denotes the highest frequency in the wave content being simulated.
Assuming that the spatial grid spacing $\Delta x$ is determined by prescribing $N_\text{ppw}$ grid points to resolve the minimum wavelength $\lambda_{\min}$, we have that
\begin{equation}
\label{appendix_Delta_x}
\Delta x \ = \ \frac{\lambda_{\min}}{N_\text{ppw}} \ = \ \frac{c_{\min}}{f \cdot N_\text{ppw}}\, .
\end{equation}
Further assuming that the temporal step size $\Delta t$ is determined by the CFL constraint, we can write 
\begin{equation}
\label{appendix_Delta_t}
\Delta t \ = \ \mathbb{C}_{\text{CFL}} \cdot \frac{\Delta x}{c_{\max}}
\ = \ \frac{\mathbb{C}_{\text{CFL}}}{f \cdot N_\text{ppw}} \cdot \frac{c_{\min}}{c_{\max}}\, ,
\end{equation}
where $\mathbb{C}_{\text{CFL}}$ is a constant that depends on the spatial and temporal discretization schemes, but is independent of the underlying medium.

Regarding $f$, $N_\text{ppw}$, and $\mathbb{C}_{\text{CFL}}$ as fixed constants in the above formulas, 
we observe from equation \eqref{appendix_Delta_x} that the spatial grid spacing $\Delta x$ is proportional to the minimum wave speed $c_{\min}$ of the underlying medium;
moreover, we observe from equation \eqref{appendix_Delta_t} that the temporal step size $\Delta t$ is constrained by the maximum contrast in the wave speed, i.e., $\nicefrac{c_{\max}}{c_{\min}}$, of the underlying medium.

On the other hand, if the same simulation domain is separated into $N_\text{L}$ layers, 
each discretized uniformly by prescribing $N_\text{ppw}$ grid points to resolve the minimum wavelength within the layer, 
the spatial grid spacing for each layer can be expressed as 
\begin{equation}
\label{appendix_Delta_x_i}
\Delta x^i \ = \ \frac{\lambda^i_{\min}}{N_\text{ppw}} \ = \ \frac{c^i_{\min}}{f \cdot N_\text{ppw}}\, , \quad i = 1 \cdots N_\text{L}\, ,
\end{equation}
where $\Delta x^i$, $\lambda^i_{\min}$, and $c^i_{\min}$ stand for the spatial grid spacing, minimum wavelength, and minimum wave speed of the $i$th layer. 
In other words, the spatial grid spacing is now only limited by the minimum wave speed within the layer, instead of the minimum wave speed of the entire simulation domain. 
Moreover, assuming that constant $\mathbb{C}_{\text{CFL}}$ remains unchanged, the temporal step size is now determined by the most stringent CFL constraint among all layers, i.e., 
\begin{equation}
\label{appendix_Delta_t_i}
\Delta t 
\ = \ 
\mathbb{C}_{\text{CFL}} \cdot \min_{i = 1 \cdots N_\text{L}} \left\{ \frac{\Delta x^i}{c^i_{\max}} \right\}
\ = \ 
\frac{\mathbb{C}_{\text{CFL}}}{f \cdot N_\text{ppw}} \cdot \min_{i = 1 \cdots N_\text{L}} \left\{ \frac{c^i_{\min}}{c^i_{\max}} \right\} \, ,
\end{equation}
which is likely to be relaxed from, at least not worse than, that of equation \eqref{appendix_Delta_t}.

\begin{figure}[b!]
\captionsetup{width=0.9\textwidth, font=small,labelfont=small}
\centering
\begin{subfigure}[b]{0.3\textwidth}
\captionsetup{width=1\textwidth, font=small,labelfont=small}
\centering\includegraphics[scale=0.0225]{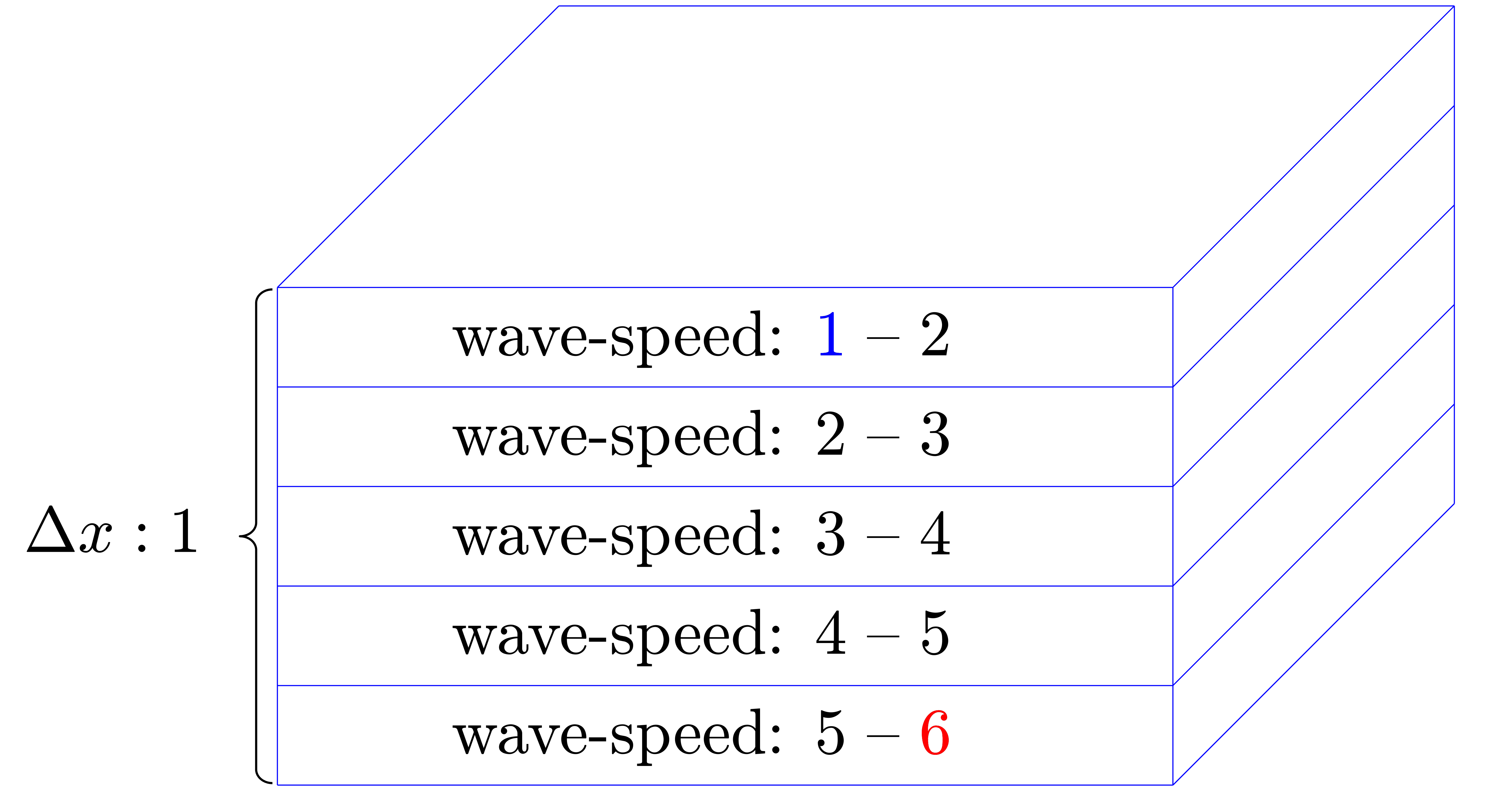}
\caption{} 
\label{fig_layer_model_1}
\end{subfigure}
\hspace*{\fill}
\begin{subfigure}[b]{0.3\textwidth}
\captionsetup{width=1\textwidth, font=small,labelfont=small}
\centering\includegraphics[scale=0.0225]{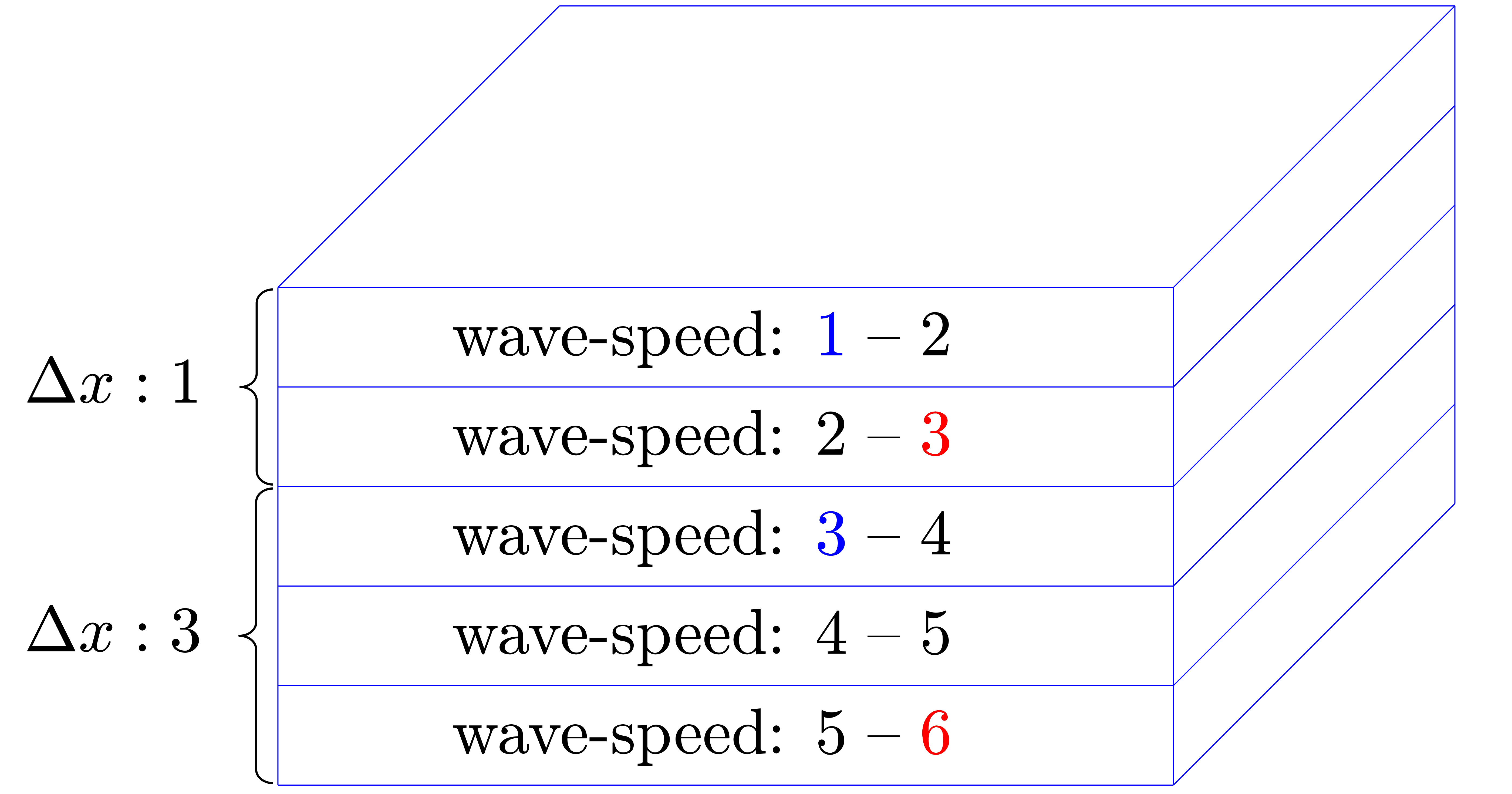}
\caption{} 
\label{fig_layer_model_2}
\end{subfigure} 
\hspace*{\fill}
\begin{subfigure}[b]{0.3\textwidth}
\captionsetup{width=1\textwidth, font=small,labelfont=small}
\centering\includegraphics[scale=0.0225]{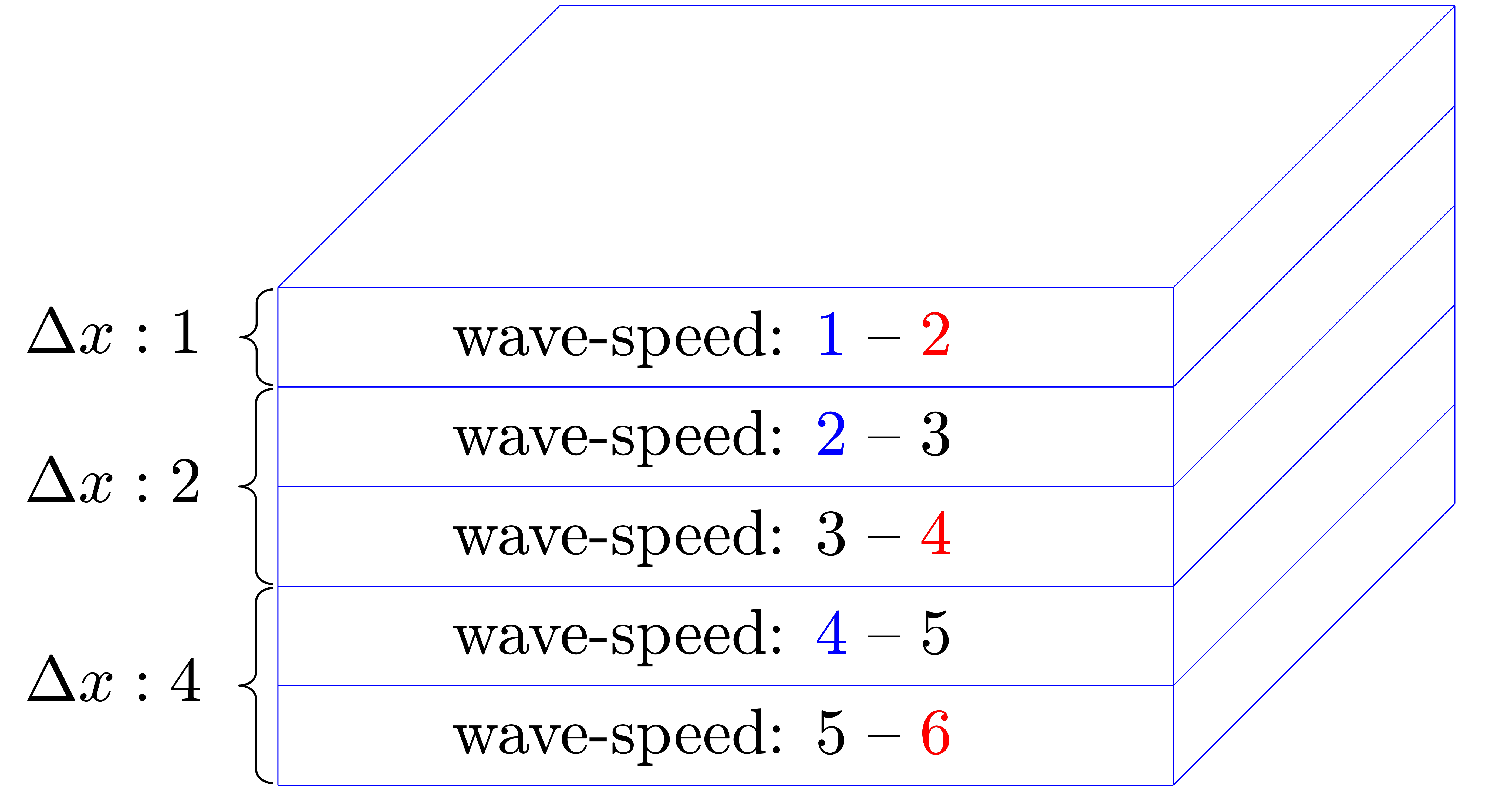}
\caption{} 
\label{fig_layer_model_3}
\end{subfigure} 
\hspace*{\fill}
\caption{
A simple model to illustrate the potential benefit of layer-wise uniform discretization over uniform discretization for wave simulations in earth media. 
The same simulation domain is considered in all three cases: 
a) the entire simulation domain is discretized uniformly;
b) the entire simulation domain is separated into two layers, each discretized uniformly;
c) the entire simulation domain is separated into three layers, each discretized uniformly. 
The normalized minimum and maximum wave speeds of each uniform discretization region are colored in blue and red, respectively; 
the former determines the spatial grid spacing $\Delta x$, while their ratio constrains the temporal step size $\Delta t$.
} 
\label{fig_layer_model}
\end{figure}

To give a concrete example, consider the simulation domain depicted in Figure~\ref{fig_layer_model}, 
which consists of five regions of equal depth with progressively increasing wave speeds.
The normalized minimum and maximum wave speeds of each region are displayed in the figure.

Figure~\ref{fig_layer_model_1} depicts the case that the entire simulation domain is discretized uniformly, which is used as the reference for comparison. 
The total number of its space-time discretization grid points is counted as 
\begin{equation}
\label{total_grid_points_layer_model_1}
\left(1 + 1 + 1 + 1 + 1\right) \times 6 \ = \ 30 \text{ units,}
\end{equation}
where the quantity inside the parentheses corresponds to the total spatial discretization grid points, while the number 6 is related to the temporal discretization grid points, stemming from the maximum contrast in the wave speed of the entire simulation domain; see equation \eqref{appendix_Delta_t}.

Figure~\ref{fig_layer_model_2} depicts the case that the entire simulation domain is separated into two layers, each discretized uniformly and with a 1:3 ratio in $\Delta x$, the total number of space-time discretization grid points becomes
\begin{equation}
\label{total_grid_points_layer_model_2}
\left(1 + 1 + \tfrac{1}{27} + \tfrac{1}{27} + \tfrac{1}{27}\right) \times 3 \ \approx \ 6.3333 \text{ units,}
\end{equation}
where the number 3 stems from the maximum contrast in the wave speed of the top layer; see equation \eqref{appendix_Delta_t_i}.
By comparing equation \eqref{total_grid_points_layer_model_2} to equation \eqref{total_grid_points_layer_model_1}, we observe that the total number of space-time discretization grid points has been reduced by a factor of $\sim\!4.7$.

Figure~\ref{fig_layer_model_3} depicts the case that the entire simulation domain is separated into three layers, each discretized uniformly and with a 1:2:4 ratio in $\Delta x$, the total number of space-time discretization grid points becomes
\begin{equation}
\label{total_grid_points_layer_model_3}
\left(1 + \tfrac{1}{8} + \tfrac{1}{8} + \tfrac{1}{64} + \tfrac{1}{64}\right) \times 2 \ = \ 2.5625 \text{ units,}
\end{equation}
where the number 2 stems from the maximum contrast in the wave speed of the top layer (or the middle layer); see equation \eqref{appendix_Delta_t_i}.
By comparing equation \eqref{total_grid_points_layer_model_3} to equation \eqref{total_grid_points_layer_model_1}, we observe that the total number of space-time discretization grid points has been reduced by a factor of $\sim\!11.7$.

We note here that for realistic earth media, it rarely happens that different layers can be separated in such a non-overlapping fashion as shown in Figure~\ref{fig_layer_model}, where the maximum wave speed in one layer matches the minimum wave speed of the next.
Instead, there are often overlaps in the wave speed ranges of two neighboring layers.
In such a situation, the cost saving factors estimated above can be overly optimistic. 
On the other hand, instead of the equal-depth model considered in Figure~\ref{fig_layer_model}, a relatively thin layer of low velocity surface material (e.g., sand and soil) often appears in realistic situations, which can severely restrict the choices of discretization parameters in uniform grid simulations.  
In such cases, the actual cost saving factors can be more pronounced than those estimated above.
Finally, addressing layer interfaces can incur computational overheads in simulations, which may also affect the actual cost saving factors. 
However, impact from such overheads tends to be insignificant for large simulations at realistic scales, as evidenced by the numerical examples in Section \ref{numerical_examples}.

\section{Compatible interpolation operators}\label{appendix_interpolation_operators}
A collection of 1D interpolation operators for a variety of grid spacing ratios are presented in the following. 
They are building blocks for the 2D interpolation operators used in Section \ref{interface_and_interpolation} to address the nonconforming interfaces.
These interpolation operators are referred to as being compatible (with each other), since, in addition to the typical approximation requirements, they also satisfy the reciprocal relationship postulated in equation \eqref{1D_interpolation_operators_properties}.

Because of this reciprocal relationship, only the interpolation operators from the coarse grid (the minus side) to the fine grid (the plus side) are presented since the interpolation operators for the opposite direction can be deduced based on equation \eqref{1D_interpolation_operators_properties}.
Furthermore, only formulas for a few exemplary grid points within a repeating cycle are shown in the following, while formulas for the other grid points can be replicated by shifting the indices accordingly.

\subsubsection*{1:2 grid spacing ratio}
\begin{figure}[H]
\captionsetup{width=0.9\textwidth, font=small,labelfont=small}
\centering
%
\centering\includegraphics[scale=0.125]{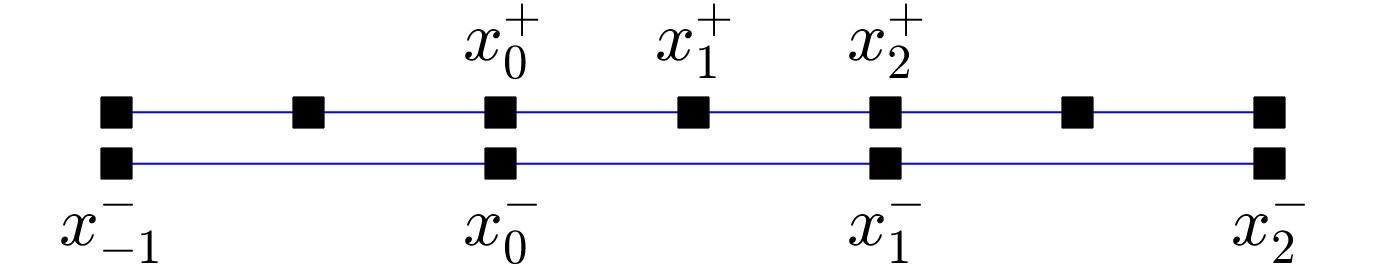}
%
%
\centering\includegraphics[scale=0.125]{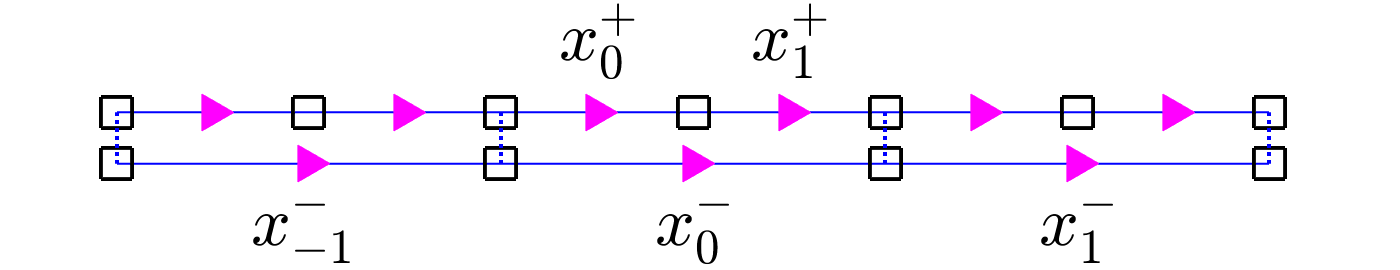}
%
%
\caption{Illustration of the grid points involved in formulas \eqref{T_N_1_2} and \eqref{T_M_1_2} for interpolation operators $T^{N}_{x^-}$ (left) and $T^{M}_{x^-}$ (right). The grid spacing ratio is 1:2. } 
\label{Grid_configuration_transfer_operator_1_2}
\end{figure}
{\footnotesize
\begin{flalign}
\label{T_N_1_2}
\indent
T^{N}_{x^-}: \ & \nonumber \\
&
\arraycolsep=2.pt\def\arraystretch{0.625}
\begin{array}{rrrrrrrrrr}
f(x^+_0) & \ \leftarrow 
				& &
				& & f(x^-_0) 
				& &
				& &
\\
f(x^+_1) & \ \leftarrow 
			 & - & \tfrac{1}{16} f(x^-_{-1}) 
                          & + & \tfrac{9}{16} f(x^-_0) 
                          & + & \tfrac{9}{16} f(x^-_1) 
                          & - & \tfrac{1}{16} f(x^-_2) 
\\
f(x^+_2) & \ \leftarrow 
				& &
				& &
				& & f(x^-_1) 
				& &
\end{array}
&
\end{flalign}
}
{\footnotesize
\begin{flalign}
\small
\label{T_M_1_2}
\indent
T^{M}_{x^-}: \ & \nonumber \\
&
\arraycolsep=2.pt\def\arraystretch{0.625}
\begin{array}{rrrrrrrrrr}
f(x^+_0) & \ \leftarrow 
			&    & \tfrac{ 5}{32} f(x^-_0) 
                         & + & \tfrac{15}{16} f(x^-_1) 
                         & -  & \tfrac{ 3}{32} f(x^-_2)
\\
f(x^+_1) & \ \leftarrow 
			& -  & \tfrac{ 3}{32} f(x^-_0) 
                         & + & \tfrac{15}{16} f(x^-_1) 
                         & + & \tfrac{ 5}{32} f(x^-_2)
\end{array}
&
\end{flalign}
}

\subsubsection*{1:3 grid spacing ratio}
\begin{figure}[H]
\captionsetup{width=0.9\textwidth, font=small,labelfont=small}
\centering
%
\centering\includegraphics[scale=0.125]{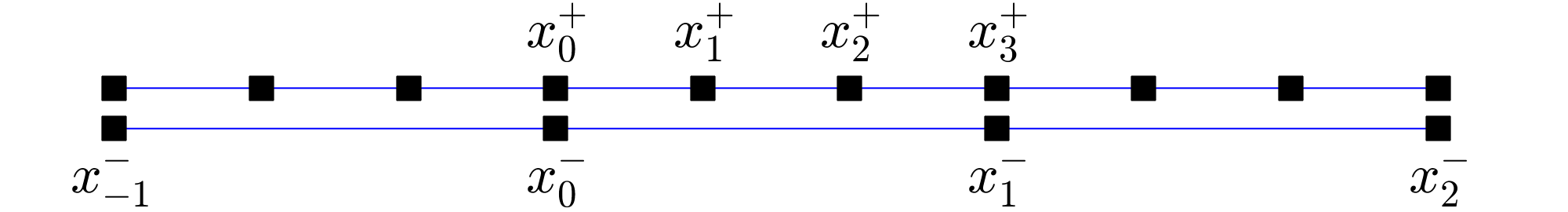}
%
\\
%
\centering\includegraphics[scale=0.125]{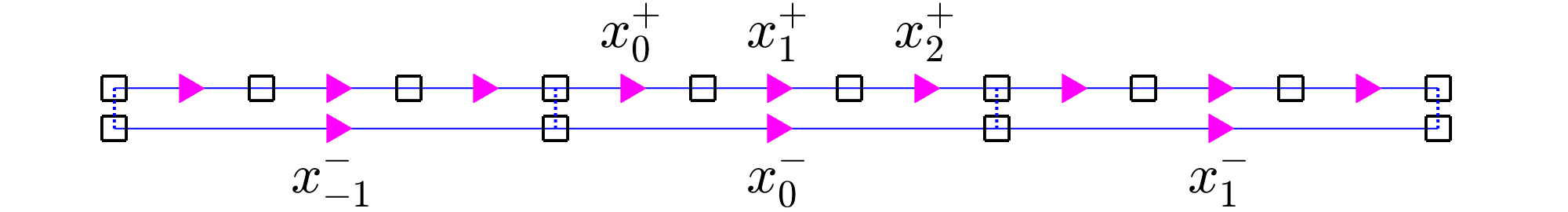}
%
%
\caption{Illustration of the grid points involved in formulas \eqref{T_N_1_3} and \eqref{T_M_1_3} for interpolation operators $T^{N}_{x^-}$ (top) and $T^{M}_{x^-}$ (bottom). The grid spacing ratio is 1:3. } 
\label{Grid_configuration_transfer_operator_1_3}
\end{figure}
{\footnotesize
\begin{flalign}
\label{T_N_1_3}
\indent
T^{N}_{x^-}: \ & \nonumber \\
&
\arraycolsep=2.pt\def\arraystretch{0.625}
\begin{array}{rrrrrrrrrr}
f(x^+_0) & \ \leftarrow 
				&&
				&&f(x^-_0) 
				&&
				&&
\\
f(x^+_1) & \ \leftarrow
			 &-&\tfrac{1}{9} f(x^-_{-1}) 
                          &+&\tfrac{8}{9} f(x^-_0) 
                          &+&\tfrac{2}{9} f(x^-_1)
                          &&
\\
f(x^+_2) & \ \leftarrow 
			&&
			&& \tfrac{2}{9} f(x^-_0) 
                         &+&\tfrac{8}{9} f(x^-_1) 
                         &-&\tfrac{1}{9} f(x^-_2) 
\\
f(x^+_3) & \ \leftarrow 
			&&
			&&
			&& f(x^-_1) 
			&&
\end{array}
&
\end{flalign}
}
{\footnotesize
\begin{flalign}
\label{T_M_1_3}
\indent
T^{M}_{x^-}: \ & \nonumber \\
&
\arraycolsep=2.pt\def\arraystretch{0.625}
\begin{array}{rrrrrrrrrr}
f(x^+_0) & \ \leftarrow 
			 && \tfrac{2}{9} f(x^-_{-1}) 
                          &+&\tfrac{8}{9} f(x^-_0) 
                          &-&\tfrac{1}{9} f(x^-_1)
\\
f(x^+_1) & \ \leftarrow 
			&&
			&& f(x^-_0) 
			&&
\\
f(x^+_2) & \ \leftarrow 
			 &-&\tfrac{1}{9} f(x^-_{-1}) 
                          &+&\tfrac{8}{9} f(x^-_0) 
                          &+&\tfrac{2}{9} f(x^-_1) 
\end{array}
&
\end{flalign}
}

\subsubsection*{2:3 grid spacing ratio}
\begin{figure}[H]
\captionsetup{width=0.9\textwidth, font=small,labelfont=small}
\centering
%
\centering\includegraphics[scale=0.125]{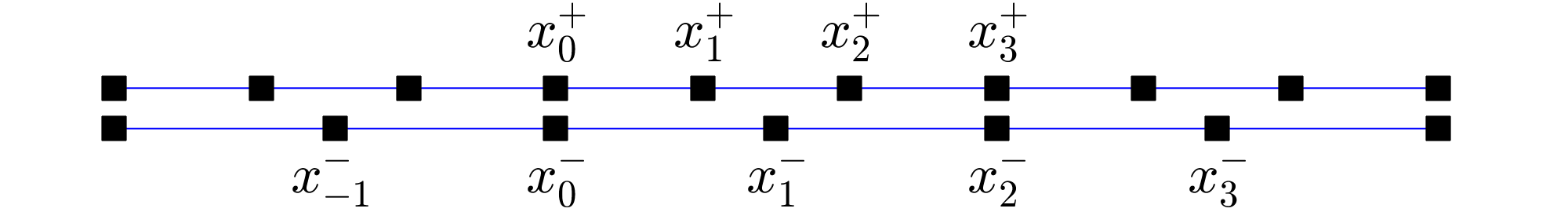}
%
\\
%
\centering\includegraphics[scale=0.125]{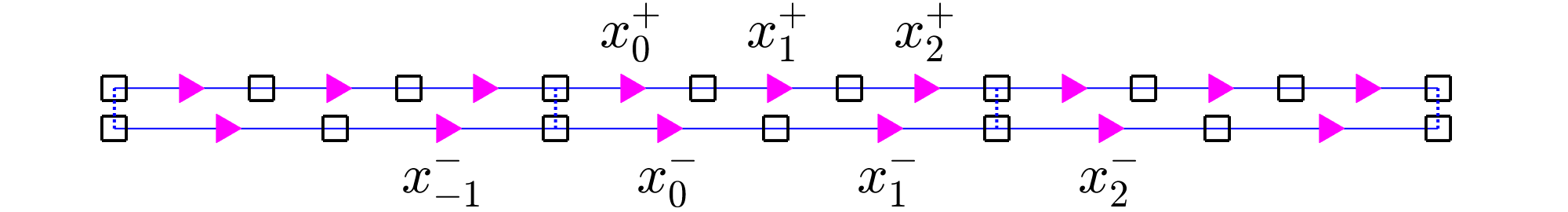}
%
%
\caption{Illustration of the grid points involved in formulas \eqref{T_N_2_3} and \eqref{T_M_2_3} for interpolation operators $T^{N}_{x^-}$ (top) and $T^{M}_{x^-}$ (bottom). The grid spacing ratio is 2:3.} 
\label{Grid_configuration_transfer_operator_2_3}
\end{figure}
{\footnotesize
\begin{flalign}
\label{T_N_2_3}
\indent
T^{N}_{x^-}: \ & \nonumber \\
&
\arraycolsep=2.pt\def\arraystretch{0.625}
\begin{array}{rrrrrrrrrrrr}
f(x^+_0) & \ \leftarrow 
			&& 
			&& f(x^-_0) 
			&&
			&&
			&& 
\\
f(x^+_1) & \ \leftarrow 
			 &-&\tfrac{ 11}{288} f(x^-_{-1})
                          &+&\tfrac{  1}{  3} f(x^-_{ 0})
                          &+&\tfrac{113}{144} f(x^-_{ 1})
                          &-&\tfrac{  1}{ 12} f(x^-_{ 2})
                          &+&\tfrac{  1}{288} f(x^-_{ 3}) 
\\
f(x^+_2) & \ \leftarrow 
			 & & \tfrac{  1}{288} f(x^-_{-1})
                          &-&\tfrac{  1}{ 12} f(x^-_{ 0})
                          &+&\tfrac{113}{144} f(x^-_{ 1})
                          &+&\tfrac{  1}{  3} f(x^-_{ 2})
                          &-&\tfrac{ 11}{288} f(x^-_{ 3})
\\
f(x^+_3) & \ \leftarrow 
			&&
			&&
			&&
			&& f(x^-_2) 
			&&
\end{array}
&
\end{flalign}
}
{\footnotesize
\begin{flalign}
\label{T_M_2_3}
\indent
T^{M}_{x^-}: \ & \nonumber \\
&
\arraycolsep=2.pt\def\arraystretch{0.625}
\begin{array}{rrrrrrrrrrrr}
f(x^+_0) & \ \leftarrow 
			&& \tfrac{ 101}{1152} f(x^-_{-1})
                         &+&\tfrac{1153}{1152} f(x^-_{ 0})
                         &-&\tfrac{ 113}{1152} f(x^-_{ 1})
                         &+&\tfrac{  11}{1152} f(x^-_{ 2})
\\
f(x^+_1) & \ \leftarrow 
			 &-&\tfrac{1}{16} f(x^-_{-1})
                          &+& \tfrac{9}{16} f(x^-_{ 0})
                          &+& \tfrac{9}{16} f(x^-_{ 1})
                          &-& \tfrac{1}{16} f(x^-_{ 2}) 
\\
f(x^+_2) & \ \leftarrow 
			 &&  \tfrac{  11}{1152} f(x^-_{-1})
                          &-&\tfrac{ 113}{1152} f(x^-_{ 0})
                          &+&\tfrac{1153}{1152} f(x^-_{ 1})
                          &+&\tfrac{ 101}{1152} f(x^-_{ 2}) 
\end{array} 
&
\end{flalign}
}

\subsubsection*{1:4 grid spacing ratio} 
\begin{figure}[H]
\captionsetup{width=0.9\textwidth, font=small,labelfont=small}
\centering
%
\centering\includegraphics[scale=0.125]{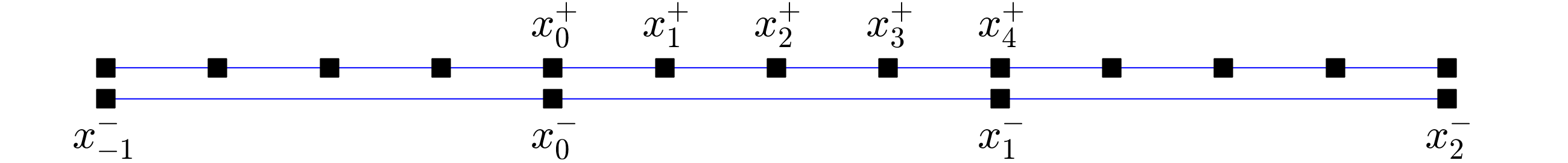}
%
\\[0.5ex]
%
\centering\includegraphics[scale=0.125]{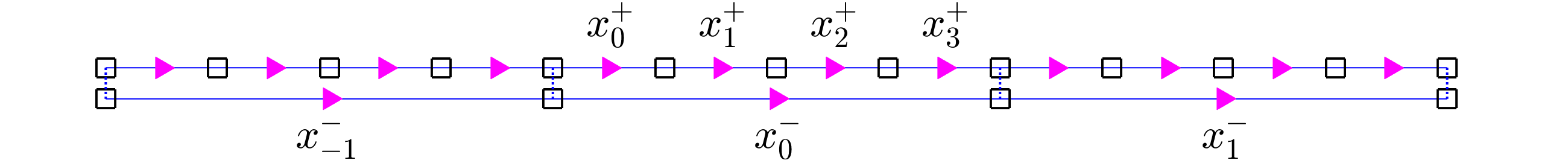}
%
%
\caption{Illustration of the grid points involved in formulas \eqref{T_N_1_4} and \eqref{T_M_1_4} for interpolation operators $T^{N}_{x^-}$ (top) and $T^{M}_{x^-}$ (bottom). The grid spacing ratio is 1:4. } 
\label{Grid_configuration_transfer_operator_1_4}
\end{figure}
{\footnotesize
\begin{flalign}
\label{T_N_1_4}
\indent
T^{N}_{x^-}: \ & \nonumber \\
&
\arraycolsep=2.pt\def\arraystretch{0.625}
\begin{array}{rrrrrrrrrrrr}
f(x^+_0) & \ \leftarrow 
				&&
				&& f(x^-_0) 
				&&
				&& 
\\
f(x^+_1) & \ \leftarrow
			 &-&\tfrac{ 3}{32} f(x^-_{-1})
                          &+&\tfrac{15}{16} f(x^-_{ 0})
                          &+&\tfrac{ 5}{32} f(x^-_{ 1}) 
                          &&
\\
f(x^+_2) & \ \leftarrow 
			  &-&\tfrac{1}{16} f(x^-_{-1}) 
                           &+&\tfrac{9}{16} f(x^-_{ 0}) 
                           &+&\tfrac{9}{16} f(x^-_{ 1})
                           &-&\tfrac{1}{16} f(x^-_{ 2}) 
\\
f(x^+_3) & \ \leftarrow 
			 &&
			 && \tfrac{ 5}{32} f(x^-_{ 0}) 
                          &+&\tfrac{15}{16} f(x^-_{ 1})
                          &-&\tfrac{ 3}{32} f(x^-_{ 2}) 
\\
f(x^+_4) & \ \leftarrow 
			&&
			&&
			&&  f(x^-_1) 
			&& 
\end{array} &
\end{flalign}
}
{\footnotesize
\begin{flalign}
\label{T_M_1_4}
\indent
T^{M}_{x^-}: \ & \nonumber \\
&
\arraycolsep=2.pt\def\arraystretch{0.625}
\begin{array}{rrrrrrrrrrrr}
f(x^+_0) & \ \leftarrow  
			 && \tfrac{33}{128} f(x^-_{-1})
                          &+&\tfrac{55}{ 64} f(x^-_{ 0})
                          &-&\tfrac{15}{128} f(x^-_{ 1}) 
\\
f(x^+_1) & \ \leftarrow  
			 && \tfrac{ 9}{128} f(x^-_{-1})
                          &+&\tfrac{63}{ 64} f(x^-_{ 0})
                          &-&\tfrac{ 7}{128} f(x^-_{ 1}) 
\\
f(x^+_2) & \ \leftarrow 
			 &-&\tfrac{ 7}{128} f(x^-_{-1})
                          &+&\tfrac{63}{ 64} f(x^-_{ 0})
                          &+&\tfrac{ 9}{128} f(x^-_{ 1})
\\
f(x^+_3) & \ \leftarrow 
			 &-&\tfrac{15}{128} f(x^-_{-1})
                          &+&\tfrac{55}{ 64} f(x^-_{ 0})
                          &+&\tfrac{33}{128} f(x^-_{ 1})
\end{array}
&
\end{flalign}
}

\subsubsection*{3:4 grid spacing ratio}
\begin{figure}[H]
\captionsetup{width=0.9\textwidth, font=small,labelfont=small}
\centering
%
\centering\includegraphics[scale=0.125]{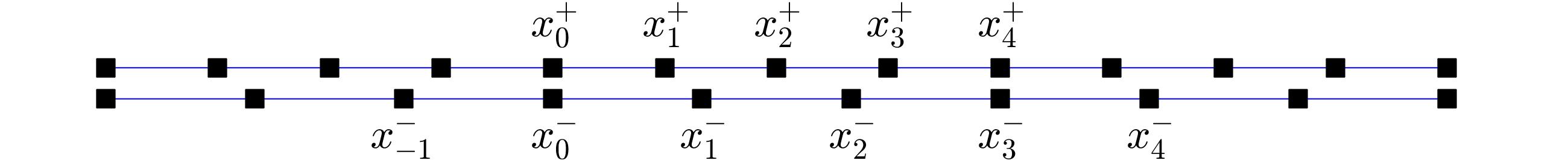}
%
\\[0.5ex]
%
\centering\includegraphics[scale=0.125]{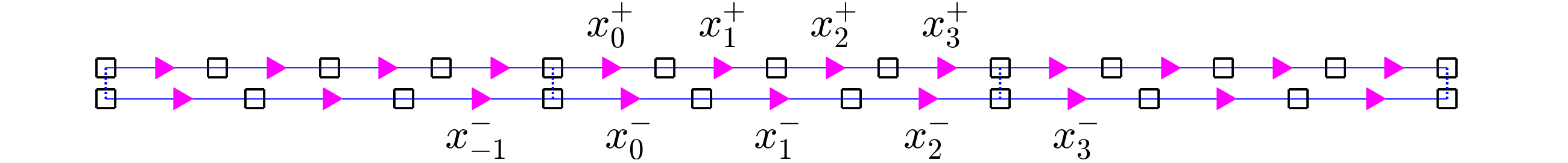}
%
%
\caption{Illustration of the grid points involved in formulas \eqref{T_N_3_4} and \eqref{T_M_3_4} for interpolation operators $T^{N}_{x^-}$ (top) and $T^{M}_{x^-}$ (bottom). The grid spacing ratio is 3:4. } 
\label{Grid_configuration_transfer_operator_3_4}
\end{figure}
{\footnotesize
\begin{flalign}
\label{T_N_3_4}
T^{N}_{x^-}: \ & \nonumber \\
&
\arraycolsep=2.pt\def\arraystretch{0.625}
\begin{array}{rrrrrrrrrrrrrrrrrr}
f(x^+_0) & \ \leftarrow  
			&&
			&& f(x^-_0) 
			&&
			&&
			&&
			&&
\\
f(x^+_1) & \ \leftarrow  
			 &-&\tfrac{119}{4320} f(x^-_{-1})
                          &+&\tfrac{ 67}{ 288} f(x^-_{ 0})
                          &+&\tfrac{629}{ 720} f(x^-_{ 1})
                          &-&\tfrac{367}{4320} f(x^-_{ 2})
                          &+&\tfrac{  1}{ 160} f(x^-_{ 3}) 
                          &&
\\
f(x^+_2) & \ \leftarrow  
			 && \tfrac{   7}{2160} f(x^-_{-1})
                          &-&\tfrac{  13}{ 180} f(x^-_{ 0})
                          &+&\tfrac{1229}{2160} f(x^-_{ 1})
                          &+&\tfrac{1229}{2160} f(x^-_{ 2})
                          &-&\tfrac{  13}{ 180} f(x^-_{ 3})
                          &+&\tfrac{   7}{2160} f(x^-_{ 4})
\\
f(x^+_3) & \ \leftarrow  
			 &&
			 && \tfrac{  1}{ 160} f(x^-_{ 0}) 
                          &-&\tfrac{367}{4320} f(x^-_{ 1}) 
                          &+&\tfrac{629}{ 720} f(x^-_{ 2}) 
                          &+&\tfrac{ 67}{ 288} f(x^-_{ 3}) 
                          &-&\tfrac{119}{4320} f(x^-_{ 4})
\\
f(x^+_4) & \ \leftarrow  
			&&
			&&
			&&
			&&
			&& f(x^-_3) 
			&&
\end{array}
&
\end{flalign}
}
{\footnotesize
\begin{flalign}
\label{T_M_3_4}
T^{M}_{x^-}: \ & \nonumber \\
&
\arraycolsep=2.pt\def\arraystretch{0.625}
\begin{array}{rrrrrrrrrrrrrrrr}
f(x^+_0) & \ \leftarrow  
			 && \tfrac{ 35}{ 576} f(x^-_{-1})
                          &+&\tfrac{389}{ 384} f(x^-_{ 0})
                          &-&\tfrac{  1}{  12} f(x^-_{ 1})
                          &+&\tfrac{ 11}{1152} f(x^-_{ 2}) 
                          &&
\\
f(x^+_1) & \ \leftarrow 
			 &-&\tfrac{  77}{1728} f(x^-_{-1})
                          &+&\tfrac{1325}{3456} f(x^-_{ 0})
                          &+&\tfrac{   3}{   4} f(x^-_{ 1})
                          &-&\tfrac{ 335}{3456} f(x^-_{ 2})
                          &+&\tfrac{   7}{ 864} f(x^-_{ 3}) 
\\
f(x^+_2) & \ \leftarrow 
			 && \tfrac{   7}{ 864} f(x^-_{-1})
                          &-&\tfrac{ 335}{3456} f(x^-_{ 0})
                          &+&\tfrac{   3}{   4} f(x^-_{ 1})
                          &+&\tfrac{1325}{3456} f(x^-_{ 2})
                          &-&\tfrac{  77}{1728} f(x^-_{ 3})
\\
f(x^+_3) & \ \leftarrow 
			&&
			&& \tfrac{ 11}{1152} f(x^-_{ 0})
                         &-&\tfrac{  1}{  12} f(x^-_{ 1})
                         &+&\tfrac{389}{ 384} f(x^-_{ 2})
                         &+&\tfrac{ 35}{ 576} f(x^-_{ 3})
\end{array} 
&
\end{flalign}
}

\indent We note here that the above interpolation formulas are designed such that the leading error terms in the approximations are of at least third order.
However, even under this condition, they are not the unique choice for compatible interpolation operators. 
Taking the 2:3 grid spacing ratio for example, the interpolation operators shown in \eqref{T_N_2_3_alternative} and \eqref{T_M_2_3_alternative} are also suitable choices.

\begin{figure}[H]
\captionsetup{width=0.9\textwidth, font=small,labelfont=small}
\centering
%
\centering\includegraphics[scale=0.125]{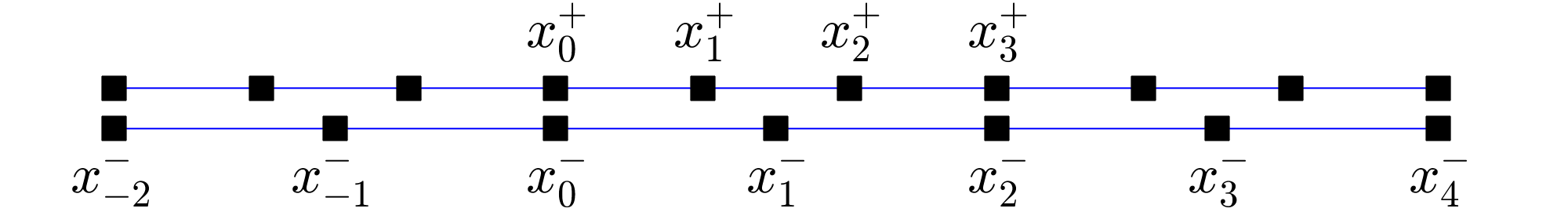}
%
\\[0.5ex]
%
\centering\includegraphics[scale=0.125]{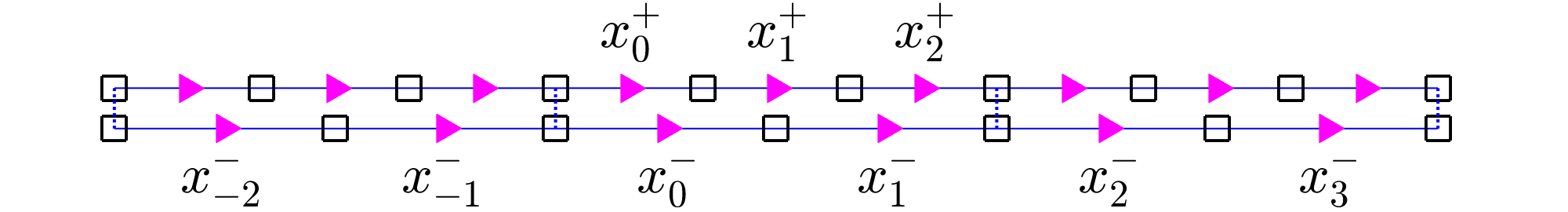}
%
%
\caption{Illustration of the grid points involved in formulas \eqref{T_N_2_3_alternative} and \eqref{T_M_2_3_alternative} for interpolation operators $T^{N}_{x^-}$ (top) and $T^{M}_{x^-}$ (bottom). The grid spacing ratio is 2:3. } 
\label{Grid_configuration_transfer_operator_2_3_alternative}
\end{figure}
{\footnotesize
\begin{flalign}
\label{T_N_2_3_alternative}
T^{N}_{x^-}: \ & \nonumber \\
&
\arraycolsep=2.pt\def\arraystretch{0.625}
\begin{array}{rrrrrrrrrrrrrrrrr}
f(x^+_0) & \ \leftarrow 
			    &-&\tfrac{ 1}{96} f(x^-_{-2})
          &+&\tfrac{ 1}{24} f(x^-_{-1})
          &+&\tfrac{15}{16} f(x^-_{ 0})
          &+&\tfrac{ 1}{24} f(x^-_{ 1})
          &-&\tfrac{ 1}{96} f(x^-_{ 2})
          & &
          & &
\\
f(x^+_1) & \ \leftarrow 
          & &
			    &-&\tfrac{ 13}{288} f(x^-_{-1})
          &+&\tfrac{103}{288} f(x^-_{ 0})
          &+&\tfrac{217}{288} f(x^-_{ 1})
          &-&\tfrac{ 19}{288} f(x^-_{ 2})
          & &
          & &
\\
f(x^+_2) & \ \leftarrow 
			    & & 
          & & 
          &-&\tfrac{ 19}{288} f(x^-_{ 0})
          &+&\tfrac{217}{288} f(x^-_{ 1})
          &+&\tfrac{103}{288} f(x^-_{ 2})
          &-&\tfrac{ 13}{288} f(x^-_{ 3})
          & & 
\\
f(x^+_3) & \ \leftarrow 
          & &
          & &
          &-&\tfrac{ 1}{96} f(x^-_{ 0})
          &+&\tfrac{ 1}{24} f(x^-_{ 1})
          &+&\tfrac{15}{16} f(x^-_{ 2})
          &+&\tfrac{ 1}{24} f(x^-_{ 3})
          &-&\tfrac{ 1}{96} f(x^-_{ 4})
\end{array}
&
\end{flalign}
}
{\footnotesize
\begin{flalign}
\label{T_M_2_3_alternative}
T^{M}_{x^-}: \ & \nonumber \\
&
\arraycolsep=2.pt\def\arraystretch{0.625}
\begin{array}{rrrrrrrrrrrrrrr}
f(x^+_0) & \ \leftarrow 
			    &-&\tfrac{ 11}{576} f(x^-_{-2}) 
          &+&\tfrac{ 89}{576} f(x^-_{-1}) 
          &+&\tfrac{527}{576} f(x^-_{ 0}) 
          &-&\tfrac{ 29}{576} f(x^-_{ 1}) 
          & &
          & &
\\
f(x^+_1) & \ \leftarrow 
          & &
			    &-&\tfrac{1}{16} f(x^-_{-1})
          &+&\tfrac{9}{16} f(x^-_{ 0})
          &+&\tfrac{9}{16} f(x^-_{ 1})
          &-&\tfrac{1}{16} f(x^-_{ 2})
          & &
\\
f(x^+_2) & \ \leftarrow 
          & &
          & &
			    &-&\tfrac{ 29}{576} f(x^-_{ 0})
          &+&\tfrac{527}{576} f(x^-_{ 1})
          &+&\tfrac{ 89}{576} f(x^-_{ 2})
          &-&\tfrac{ 11}{576} f(x^-_{ 3})
\end{array} 
&
\end{flalign}
}

\section{1D SBP finite difference operators}\label{appendix_SBP_operators}
The following 1D SBP finite difference operators in equations \eqref{SBP_matrices_1D_D} and \eqref{SBP_matrices_1D_A} are used as the building blocks for the 3D SBP operators presented in section \eqref{SBP_operators_3D}. 
We note here that these operators correspond to the case of unit grid spacing, i.e., $\Delta x=1$. 
For general cases, $\mathcal D^N$ and $\mathcal D^M$ need to be scaled by $\nicefrac{1}{\Delta x}$ while $\mathcal A^N$ and $\mathcal A^M$ need to be scaled by $\Delta x$.

\begin{subequations}
\label{SBP_matrices_1D_D}
\footnotesize
\begin{align}
\label{SBP_matrices_1D_DN}
\mathcal D^N & = \left[ 
\arraycolsep=3.2pt\def\arraystretch{0.625}
\begin{array}{r r r r r r r r r r r r r r r r r} 
-\nicefrac{79}{78}   &   \nicefrac{27}{26}   &   -\nicefrac{1}{26}   &   \nicefrac{1}{78}   &   \multicolumn{1}{r|}{0}      \\
 \nicefrac{2}{21}     &  -\nicefrac{9}{7}       &    \nicefrac{9}{7}     &  -\nicefrac{2}{21}   &   \multicolumn{1}{r|}{0}      \\
 \nicefrac{1}{75}     &    0            &   -\nicefrac{27}{25}  &   \nicefrac{83}{75} &  \multicolumn{1}{r|}{ -\nicefrac{1}{25} } \\ \cline{1-5}
   &     &  \nicefrac{1}{24}  &   -\nicefrac{9}{8}   &   \nicefrac{9}{8}   &  -\nicefrac{1}{24} &  \\
   &     &    &  \nicefrac{1}{24}  &  -\nicefrac{9}{8}   &  \nicefrac{9}{8}   &  -\nicefrac{1}{24} \\
& & & & \ddots & \ddots & \ddots & \ddots & \\  
& & & & & \nicefrac{1}{24}  &  -\nicefrac{9}{8}   &  \nicefrac{9}{8}   &  -\nicefrac{1}{24} \\
& & & & & & \nicefrac{1}{24}  &  -\nicefrac{9}{8}   &  \nicefrac{9}{8}   &  -\nicefrac{1}{24} \\ \cline{8-12}
& & & & & & & \multicolumn{1}{|r}{ \nicefrac{1}{25} }  &   -\nicefrac{83}{75}   &   \nicefrac{27}{25}   &    0  &   -\nicefrac{1}{75}   \\
& & & & & & & \multicolumn{1}{|r}{0}   &    \nicefrac{2}{21}     &  -\nicefrac{9}{7}       &    \nicefrac{9}{7}     &   -\nicefrac{2}{21}   \\  
& & & & & & & \multicolumn{1}{|r}{0}   &   -\nicefrac{1}{78}     &   \nicefrac{1}{26}     &   -\nicefrac{27}{26} &    \nicefrac{79}{78} 
\end{array} 
\right];
\allowdisplaybreaks[4]
\\
\label{SBP_matrices_1D_DM}
\mathcal D^M & = \left[ 
\arraycolsep=3.2pt\def\arraystretch{0.625}
\begin{array}{r r r r r r r r r r r r r r r r r} 
-2       &     3    &    -1    &    0    &   \multicolumn{1}{r|}{0}      \\
-1       &     1    &      0   &    0    &   \multicolumn{1}{r|}{0}     \\
\nicefrac{1}{24}  &  -\nicefrac{9}{8}    &     \nicefrac{9}{8} &   -\nicefrac{1}{24} &  \multicolumn{1}{r|}{0}  &  \\
-\nicefrac{1}{71} &   \nicefrac{6}{71}  &  -\nicefrac{83}{71}  &  \nicefrac{81}{71} &  \multicolumn{1}{r|}{ -\nicefrac{3}{71} } &  \\ \cline{1-5}
   &     &  \nicefrac{1}{24}  &   -\nicefrac{9}{8}   &   \nicefrac{9}{8}   &  -\nicefrac{1}{24} &  \\
   &     &    &  \nicefrac{1}{24}  &  -\nicefrac{9}{8}   &  \nicefrac{9}{8}   &  -\nicefrac{1}{24} \\
& & & & \ddots & \ddots & \ddots & \ddots & \\  
& & & & & \nicefrac{1}{24}  &  -\nicefrac{9}{8}   &  \nicefrac{9}{8}   &  -\nicefrac{1}{24} \\
& & & & & & \nicefrac{1}{24}  &  -\nicefrac{9}{8}   &  \nicefrac{9}{8}   &  -\nicefrac{1}{24} \\ \cline{8-12}
& & & & & & & \multicolumn{1}{|r}{ \nicefrac{3}{71} }  &  -\nicefrac{81}{71} & \nicefrac{83}{71} & -\nicefrac{6}{71}  & \nicefrac{1}{71} \\
& & & & & & & \multicolumn{1}{|r}{0} & \nicefrac{1}{24}  &  -\nicefrac{9}{8}  & \nicefrac{9}{8} & -\nicefrac{1}{24}  \\
& & & & & & & \multicolumn{1}{|r}{0} & 0            &  0               &            -1   &   1       \\
& & & & & & & \multicolumn{1}{|r}{0} & 0            &  1           &                  -3   &    2   
\end{array} 
\right].
\end{align}
\end{subequations}

\begin{subequations}
\label{SBP_matrices_1D_A}
\footnotesize
\begin{align}
\label{SBP_matrices_1D_AN}
\mathcal A^N & = \left[ 
\arraycolsep=5pt\def\arraystretch{0.625}
\begin{array}{r r r r r r r r r r r r r r r r r} 
\nicefrac{7}{18}   &    &    &    \multicolumn{1}{r|}{}      \\
               &  \nicefrac{9}{8}    &     &   \multicolumn{1}{r|}{}     \\
        & &  1   &    \multicolumn{1}{r|}{}     \\               
        &    &    & \multicolumn{1}{r|}{ \nicefrac{71}{72} }  \\ \cline{1-4}
&    &    &    &  1  &    \\
&    &    &    &    &  1  &    \\
&    &    &    &    &      &  \ddots   \\
&    &    &    &    &    &    & 1 &    \\
&    &    &    &    &    &    &    & 1 &    \\ \cline{10-13}
&    &    &    &    &    &    &    &    &  \multicolumn{1}{|r}{ \nicefrac{71}{72} } \\
&    &    &    &    &    &    &    &    &  \multicolumn{1}{|r}{}   &  1  \\
&    &    &    &    &    &    &    &    &  \multicolumn{1}{|r}{}   &    &  \nicefrac{9}{8}  \\
&    &    &    &    &    &    &    &    &  \multicolumn{1}{|r}{}   &    &    &  \nicefrac{7}{18} 
\end{array} 
\right];
\allowdisplaybreaks[4]
\\
\label{SBP_matrices_1D_AM}
\mathcal A^M & = \left[ 
\arraycolsep=5pt\def\arraystretch{0.625}
\begin{array}{r r r r r r r r r r r r r r r r r} 
\nicefrac{13}{12}   &    &     \multicolumn{1}{r|}{}      \\
               &  \nicefrac{7}{8}    &    \multicolumn{1}{r|}{}     \\
        &    &    \multicolumn{1}{r|}{ \nicefrac{25}{24} }  \\ \cline{1-3}
&    &    & 1 &    \\
&    &    &    & 1 &    \\
&    &    &    &    &  \ddots   \\
&    &    &    &    &    & 1 &    \\
&    &    &    &    &    &    & 1 &    \\ \cline{9-11}
&    &    &    &    &    &    &    &  \multicolumn{1}{|r}{ \nicefrac{25}{24} } \\
&    &    &    &    &    &    &    &  \multicolumn{1}{|r}{}    &  \nicefrac{7}{8}  \\
&    &    &    &    &    &    &    &  \multicolumn{1}{|r}{}    &    &  \nicefrac{13}{12} 
\end{array} 
\right].
\end{align}
\end{subequations}

With the above operators, matrix $Q$, defined as $\mathcal A^N \mathcal D^M \ + \ ( \mathcal A^M \mathcal D^N )^T$ in equation \eqref{1D_Q_definition} takes the following form:
\begin{equation}
\label{Q_matrix_1D}
\footnotesize
\left[ 
\def\arraystretch{0.625}
\begin{array}{r r r r r r r r r r r r r r r r} 
-\nicefrac{15}{8}  & \nicefrac{5}{4}  &  -\nicefrac{3}{8}      \\
&    &    &    &    &    \\
&    &    &    &    &    \\
&    &    &    &    &    &    &    &      \\
&    &    &    &    &    &    &    &      \\
&    &    &    &    &    &    &    &     \nicefrac{3}{8}  &  -\nicefrac{5}{4}  &  \nicefrac{15}{8}
\end{array} 
\right],
\end{equation}
which can be expressed as $\mathcal E^L \left(\mathcal P^L\right)^T + \, \mathcal E^R \left(\mathcal P^R\right)^T$
with $\mathcal E^L$, $\mathcal E^R$, $\mathcal P^L$, and $\mathcal P^R$ given by 
\begin{equation}
\label{Vectors_E_and_P}
\footnotesize
\mathcal E^L = 
\left[ \def\arraystretch{0.625}\begin{array}{c} 1 \\ 0 \\ \vdots \\ 0 \end{array} \right], \quad
\mathcal E^R = 
\left[ \def\arraystretch{0.625}\begin{array}{c} 0 \\ \vdots \\ 0 \\ 1 \end{array} \right], \quad
\mathcal P^L = 
\left[ \def\arraystretch{0.625}\begin{array}{c} \nicefrac{15}{8} \\ - \nicefrac{5}{4} \\ \nicefrac{3}{8} \\ 0 \\ \vdots \\ 0 \end{array} \right], \quad
\mathcal P^R = 
\left[ \def\arraystretch{0.625}\begin{array}{c} 0 \\ \vdots \\ 0 \\ \nicefrac{3}{8} \\ - \nicefrac{5}{4} \\ \nicefrac{15}{8} \end{array} \right],
\end{equation}
respectively.

\addtolength{\bibsep}{-0.5em}
\renewcommand{\bibfont}{\normalfont\small}
\bibliographystyle{./bst_base/abbrv.bst}
\bibliography{refs}

\begin{thebibliography}{10}

\bibitem{ogilvie1996effects}
J.~S. Ogilvie and G.~W. Purnell.
\newblock Effects of salt-related mode conversions on subsalt prospecting.
\newblock {\em Geophysics}, 61(2):331--348, 1996.

\bibitem{miller1955partition}
G.~Miller and H.~Pursey.
\newblock On the partition of energy between elastic waves in a semi-infinite
  solid.
\newblock {\em Proceedings of the Royal Society of London Series A},
  233:55--69, 1955.

\bibitem{thierry1987acoustics}
T.~Bourbi{\'e}, O.~Coussy, and B.~Zinszner.
\newblock {\em Acoustics of Porous Media}.
\newblock Institut Fran{\c{c}}ais du P{\'e}trole Publications. Editions
  Technip, 1987.

\bibitem{lysmer1972finite}
J.~Lysmer and L.~A. Drake.
\newblock {\em A Finite Element Method for Seismology. {In}: {Methods in
  Computational Physics: Advances in Research and Applications}}.
\newblock Academic Press, 1972.

\bibitem{BaoBielakGhattasEtAl98}
H.~Bao, J.~Bielak, O.~Ghattas, L.~F. Kallivokas, D.~R. O'Hallaron, J.~R.
  Shewchuk, and J.~Xu.
\newblock Large-scale simulation of elastic wave propagation in heterogeneous
  media on parallel computers.
\newblock {\em Computer Methods in Applied Mechanics and Engineering},
  152(1--2):85--102, January 1998.

\bibitem{komatitsch1998spectral}
D.~Komatitsch and J.-P. Vilotte.
\newblock The spectral element method: {An} efficient tool to simulate the
  seismic response of {2D} and {3D} geological structures.
\newblock {\em Bulletin of the Seismological Society of America},
  88(2):368--392, 1998.

\bibitem{komatitsch1999introduction}
D.~Komatitsch and J.~Tromp.
\newblock Introduction to the spectral element method for three-dimensional
  seismic wave propagation.
\newblock {\em Geophysical Journal International}, 139(3):806--822, 1999.

\bibitem{kaser2006arbitrary}
M.~K{\"a}ser and M.~Dumbser.
\newblock An arbitrary high-order discontinuous {Galerkin} method for elastic
  waves on unstructured meshes -- {I.} {The} two-dimensional isotropic case
  with external source terms.
\newblock {\em Geophysical Journal International}, 166(2):855--877, 2006.

\bibitem{dumbser2006arbitrary}
M.~Dumbser and M.~K{\"a}ser.
\newblock An arbitrary high-order discontinuous {Galerkin} method for elastic
  waves on unstructured meshes -- {II.} {The} three-dimensional isotropic case.
\newblock {\em Geophysical Journal International}, 167(1):319--336, 2006.

\bibitem{chaljub2007spectral}
E.~Chaljub, D.~Komatitsch, J.-P. Vilotte, Y.~Capdeville, B.~Valette, and
  G.~Festa.
\newblock Spectral-element analysis in seismology.
\newblock {\em Advances in Geophysics}, 48:365--419, 2007.

\bibitem{WilcoxStadlerBursteddeEtAl10}
L.~C. Wilcox, G.~Stadler, C.~Burstedde, and O.~Ghattas.
\newblock A high-order discontinuous {G}alerkin method for wave propagation
  through coupled elastic-acoustic media.
\newblock {\em Journal of Computational Physics}, 229(24):9373--9396, 2010.

\bibitem{bui2012analysis}
T.~Bui-Thanh and O.~Ghattas.
\newblock Analysis of an hp-nonconforming discontinuous {Galerkin} spectral
  element method for wave propagation.
\newblock {\em SIAM Journal on Numerical Analysis}, 50(3):1801--1826, 2012.

\bibitem{jastram1994elastic}
C.~Jastram and E.~Tessmer.
\newblock Elastic modelling on a grid with vertically varying spacing.
\newblock {\em Geophysical prospecting}, 42(4):357--370, 1994.

\bibitem{aoi19993d}
S.~Aoi and H.~Fujiwara.
\newblock {3D} finite-difference method using discontinuous grids.
\newblock {\em Bulletin of the Seismological Society of America},
  89(4):918--930, 1999.

\bibitem{hayashi2001discontinuous}
K.~Hayashi, D.~R. Burns, and M.~N. Toks{\"o}z.
\newblock Discontinuous-grid finite-difference seismic modeling including
  surface topography.
\newblock {\em Bulletin of the Seismological Society of America},
  91(6):1750--1764, 2001.

\bibitem{kristek2010stable}
J.~Kristek, P.~Moczo, and M.~Galis.
\newblock Stable discontinuous staggered grid in the finite-difference
  modelling of seismic motion.
\newblock {\em Geophysical Journal International}, 183(3):1401--1407, 2010.

\bibitem{zhang2013stable}
Z.~Zhang, W.~Zhang, H.~Li, and X.~Chen.
\newblock Stable discontinuous grid implementation for collocated-grid
  finite-difference seismic wave modelling.
\newblock {\em Geophysical Journal International}, 192(3):1179--1188, 2013.

\bibitem{nie2017fourth}
S.~Nie, Y.~Wang, K.~B. Olsen, and S.~M. Day.
\newblock Fourth-order staggered-grid finite-difference seismic wavefield
  estimation using a discontinuous mesh interface {(WEDMI)}.
\newblock {\em Bulletin of the Seismological Society of America},
  107(5):2183--2193, 2017.

\bibitem{gao2018long}
L.~Gao, D.~Ketcheson, and D.~Keyes.
\newblock On long-time instabilities in staggered finite difference simulations
  of the seismic acoustic wave equations on discontinuous grids.
\newblock {\em Geophysical Journal International}, 212(2):1098--1110, 2018.

\bibitem{fernandez2014review}
D.~Fern{\'a}ndez, J.~Hicken, and D.~Zingg.
\newblock Review of summation-by-parts operators with simultaneous
  approximation terms for the numerical solution of partial differential
  equations.
\newblock {\em Computers \& Fluids}, 95:171--196, 2014.

\bibitem{svard2014review}
M.~Sv{\"a}rd and J.~Nordstr{\"o}m.
\newblock Review of summation-by-parts schemes for initial--boundary-value
  problems.
\newblock {\em Journal of Computational Physics}, 268:17--38, 2014.

\bibitem{appelo2009stable}
D.~Appel{\"o} and N.~A. Petersson.
\newblock A stable finite difference method for the elastic wave equation on
  complex geometries with free surfaces.
\newblock {\em Communications in Computational Physics}, 5(1):84--107, 2009.

\bibitem{sjogreen2012fourth}
B.~Sj{\"o}green and N.~A. Petersson.
\newblock A fourth order accurate finite difference scheme for the elastic wave
  equation in second order formulation.
\newblock {\em Journal of Scientific Computing}, 52(1):17--48, 2012.

\bibitem{petersson2015wave}
N.~A. Petersson and B.~Sj{\"o}green.
\newblock Wave propagation in anisotropic elastic materials and curvilinear
  coordinates using a summation-by-parts finite difference method.
\newblock {\em Journal of Computational Physics}, 299:820--841, 2015.

\bibitem{wang2016high}
S.~Wang, K.~Virta, and G.~Kreiss.
\newblock High order finite difference methods for the wave equation with
  non-conforming grid interfaces.
\newblock {\em Journal of Scientific Computing}, 68(3):1002--1028, 2016.

\bibitem{o2017energy}
O.~O'Reilly, T.~Lundquist, E.~M. Dunham, and J.~Nordstr{\"o}m.
\newblock Energy stable and high-order-accurate finite difference methods on
  staggered grids.
\newblock {\em Journal of Computational Physics}, 346:572--589, 2017.

\bibitem{gao2019combining}
L.~Gao and D.~Keyes.
\newblock Combining finite element and finite difference methods for isotropic
  elastic wave simulations in an energy-conserving manner.
\newblock {\em Journal of Computational Physics}, 378:665--685, 2019.

\bibitem{gao2019sbp}
L.~Gao, D.~Fern{\'a}ndez, M.~Carpenter, and D.~Keyes.
\newblock {SBP}-{SAT} finite difference discretization of acoustic wave
  equations on staggered block-wise uniform grids.
\newblock {\em Journal of Computational and Applied Mathematics}, 348:421--444,
  2019.

\bibitem{gao2020explicit}
L.~Gao and D.~Keyes.
\newblock Explicit coupling of acoustic and elastic wave propagation in
  finite-difference simulations.
\newblock {\em Geophysics}, 85(5):T293--T308, 2020.

\bibitem{gao2020simultaneous}
L.~Gao and D.~Keyes.
\newblock Simultaneous approximation terms for elastic wave equations on
  nonuniform grids.
\newblock In R.~Haynes, S.~MacLachlan, X.-C. Cai, L.~Halpern, H.~H. Kim,
  A.~Klawonn, and O.~Widlund, editors, {\em Domain Decomposition Methods in
  Science and Engineering XXV}, pages 125--133. Springer International
  Publishing, 2020.

\bibitem{yee1966numerical}
K.~Yee.
\newblock Numerical solution of initial boundary value problems involving
  {Maxwell's} equations in isotropic media.
\newblock {\em IEEE Transactions on Antennas and Propagation}, 14(3):302--307,
  1966.

\bibitem{virieux1986p}
J.~Virieux.
\newblock {P-SV} wave propagation in heterogeneous media: {Velocity}-stress
  finite-difference method.
\newblock {\em Geophysics}, 51(4):889--901, 1986.

\bibitem{levander1988fourth}
A.~Levander.
\newblock Fourth-order finite-difference {P-SV} seismograms.
\newblock {\em Geophysics}, 53(11):1425--1436, 1988.

\bibitem{kreiss1974finite}
H.-O. Kreiss and G.~Scherer.
\newblock {\em Finite element and finite difference methods for hyperbolic
  partial differential equations. {In}: {Mathematical} aspects of finite
  elements in partial differential equations}.
\newblock Academic Press, 1974.

\bibitem{carpenter1994time}
M.~Carpenter, D.~Gottlieb, and S.~Abarbanel.
\newblock Time-stable boundary conditions for finite-difference schemes solving
  hyperbolic systems: methodology and application to high-order compact
  schemes.
\newblock {\em Journal of Computational Physics}, 111(2):220--236, 1994.

\bibitem{hesthaven2008nodal}
J.~Hesthaven and T.~Warburton.
\newblock {\em Nodal discontinuous {Galerkin} methods: algorithms, analysis,
  and applications}.
\newblock Springer Science \& Business Media, 2008.

\bibitem{hicken2013summation}
J.~Hicken and D.~Zingg.
\newblock Summation-by-parts operators and high-order quadrature.
\newblock {\em Journal of Computational and Applied Mathematics},
  237(1):111--125, 2013.

\bibitem{strand1994summation}
B.~Strand.
\newblock Summation by parts for finite difference approximations for d/dx.
\newblock {\em Journal of Computational Physics}, 110(1):47--67, 1994.

\bibitem{stein2009introduction}
S.~Stein and M.~Wysession.
\newblock {\em An introduction to seismology, earthquakes, and earth
  structure}.
\newblock John Wiley \& Sons, 2009.

\bibitem{aminzadeh19973}
F.~Aminzadeh, B.~Jean, and T.~Kunz.
\newblock {\em {3-D} Salt and Overthrust Models}.
\newblock Society of Exploration Geophysicists, 1997.

\bibitem{gregory1976fluid}
A.~Gregory.
\newblock Fluid saturation effects on dynamic elastic properties of sedimentary
  rocks.
\newblock {\em Geophysics}, 41(5):895--921, 1976.

\bibitem{hamilton1979v}
E.~L. Hamilton.
\newblock $v_p$/$v_s$ and {Poisson’s} ratios in marine sediments and rocks.
\newblock {\em The Journal of the Acoustical Society of America},
  66(4):1093--1101, 1979.

\end{thebibliography}

\newpage
\supp

\begin{center}
{\Large Supplementary Materials}
\end{center}

\section{Additional plots related to the example of Layered model}\label{supp_additional_plots}

\begin{figure}[H]
\captionsetup{font=small,labelfont=small}
\centering
\begin{subfigure}[b]{0.495\textwidth}
%
\captionsetup{width=1\textwidth, font=small,labelfont=small}
\centering\includegraphics[scale=0.125]{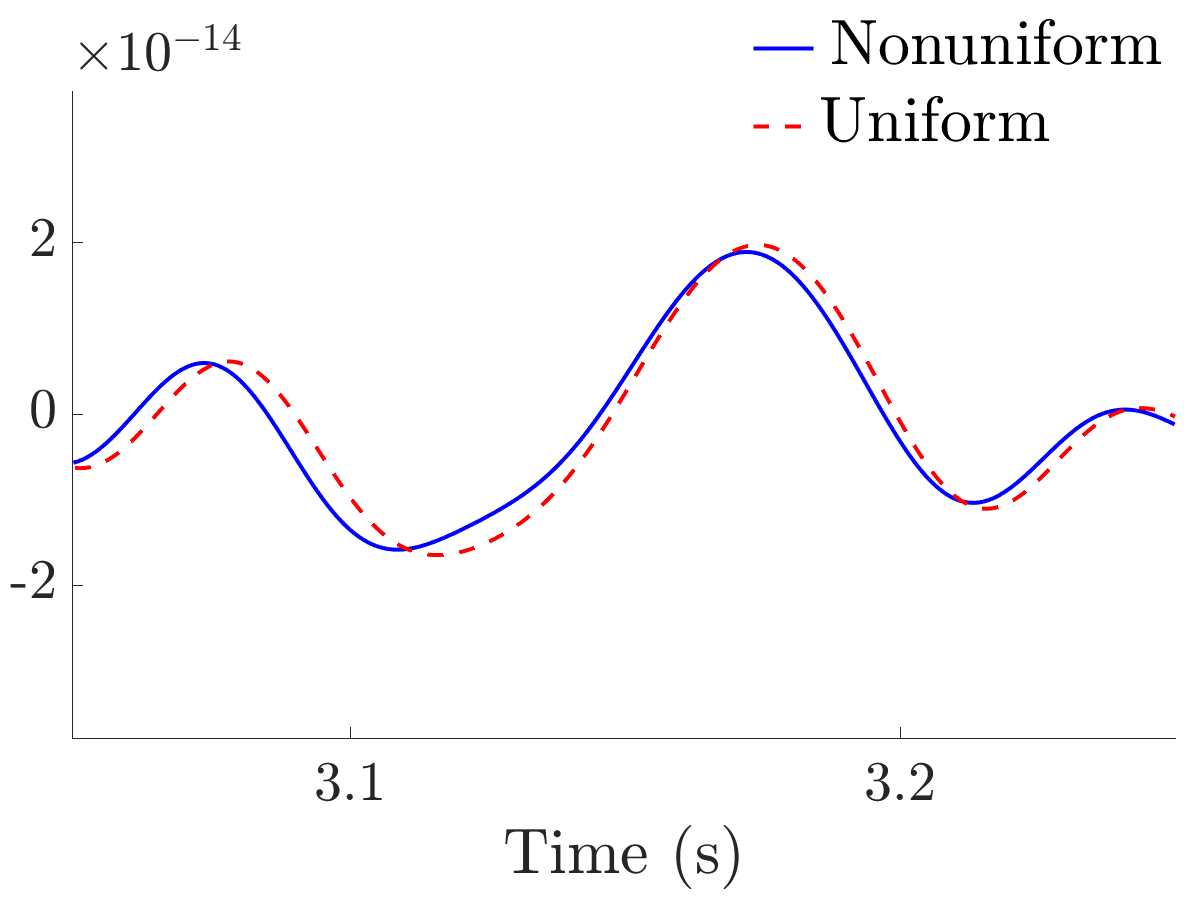}
%
\end{subfigure}
\begin{subfigure}[b]{0.495\textwidth}
\captionsetup{width=1\textwidth, font=small,labelfont=small}
\centering\includegraphics[scale=0.125]{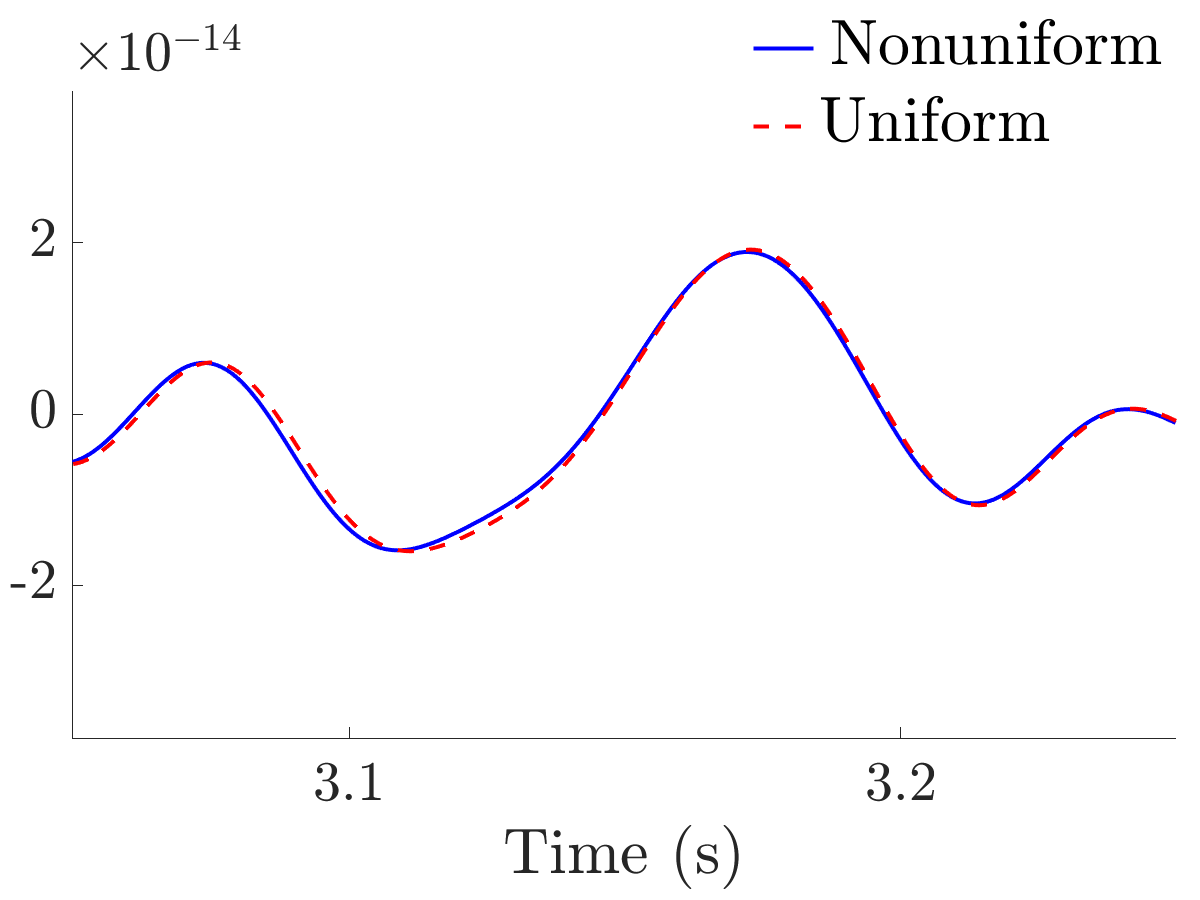}
%
\end{subfigure} 
\caption[]{
Segments of the seismograms between 3.05 s and 3.25 s. The left figure corresponds to Figures~\ref{fig_layered_seismogram} of the manuscript; the right figure corresponds to Figure~\ref{fig_layered_seismogram_3x} of the manuscript, which uses three times finer grid spacing. The improved agreement can be observed more easily here than by comparing Figures~\ref{fig_layered_seismogram} and~\ref{fig_layered_seismogram_3x} of the manuscript.}
\label{fig_seismograms_small_segment}
\end{figure}

\section{Additional expenditure information on computing resources}\label{supp_additional_tables}

\begin{table}[H]
\captionsetup{font=small,labelfont=small}
\caption[]{Expenditure on computing resources corresponding to Figure~\ref{fig_layered_seismogram_3x} of the manuscript. The expected ratio in total computational resources, i.e., core hours, is $\sim$$19.1$, based on the amount of spatial-temporal discretization points involved. The measured ratio is $\sim$$18.5$, which suggests a $\sim$$3\%$ overhead on a per core basis in our implementation.}
\label{computing_resources_figure_5}
\arraycolsep=2.pt\def\arraystretch{1}
\centering
\small
\begin{tabular}{| r | c | c | c |}
\hline
	     	  & core count & elapsed time & total core hours \\ \hline
uniform 	  & 23328   &  63482 seconds  & $4.114 \times 10^5$  \\ \hline
nonuniform & 6785    &  11805 seconds &  $2.225 \times 10^4$  \\ \hline
\end{tabular}
\end{table}

\begin{table}[H]
\captionsetup{font=small,labelfont=small}
\caption[]{Expenditure on computing resources corresponding to Figure~\ref{fig_overthrust_seismogram} of the manuscript. The expected ratio in total computational resources, i.e., core hours, is $\sim$$3.153$, based on the amount of spatial-temporal discretization points involved. The measured ratio is $\sim$$3.095$, which suggests a $\sim$$2\%$ overhead on a per core basis in our implementation.}
\label{computing_resources_figure_9}
\arraycolsep=2.pt\def\arraystretch{1}
\centering
\small
\begin{tabular}{| r | c | c | c |}
\hline
	     	  & core count & elapsed time & total core hours \\ \hline
uniform 	  & 19200   &  12787 seconds  & $6.820 \times 10^4$  \\ \hline
nonuniform & 7300    &  10860 seconds &  $2.202 \times 10^4$  \\ \hline
\end{tabular}
\end{table}

\section{Supplementary performance tests}\label{supp_performance_tests}

In the comparison of the elapsed wall clock times shown in section Numerical examples, each processor core is assigned the same number of grid points for both nonuniform and uniform simulations. 
This ensures that, on each processor core, the pressure on performance-critical resources are similar and that the comparison is fair on a per core basis.
However, this implies that the total number of processor cores used in the uniform simulation is larger, which may lead to an unfair comparison if a larger number of processor cores incurs considerable communication overhead.
(Such a situation would unfairly favor the proposed nonuniform simulation.) 
%
%
In the following, we provide some additional test results to further validate our algorithmic intuition that the nonuniform simulations incur little overhead. The setup of the physical problem is the same as in Example Layered model. 

First, we run the uniform simulation using a number of processor cores (6480 = $12 \times 15 \times 36$) similar to that used in the nonuniform case (6785 cores) presented in Example Layered model. 
Each processor core is assigned $90 \times 72 \times 30$ grid points.
The elapsed time (averaged over 5 runs) is $\sim$3363 s. 
In terms of the total resources consumed (i.e., number of processor cores multiplied by elapsed time), this simulation setting is actually $\sim$7\% more costly than the uniform simulation setting used in Example Layered model. 
(Plausible explanations for the increase of consumed resources include heavier memory pressure on each processor core, longer runtime and therefore greater susceptibility to interference from the shared hardware environment, and finally, natural performance variances on modern large-scale computing systems.
Nevertheless, this test suggests that the performance comparison presented in Example Layered model did not unfairly favor the proposed nonuniform simulation due to the difference in core numbers.)
%

Additionally, some weak scaling results are provided in Tables~\ref{weak_scaling_nonuniform} and~\ref{weak_scaling_uniform}, where the number of cores in the $z$ direction doubles from column to column.
The grid points per processor core are kept the same at $60 \times 60 \times 15$. 
All simulations are performed for 2000 time steps.
The elapsed times shown in Tables~\ref{weak_scaling_nonuniform} and~\ref{weak_scaling_uniform} are averaged over 5 runs.
From both tables, we observe that increasing the number of cores does not lead to an increase in simulation time for both the nonuniform and uniform simulations, suggesting that both weak-scale well.

\begin{table}[H]
\caption[]{Weak scaling (nonuniform).}
\label{weak_scaling_nonuniform}
\setlength{\tabcolsep}{12pt}
\def\arraystretch{1}
\centering
\small
\begin{tabular}{| c | c | c | c | c | c |}
\hline
Number of cores & 6785 & 13570 & 27140 & 54280 & 108560 \\ \hline
Elapsed time (seconds) & 45.7 & 45.4 & 44.2 & 43.5 & 43.2 \\ \hline
\end{tabular}
\end{table}
\begin{table}[H]
\caption[]{Weak scaling (uniform).}
\label{weak_scaling_uniform}
\setlength{\tabcolsep}{12pt}
\def\arraystretch{1}
\centering
\small
\begin{tabular}{| c | c | c | c | c |}
\hline
Number of cores & 11664 & 23328 & 46656 & 93312 \\ \hline
Elapsed time (seconds) & 46.4 & 44.6 & 43.5 & 43.3 \\ \hline
\end{tabular}
\end{table}

Together, these tests suggest that for the comparison made in Example Layered model, difference in the sizes of the processor network 
is unlikely to be a major contributing factor.

\end{document}